\documentclass[aop,MSNbibl,citesort,dvips]{arximspdf}
\usepackage{mathbh}
\usepackage{subenv}
\usepackage{mathrsfs}
\usepackage{stmaryrd}
\usepackage{stfloats}
\usepackage{graphicx}

\doi{10.1214/11-AOP675}
\volume{40}
\issue{5}
\pubyear{2012}
\firstpage{1897}
\lastpage{1944}

\makeatletter

\fnbelowfloat

\newcommand{\bs}{\setminus}
\newcommand{\eps}{\varepsilon}
\newcommand{\lbracket}{[\![}
\newcommand{\rbracket}{]\!]}
\newcommand{\de}{:=}

\newcommand{\xrightarrow}{\mathop{\longrightarrow}}
\newcommand{\ton}{\xrightarrow_{n\to\infty}}
\newcommand{\tolk}{\xrightarrow_{k \to\infty}^{(d)}}

\newcommand{\MM}{\mathbb{M}}
\newcommand{\N}{\mathbb{N}}
\newcommand{\Pb}{\mathbb{P}}
\newcommand{\R}{\mathbb{R}}
\newcommand{\TTT}{\mathbb{T}}
\newcommand{\Z}{\mathbb{Z}}

\newcommand{\C}{\mathcal{C}}
\newcommand{\EE}{\mathcal{E}}
\newcommand{\F}{\mathcal{F}}
\newcommand{\I}{\mathcal{I}}
\newcommand{\K}{\mathcal{K}}
\newcommand{\Q}{\mathcal{Q}}
\newcommand{\Qb}{\mathcal{Q}^{\bullet}}
\newcommand{\cS}{\mathcal{S}}
\newcommand{\T}{\mathcal{T}}
\newcommand{\U}{\mathcal{U}}
\newcommand{\V}{\mathcal{V}}
\newcommand{\W}{\mathcal{W}}
\newcommand{\X}{\mathcal{X}}

\newcommand{\Leb}{\mathscr{L}}
\newcommand{\loo}{\mathscr{L}}
\newcommand{\TT}{\mathscr{T}}

\newcommand{\aaa}{{\mathfrak{a}}}
\newcommand{\CC}{{\mathfrak{C}}}
\newcommand{\fF}{{\mathfrak{F}}}
\newcommand{\ii}{{\mathfrak{i}}}
\newcommand{\Lab}{{\mathfrak{L}}}
\newcommand{\m}{{\mathfrak{m}}}
\newcommand{\qi}{\mathfrak{q}_{\infty}}
\newcommand{\s}{{\mathfrak{s}}}
\newcommand{\Sg}{{\mathfrak{S}}}
\newcommand{\tr}{{\mathfrak{t}}}

\newcommand{\oo}{\mathbf{0}}
\newcommand{\eee}{\mathbf{e}}
\newcommand{\sss}{\mathbf{s}}
\newcommand{\CCC}{\mathbf{C}}
\newcommand{\LLL}{\mathbf{L}}
\newcommand{\llll}{\mathbf{l}}

\newcommand{\cov}{\operatorname{cov}}
\newcommand{\law}{\stackrel{(d)}{=}}

\newcommand{\ori}{\check}
\newcommand{\suc}{\operatorname{succ}}
\newcommand{\fISE}{f_{\mathrm{ISE}}}
\newcommand{\IP}{\operatorname{IP}}
\newcommand{\fl}{\mathit{fl}}
\newcommand{\of}{\mathit{of}}
\newcommand{\dish}{\delta_{\mathcal{H}}}
\newcommand{\diam}{\operatorname{diam}}
\newcommand{\dGH}{d_{\mathrm{GH}}}
\newcommand{\di}{d_{\infty}}
\newcommand{\PM}{\mathrm{P}\mathbb{M}}
\newcommand{\rr}{\dot}
\newcommand{\g}{\gamma n^{1/4}}
\newcommand{\gi}{\gamma^{-1} n^{-1/4}}
\newcommand{\ga}{p}
\newcommand{\ks}{{2(6g-3)}}

\newcommand{\lhb}{[[}
\newcommand{\biglhb}{[[}
\newcommand{\rhb}{]]}
\newcommand{\bigrhb}{]]}
\newcommand{\lf}{\lfloor}
\newcommand{\rf}{\rfloor}
\newcommand{\lc}{\lceil}
\newcommand{\rc}{\rceil}
\newcommand{\lt}{|}
\newcommand{\rt}{|}
\newcommand{\la}{\langle\!\langle}
\newcommand{\ra}{\rangle\!\rangle}

\newtheorem{theorem}{Theorem}
\newtheorem{lem}[theorem]{Lemma}
\newtheorem{prop}[theorem]{Proposition}

\newproclaim{defi}{Definition}
\newproclaim{rem}{Remark}

\makeatother

\begin{document}
\begin{frontmatter}

\title{The topology of scaling limits of positive genus random quadrangulations\thanksref{T1}}
\runtitle{Topology of the genus $g$ {B}rownian map}
\thankstext{T1}{Supported in part by ANR-08-BLAN-0190.}

\begin{aug}
\author[A]{\fnms{J\'{e}r\'{e}mie} \snm{Bettinelli}\corref{}\ead[label=e1]{jeremie.bettinelli@normalesup.org}\ead[label=u1,url]{www.normalesup.org/\textasciitilde bettinel}}

\runauthor{J. Bettinelli}
\affiliation{Universit\'{e} Paris-Sud 11}
\address[A]{Laboratoire de Math\'{e}matiques\\
Universit\'{e} Paris-Sud 11\\
F-91405 Orsay Cedex\\
France\\
\printead{e1}\\
\printead{u1}} 
\end{aug}

\received{\smonth{12} \syear{2010}}
\revised{\smonth{4} \syear{2011}}

%
\begin{abstract}
We discuss scaling limits of large bipartite quadrangulations of
positive genus. For a given $g$, we consider, for every $n \ge1$, a
random quadrangulation $\mathfrak{q}_n$ uniformly distributed over the
set of all
rooted bipartite quadrangulations of genus $g$ with $n$ faces. We view
it as a~metric space by endowing its set of vertices with the graph
metric. As $n$ tends to infinity, this metric space, with distances
rescaled by the factor $n^{-1/4}$, converges in distribution, at least
along some subsequence, toward a limiting random metric space. This
convergence holds in the sense of the Gromov--Hausdorff topology on
compact metric spaces. We show that, regardless of the choice of the
subsequence, the limiting space is almost surely homeomorphic to the
genus $g$-torus.
\end{abstract}

%

\end{frontmatter}

\section{Introduction}

\subsection{Motivation}

The present work is a sequel to a work by Bettinelli
\cite{bettinelli10slr}, whose aim is to investigate the topology of scaling
limits for random maps of arbitrary genus. A map is a cellular
embedding of a finite graph (possibly with multiple edges and loops)
into a compact connected orientable surface without boundary,
considered up to orientation-preserving homeomorphisms. By \textit
{cellular}, we mean that the faces of the map---the connected
components of the complement of edges---are all homeomorphic to disks.
The genus of the map is defined as the genus of the surface into which
it is embedded. For technical reasons, it will be convenient to deal
with rooted maps, meaning that one of the half-edges---or oriented
edges---is distinguished.

We will particularly focus on bipartite quadrangulations: a map is a
quadrangulation if all its faces have degree $4$; it is bipartite if
each vertex can be colored in black or white,\vadjust{\goodbreak} in such a way that no
edge links two vertices that have the same color. Although in genus
$g=0$, all quadrangulations are bipartite, this is no longer true in
positive genus $g \ge1$.

A natural way to generate a large random bipartite quadrangulation of
genus $g$ is to choose it uniformly at random from the set $\Q_n$ of
all rooted bipartite quadrangulations of genus $g$ with $n$ faces, and
then consider the limit as $n$ goes to infinity. A natural setting for
this problem is to consider quadrangulations as metric spaces endowed
with their graph metric, properly rescaled by the factor $n^{-1/4}$
\cite{marckert06limit} and to study their limit in the
Gromov--Hausdorff topology~\cite{gromov99msr}. From this point of view,
the planar case $g=0$ has largely been studied during the last decade.
Le Gall~\cite{legall07tss} showed the convergence of these metric
spaces along some subsequence. It is believed that the convergence
holds without the ``along some subsequence'' part in the last sentence,
and Le Gall gave a conjecture for a limiting space to this
sequence~\cite{legall07tss}. Although the whole convergence is yet to
be proved, some information is available on the accumulation points of
this sequence. Le Gall and Paulin~\cite{legall08slb} proved that every
possible limiting metric space is almost surely homeomorphic to the
two-dimensional sphere. Miermont~\cite{miermont08sphericity} later gave
a variant proof of this fact.

We showed in~\cite{bettinelli10slr} that the convergence along some
subsequence still holds in any fixed positive genus $g$. In this work,
we show that the topology of every possible limiting space is that of
the genus $g$-torus $\TTT_g$.

\subsection{Main results}

We will work in fixed genus $g$. On the whole, we will not let it
figure in the notation, in order to lighten them. As the case $g=0$ has
already been studied, we suppose $g \ge1$.

Recall that the Gromov--Hausdorff distance between two compact metric
spaces $(\X,\delta)$ and $(\X',\delta')$ is defined by
\[
\dGH((\X,\delta),(\X',\delta'))\de\inf\{ \dish
(\varphi(\X
),\varphi'(\X')) \},
\]
where the infimum is taken over all isometric embeddings $\varphi
\dvtx \X
\to\X''$ and $\varphi'\dvtx\X'\to\X''$ of $\X$ and $\X'$ into
the same
metric space $(\X'', \delta'')$, and $\dish$ stands for the usual
Hausdorff distance between compact subsets of $\X''$. This defines a
metric on the set $\MM$ of isometry classes of compact metric spaces
\cite{burago01cmg}, Theorem~7.3.30, making it a Polish space.\setcounter{footnote}{1}\footnote
{This is a simple consequence of Gromov's compactness theorem
\cite{burago01cmg}, Theorem 7.4.15.}

For any map $\m$, we call $V(\m)$ its set of vertices. There exists on
$V(\m)$ a natural graph metric $d_\m$: for any vertices~$a$ and $b\in
V(\m)$,
the distance $d_\m(a,b)$ is defined as the number
of edges of any shortest path linking~$a$ to~$b$. The main result
of~\cite{bettinelli10slr} is the following.
\begin{prop}\label{cvq}
Let $\mathfrak{q}_n$ be uniformly distributed over the set $\Q_n$ of all
bipartite quadrangulations of genus $g$ with $n$ faces. Then, from any
increasing sequence of integers, we may extract a subsequence\vadjust{\goodbreak}
$(n_k)_{k\ge0}$ such that there exists a metric space $(\mathfrak
{q}_\infty
,d_\infty)$ satisfying
\[
\biggl( V(\mathfrak{q}_{n_k}),\frac1 {\gamma n_k^{1/4}} d_{\mathfrak
{q}_{n_k}} \biggr) \tolk
(\mathfrak{q}_\infty,d_\infty)
\]
in the sense of the Gromov--Hausdorff topology, where
\[
\gamma\de\bigl(\tfrac{8}9 \bigr)^{1/4}.
\]

Moreover, the Hausdorff dimension of the limit space $(\mathfrak
{q}_\infty
,d_\infty)$ is almost surely equal to $4$, regardless of the choice of
the sequence of integers.
\end{prop}

Remark that the constant $\gamma$ is not necessary in this statement
(simply change $\di$ into $\gamma\di$). We kept it for the sake of
consistency with~\cite{bettinelli10slr}, and because of our definition
of $\di$ later in the paper, although it is irrelevant for the moment.
Note also that, a priori, the metric space $(\mathfrak{q}_\infty
,d_\infty)$
depends on the subsequence $(n_k)_{k\ge0}$. Similarly to the planar
case, we believe that the extraction in Proposition~\ref{cvq} is not
necessary, and we conjecture the space $(\mathfrak{q}_\infty
,d^*_\infty)$ for the
limit, where $d^*_\infty$ was defined at the end of Section 6.3
in \cite
{bettinelli10slr}. We also believe that the space $(\mathfrak
{q}_\infty,d^*_\infty
)$ is somewhat universal, in the sense that we conjecture it as the
scaling limit of more general classes of random maps. More precisely,
we think that Proposition~\ref{cvq} still holds while replacing the
class of quadrangulations with some other ``reasonable'' class of maps,
as well as the constant $\gamma$, which is inherent to the class of
quadrangulations, with the appropriate constant. In particular, our
approach can be generalized to the case of $2p$-angulations, $p\ge2$,
by following the same lines as Le Gall in~\cite{legall07tss}.

We may now state our main result, which identifies the topology of
$(\mathfrak{q}
_\infty,d_\infty)$, regardless of the subsequence $(n_k)_{k\ge0}$.
\begin{theorem}\label{cvq2}
The metric space $(\mathfrak{q}_\infty,d_\infty)$ is a.s.
homeomorphic to the
$g$-torus~$\TTT_g$.
\end{theorem}

In the general picture, we rely on the same techniques as in the planar
case. The starting point is to use a bijection due to Chapuy, Marcus
and Schaeffer~\cite{chapuy07brm} between bipartite quadrangulations of
genus $g$ with $n$ faces and so-called well-labeled $g$-trees with $n$
edges. The study of the scaling limit as $n \to\infty$ of uniform
random well-labeled $g$-trees with $n$ edges was the major purpose
of~\cite{bettinelli10slr}. This study leads to the construction of a
continuum random $g$-tree, which generalizes Aldous's CRT \cite
{aldous91crt,aldous93crt}. The first step of our proof is to carry out
the analysis of Le Gall~\cite{legall07tss} in the nonplanar case and
see the space $(\qi,\di)$ as a quotient of this continuum random
$g$-tree via an equivalence relation defined in terms of Brownian
labels on it. We then adapt Miermont's approach \cite
{miermont08sphericity}, and use the notion of 1-regularity introduced
by Whyburn~\cite{whyburn35rcm} and studied by Whyburn and Begle
\cite{begle44rc,whyburn35rcm} in order to see that the genus remains the
same in the limit.\vadjust{\goodbreak}

%
%
\begin{figure}[b]
\begin{tabular}{cc}

\includegraphics{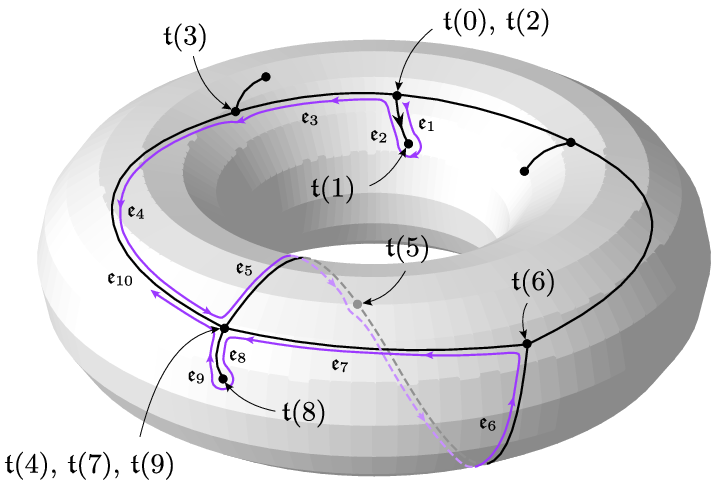}
 & \includegraphics{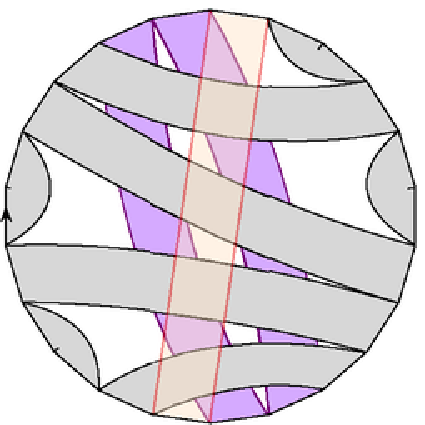}
\end{tabular}
\caption{\textup{Left.} The facial order and facial sequence of a
$g$-tree. \textup{Right.} Its representation as a polygon whose edges
are pairwise identified.}
\label{gon}
\end{figure}

Finally, we deduce the technical estimates we need from the planar case
thanks to a bijection due to Chapuy~\cite{chapuy08sum} between
well-labeled $g$-trees and well-labeled plane trees with $g$
distinguished triples of vertices.

We will use the background provided in~\cite{bettinelli10slr}. We
briefly recall it in Section~\ref{secprel}. In Section~\ref{secrgt}, we
define real $g$-trees and explain how we may see $(\mathfrak{q}_\infty
,d_\infty)$
as a quotient of such objects. Theorem~\ref{ip} in Section~\ref{secip}
gives a criteria telling which points are identified in this quotient,
and Section~\ref{secsurf} is dedicated to the proof of Theorem \ref
{cvq2}. Finally, we expose in Section~\ref{secprle} Chapuy's bijection,
and use it to prove four technical lemmas stated during Section~\ref{secip}.


\section{Preliminaries}\label{secprel}

In this section, we recall the notation, settings and results
from
\cite{bettinelli10slr} that we will need for this work. We refer the reader
to~\cite{bettinelli10slr} for more details.

We use the following formalism for maps. For any map $\m$, we denote by
$V(\m)$ and $E(\m)$, respectively,\vspace*{2pt} its sets of vertices and edges. We
also call $\vec E(\m)$ its set of half-edges, and $\mathfrak{e}_*\in
\vec E(\m)$
its root. For\vspace*{1pt} any half-edge $\mathfrak{e}$, we write $\bar\mathfrak
{e}$ its reverse---so
that $E(\m) = \{ \{ \mathfrak{e},\bar\mathfrak{e}\} \dvtx
\mathfrak
{e}\in\vec E(\m)\}$---as well as\vspace*{1pt}
$\mathfrak{e}^-$ and $\mathfrak{e}^+$ its origin and end. Finally, we
say that $\ori E(\m)
\subset\vec E(\m)$ is an orientation of the half-edges if for every
edge $\{ \mathfrak{e},\bar\mathfrak{e}\} \in E(\m)$ exactly one
of $\mathfrak{e}$ or $\bar\mathfrak{e}$
belongs to~$\ori E(\m)$.\looseness=-1

\subsection{The Chapuy--Marcus--Schaeffer bijection}\label{seccms}

The first main tool we will need consists of the
Chapuy--Marcus--Schaeffer bijection
\cite{chapuy07brm}, Corollary 2 to Theorem 1, which allows us to code (rooted)
quadrangulations by so-called well-labeled (rooted) $g$-trees.

A $g$-tree is a map of genus $g$ with only one face. This notion
naturally generalizes the notion of plane tree: in particular,
$0$-trees are plane trees. It may be convenient to represent a $g$-tree
$\tr$ with $n$ edges by a $2n$-gon whose edges are pairwise identified
(see Figure~\ref{gon}).\vadjust{\goodbreak} We note $\mathfrak{e}_1\de\mathfrak{e}_*$,
$\mathfrak{e}_2, \ldots, \mathfrak{e}
_{2n}$ the half-edges of $\tr$ arranged according to the clockwise
order around this $2n$-gon. The half-edges are said to be arranged
according to the \textit{facial order} of $\tr$. Informally, for $2
\le
i \le2n$, $\mathfrak{e}_i$ is the ``first half-edge to the left after~$\mathfrak{e}
_{i-1}$.'' We call \textit{facial sequence} of $\tr$ the sequence
$\tr
(0)$, $\tr(1), \ldots, \tr(2n)$ defined by $\tr(0)=\tr(2n) =
\mathfrak{e}_1^- =
\mathfrak{e}_{2n}^+$ and for $1 \le i \le2n-1$, $\tr(i)=\mathfrak
{e}_i^+=\mathfrak{e}_{i+1}^-$.
Imagine a fly flying along the boundary of the unique face of $\tr$.
Let it start at time~$0$ by following the root $\mathfrak{e}_*$, and
let it take
one unit of time to follow each half-edge, then $\tr(i)$ is the vertex
where the fly is at time $i$.

Let $\tr$ be a $g$-tree. Two vertices $u,v\in V(\tr)$ are said to be
\textit{neighbors}, and we write $u\sim v$, if there is an edge
linking them.
\begin{defi}
A \textit{well-labeled $g$-tree} is a pair $(\tr,\mathfrak{l})$
where $\tr$ is
a $g$-tree and $\mathfrak{l}\dvtx V(\tr) \to\Z$ is a function
(thereafter called
\textit{labeling function}) satisfying:
\begin{longlist}
\item$\mathfrak{l}(\mathfrak{e}_*^-) = 0$, where $\mathfrak{e}_*$
is the root of $\tr$;
\item if $u \sim v$, then $|\mathfrak{l}(u) - \mathfrak{l}(v)| \le1$.
\end{longlist}
\end{defi}

We call $\T_n$ the set of all well-labeled $g$-trees with $n$ edges. A
\textit{pointed quadrangulation} is a pair $(\mathfrak{q},v^\bullet
)$ consisting
in a quadrangulation $\mathfrak{q}$ together with a vertex $v^\bullet
\in V(\mathfrak{q})$.
We call $\Qb_n \de\{(\mathfrak{q},v^\bullet)\dvtx \mathfrak{q}\in
\Q_n, v^\bullet\in
V(\mathfrak{q}) \}$ the set of all pointed bipartite
quadrangulations of
genus $g$ with $n$ faces.

The Chapuy--Marcus--Schaeffer bijection is a bijection between the sets
$\T_n \times\{-1,+1\}$ and $\Qb_n$. We briefly describe here the
mapping from $\T_n \times\{-1,+1\}$ onto $\Qb_n$, and we refer the
reader to~\cite{chapuy07brm} for a more precise description. Let $(\tr
,\mathfrak{l}) \in\T_n$ be a well-labeled $g$-tree with $n$ edges
and $\eps
_\pm\in\{-1,+1\}$. As above, we write $\tr(0)$, $\tr(1), \ldots,
\tr
(2n)$ its facial sequence. The pointed quadrangulation $(\mathfrak
{q},v^\bullet)$
corresponding to $((\tr,\mathfrak{l}),\eps_\pm)$ is then
constructed as
follows. First, shift all the labels in such a way that the minimal
label is equal to~$1$. Let us call $\tilde\mathfrak{l}\de\mathfrak
{l}-\min\mathfrak{l}+1$
this shifted labeling function. Then, add an extra vertex $v^\bullet$
carrying the label $\tilde\mathfrak{l}(v^\bullet)\de0$ inside the
only face
of $\tr$. Finally, following the facial sequence, for every $0\le i
\le
2n-1$, draw an arc---without crossing any edge of $\tr$ or arc already
drawn---between $\tr(i)$ and $\tr(\suc(i))$, where $\suc(i)$ is the
\textit{successor} of $i$, defined by
%
%
\begin{equation}\label{suc}
\suc(i) \de
\cases{
\inf\{k\ge i \dvtx \tilde\mathfrak{l}(\tr(k))= \ell\}, &\quad if
$\{k\ge i
\dvtx \tilde\mathfrak{l}(\tr(k))= \ell\} \neq\varnothing$,\cr
\inf\{k\ge1 \dvtx \tilde\mathfrak{l}(\tr(k))=\ell\}, &\quad otherwise,}
\end{equation}
where $\ell=\tilde\mathfrak{l}(\tr(i)) - 1$, and
with the conventions $\inf\varnothing= \infty$, and $\tr(\infty) =
v^\bullet$.

\begin{figure}

\includegraphics{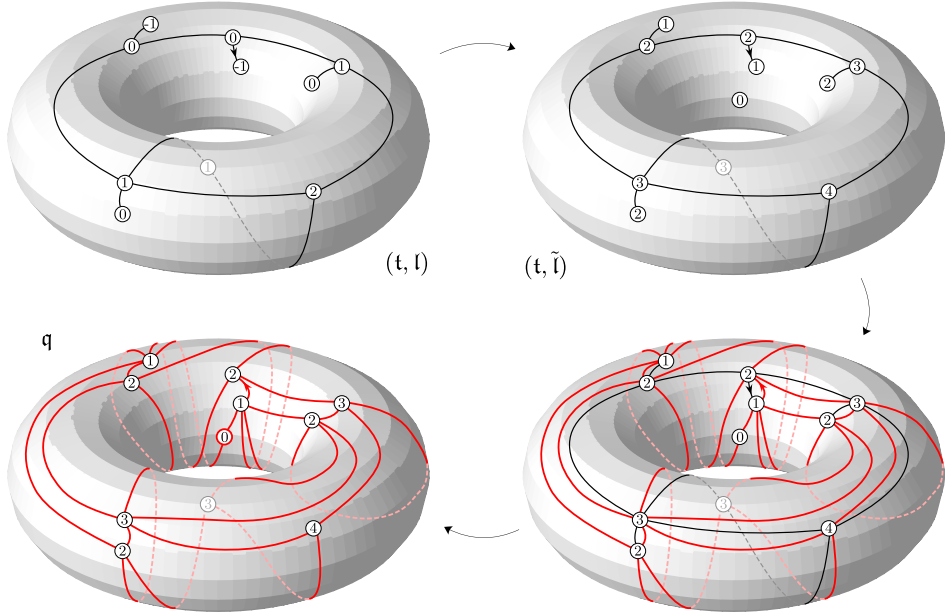}

\caption{The Chapuy--Marcus--Schaeffer bijection. In this example,
$\eps_\pm= +1$. On the bottom--left picture, the vertex $v^\bullet$ has
a thicker (red) borderline.}
\label{cms}
\end{figure}

The quadrangulation $\mathfrak{q}$ is then defined as the map whose
set of
vertices is $V(\tr) \cup\{v^\bullet\}$, whose edges are the arcs we
drew and whose root is the first arc drawn, oriented \textit{from}
$\tr
(0)$ if $\eps_\pm=-1$ or \textit{toward} $\tr(0)$ if $\eps_\pm
=+1$; see
Figure~\ref{cms}.

%
%

Because of the way we drew the arcs of $\mathfrak{q}$, we see that for
any vertex
$v\in V(\mathfrak{q})$, $\tilde\mathfrak{l}(v)=d_\mathfrak
{q}(v^\bullet,v)$. When seen as a vertex
in $V(\mathfrak{q})$, we write $\mathfrak{q}(i)$ instead of $\tr
(i)$. In particular, $\{\mathfrak{q}
(i),0 \le i \le2n\} = V(\mathfrak{q})\bs\{v^\bullet\}$.

We end this section by giving an upper bound for the distance between
two vertices $\mathfrak{q}(i)$ and $\mathfrak{q}(j)$, in terms of the
labeling function
$\mathfrak{l}$:
%
%
\begin{eqnarray}\label{dlemme}
&&
d_\mathfrak{q}(\mathfrak{q}(i),\mathfrak{q}(j)) \nonumber\\[-8pt]\\[-8pt]
&&\qquad\le\mathfrak
{l}(\tr(i)) + \mathfrak{l}(\tr(j)) - 2 \max\Bigl(\min_{\mbox{\fontsize{8.36pt}{8.36pt}\selectfont{$k \mbox{$\in$}
\overrightarrow{\lbracket i,j \rbracket}$}}} \mathfrak{l}(\tr
(k)),\min_{ \mbox{\fontsize{8.36pt}{8.36pt}\selectfont{$k \mbox{$\in$}
\overrightarrow
{\lbracket j,i \rbracket}$}}} \mathfrak{l}(\tr(k)) \Bigr)+2,\nonumber
\end{eqnarray}
where we note, for $i\le j$, $\lbracket i,j \rbracket \de[i,j] \cap
\Z= \{ i, i+1,
\ldots, j \}$, and
%
%
\begin{equation}\label{oraij}
\overrightarrow{\lbracket i,j \rbracket} \de
\cases{
\lbracket i,j \rbracket, &\quad if $i \le j$,\cr
\lbracket i,2n \rbracket \cup\lbracket0,j \rbracket, &\quad if $j < i$.}
\end{equation}
We refer the reader to~\cite{miermont09trm}, Lemma 4, for a detailed
proof of this bound.


\subsection{Decomposition of a $g$-tree}\label{decomp}

We explained in~\cite{bettinelli10slr} how to decompose a~$g$-tree into
simpler objects. Roughly speaking, a $g$-tree is a scheme (a $g$-tree
whose all vertices have degree at least $3$) in which every half-edge
is replaced by a forest.

\subsubsection{Forests}\label{secfor}

We adapt the standard formalism for plane trees---as found in
\cite{neveu86apg} for instance---to forests. Let us call $\U\de\bigcup
_{n=1}^{\infty} \N^n$, where $\N\de\{1,2,\ldots\}$. If $u \in\N^n$,
we write $|u| \de n$. For $u=(u_1,\ldots,u_n)$, \mbox{$v=(v_1,\ldots,v_p) \in
\U
$}, we let $uv \de(u_1,\ldots,u_n,v_1,\ldots,v_p)$ be the concatenation
of $u$ and $v$. If $w=uv$ for some $u,v\in\U$, we say that $u$ is
a~\textit{ancestor} of $w$ and that $w$ is a \textit{descendant} of $u$.
In the case where $v \in\N$, we may also use the terms \textit{parent}
and \textit{child} instead.
\begin{defi}
A \textit{forest} is a finite subset $\mathfrak{f}\subset\U$
satisfying the following:
\begin{longlist}
\item there is an integer $t(\mathfrak{f}) \ge1$ such that $\mathfrak
{f}\cap\N= \lbracket1,t(\mathfrak{f})+1 \rbracket$;
\item if $u \in\mathfrak{f}$, $|u| \ge2$, then its parent belongs
to $\mathfrak{f}$;
\item for every $u \in\mathfrak{f}$, there is an integer
$c_u(\mathfrak{f}) \ge0$ such
that $ui \in\mathfrak{f}$ if and only if $1 \le i \le c_u(\mathfrak{f})$;
\item$c_{t(\mathfrak{f})+1}(\mathfrak{f}) =0$.
\end{longlist}

The integer $t(\mathfrak{f})$ is called the \textit{number of trees}
of $\mathfrak{f}$.
\end{defi}

For $u=(u_1,\ldots,u_p) \in\mathfrak{f}$, we call $\aaa(u) \de u_1$
its oldest
ancestor. A \textit{tree} of $\mathfrak{f}$ is a level set for $\aaa
$: for $1\le
j \le t(\mathfrak{f})$, the $j$th tree of $\mathfrak{f}$ is the set
$\{u\in\mathfrak{f} \dvtx \aaa
(u)=j\}$. The integer $\aaa(u)$ hence records which tree $u$ belongs to.
We call $\mathfrak{f}\cap\N= \{\aaa(u), u \in\mathfrak{f}\}
$ the \textit{floor} of
the forest $\mathfrak{f}$.

\begin{figure}

\includegraphics{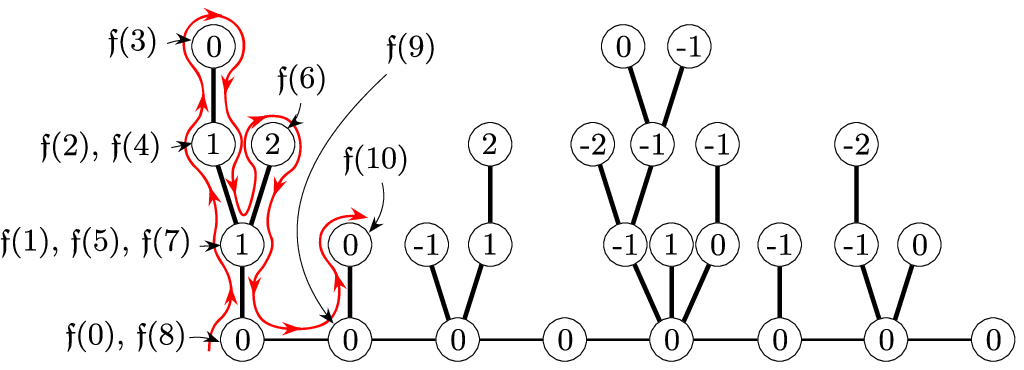}

\caption{The facial sequence of a well-labeled forest from $\fF_7^{20}$.}
\label{dfs}\vspace*{-3pt}
\end{figure}

For $u,v \in\mathfrak{f}$, we write $u \sim v$ if either $u$ is a
parent or
child of $v$, or $u,v \in\N$ and $|u-v| = 1$. It is convenient, when
representing a forest, to draw edges between $u$'s and $v$'s such that $u
\sim v$; see Figure~\ref{dfs}. We say that an edge drawn between a
%
%
parent and its child is a \textit{tree edge} whereas an edge drawn
between two consecutive tree roots, that is, between some $i$ and
$i+1$, will be called a \textit{floor edge}. We call $\F_\sigma^m\de
\{\mathfrak{f} \dvtx t(\mathfrak{f})=\sigma, |\mathfrak{f}|=
m+\sigma+1 \}$ the set of all
forests with $\sigma$ trees and $m$ tree edges.
\begin{defi}
A \textit{well-labeled forest} is a pair $(\mathfrak{f},\mathfrak
{l})$ where $\mathfrak{f}$ is a
forest, and $\mathfrak{l}\dvtx\mathfrak{f}\to\Z$ is a function satisfying:
\begin{longlist}
\item for all $u\in\mathfrak{f}\cap\N$, $\mathfrak{l}(u) = 0$;
\item if $u \sim v$, $|\mathfrak{l}(u) - \mathfrak{l}(v)| \le1$.
\end{longlist}
\end{defi}

Let $\fF_\sigma^m \de\{(\mathfrak{f},\mathfrak{l})\dvtx \mathfrak
{f}\in\F_\sigma^m \}$ be the
set of well-labeled forests with $\sigma$ trees and $m$ tree edges.\vadjust{\goodbreak}

\paragraph*{\quad Encoding by contour and spatial contour functions}

There is a very convenient way to code forests and well-labeled
forests. Let $\mathfrak{f}\in\F_\sigma^m$ be a forest. Let us begin
by defining
its \textit{facial sequence} $\mathfrak{f}(0),\mathfrak{f}(1),\ldots
,\mathfrak{f}(2m+\sigma)$ as
follows (see Figure~\ref{dfs}): $\mathfrak{f}(0) \de1$, and for
$0\le i \le
2m+\sigma-1$:
\begin{itemize}[$\diamond$]
\item[$\diamond$] if $\mathfrak{f}(i)$ has children that do not appear in the
sequence $\mathfrak{f}
(0),\mathfrak{f}(1),\ldots,\mathfrak{f}(i)$, then $\mathfrak
{f}(i+1)$ is the first of these children,
that is, $\mathfrak{f}(i+1)\de\mathfrak{f}(i)j_0$ where
\[
j_0 = \min\bigl\{j\ge1 \dvtx \mathfrak{f}(i)j \notin\{\mathfrak
{f}(0),\mathfrak{f}(1),\ldots,\mathfrak{f}(i) \}
\bigr\};
\]
\item[$\diamond$] otherwise, if $\mathfrak{f}(i)$ has a parent [i.e., $|\mathfrak
{f}(i)| \ge2$], then
$\mathfrak{f}(i+1)$ is this parent;
\item[$\diamond$] if neither of these cases occur, which implies that $|\mathfrak
{f}(i)|=1$,
then $\mathfrak{f}(i+1) \de\mathfrak{f}(i)+1$.
\end{itemize}

%
%
\begin{figure}

\includegraphics{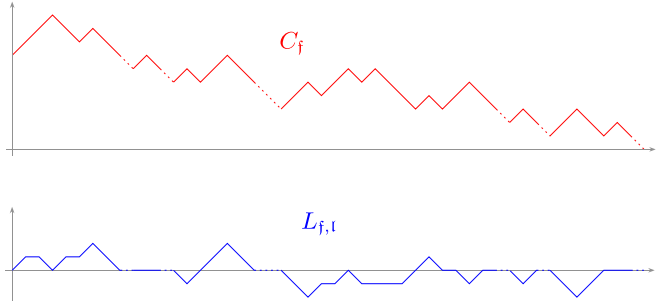}

\caption{The contour pair of the well-labeled forest appearing in
Figure \protect\ref{dfs}. The paths are dashed on the intervals corresponding
to floor edges.}
\label{wlf}\vspace*{-3pt}
\end{figure}

Each tree edge is visited exactly twice---once going from the parent to
the child, once going the other way around---whereas each floor edge is
visited only once---from some $i$ to $i+1$. As a result, $\mathfrak
{f}(2m+\sigma)
= t(\mathfrak{f}) + 1$.

The \textit{contour pair} $(C_\mathfrak{f},L_{\mathfrak{f},\mathfrak
{l}})$ of $(\mathfrak{f},\mathfrak{l})$ consists
in the \textit{contour function} $C_\mathfrak{f}\dvtx[0,\break 2m+\sigma
]\to\R
_+$ of $\mathfrak{f}$
and the \textit{spatial contour function} $L_{\mathfrak{f},\mathfrak
{l}}\dvtx[0,2m+\sigma]\to
\R$ defined by
\[
C_\mathfrak{f}(i) \de|\mathfrak{f}(i)| + t(\mathfrak{f}) -\aaa(
\mathfrak{f}(i) )\quad\mbox{and}\quad L_{\mathfrak{f},\mathfrak{l}}(i) \de
\mathfrak{l}(\mathfrak{f}(i)),\qquad 0 \le i \le2m+\sigma,
\]
and linearly interpolated between integer values (see Figure
\ref{wlf}).\vspace*{-2pt}
It entirely determines $(\mathfrak{f},\mathfrak{l})$.


\subsubsection{Decomposition of a well-labeled $g$-tree into
simpler objects}\label{secdec}

We explain here how to decompose a well-labeled $g$-tree. See
\cite{bettinelli10slr} for a more precise description.\vspace*{-2pt}
\begin{defi}
We call \textit{scheme} of genus $g$ a $g$-tree with no vertices of
degree one or two. A scheme is said to be \textit{dominant} when it
only has vertices of degree exactly three.\vspace*{-2pt}
\end{defi}

We call $\Sg$ the finite set of all schemes of genus $g$ and $\Sg^*$
the set of all dominant schemes of genus $g$.\vadjust{\goodbreak}

Let us first explain how to decompose a~$g$-tree (without labels) into
a~scheme, a~family of forests and an integer. Let $\s$ be a scheme. We
suppose that we have forests $\mathfrak{f}^\mathfrak{e}\in\F
_{\sigma^\mathfrak{e}}^{m^\mathfrak{e}}$, $\mathfrak{e}\in
\vec E(\s)$, where for all $\mathfrak{e}$, $\sigma^{\bar\mathfrak
{e}}=\sigma^\mathfrak{e}$, as well
as an integer $u \in\lbracket0$, $2m^{\mathfrak{e}_*} + \sigma
^{\mathfrak {e}_*}-1 \rbracket$, where $\mathfrak{e}
_*$ denotes the root of $\s$. We construct a $g$-tree as follows.
First, we replace\vspace*{1pt} every edge $\{\mathfrak{e},\bar\mathfrak{e}\}$ in
$\s$ with a chain of
$\sigma^\mathfrak{e}= \sigma^{\bar\mathfrak{e}}$ edges. Then, for
every half-edge $\mathfrak{e}\in
\vec E(\s)$, we replace the chain of half-edges corresponding to it
with the forest $\mathfrak{f}^\mathfrak{e}$, in such a way that its
floor matches with the
chain. In other words, we ``graft'' the forest $\mathfrak{f}^\mathfrak
{e}$ to the left
of $\mathfrak{e}$. Finally, the root of the $g$-tree is the half-edge
linking $\mathfrak{f}
^{\mathfrak{e}_*}(u)$ to $\mathfrak{f}^{\mathfrak{e}_*}(u+1)$ in the
forest grafted to the left of $\mathfrak{e}_*$.
\begin{prop}\label{decompnolab}
The above description provides a bijection between the set of all
$g$-trees and the set of all triples $(\s,(\mathfrak{f}^\mathfrak
{e})_{\mathfrak{e}\in\vec E(\s
)},u)$ where $\s\in\Sg$ is a~scheme (of genus $g$), the forests
$\mathfrak{f}
^\mathfrak{e}\in\F_{\sigma^\mathfrak{e}}^{m^\mathfrak{e}}$ are
such that $\sigma^{\bar\mathfrak{e}}=\sigma^\mathfrak{e}$
for all $\mathfrak{e}$ and $u \in\lbracket0,2m^{\mathfrak{e}_*} +
\sigma^{\mathfrak{e}_*}-1 \rbracket$.

Moreover, $g$-trees with $n$ edges correspond to triples satisfying the
condition $\sum_{\mathfrak{e}\in\vec E(\s)} (m^\mathfrak{e}+
\frac1 2 \sigma^\mathfrak{e})= n$.
\end{prop}

\begin{figure}[b]

\includegraphics{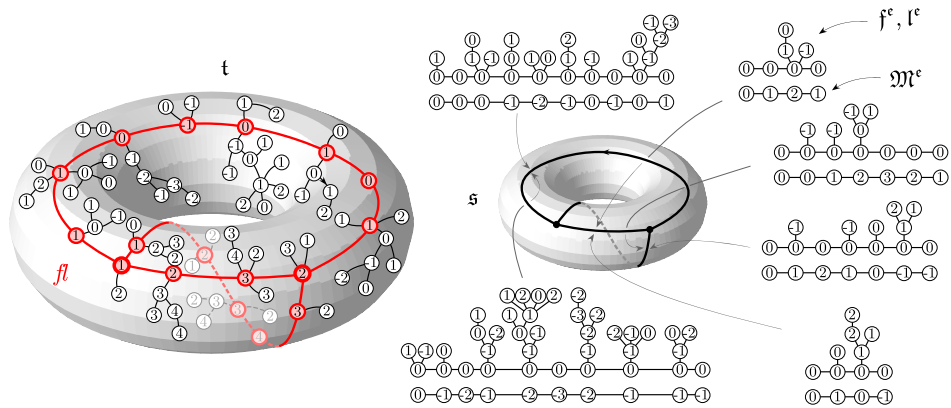}

\caption{Decomposition of a well-labeled $g$-tree $\tr$ into its
scheme $\s$, the collection of its Motzkin paths $(\mathfrak
{M}^\mathfrak{e})_{\mathfrak{e}\in\vec
E(\s)}$ and the collection of its well-labeled forests $(\mathfrak
{f}^{\mathfrak{e}},\mathfrak{l}^\mathfrak{e}
)_{\mathfrak{e}\in\vec E(\s)}$. In this example, the integer $u=
10$. The floor
of $\tr$ is more thickly outlined, and its two nodes are even more
thickly outlined.}
\label{structure}
\end{figure}

Let $\tr$ be a $g$-tree and $(\s,(\mathfrak{f}^\mathfrak
{e})_{\mathfrak{e}\in\vec E(\s)},u)$ be
the corresponding triple. We say that $\s$ is the scheme of $\tr$ and
that the forests $\mathfrak{f}^\mathfrak{e}$, $\mathfrak{e}\in\vec
E(\s)$, are its forests. The set
$V(\s)$ may be seen as a subset of $\tr$; we call \textit{nodes} its
elements. Finally, we call \textit{floor} of $\tr$ the set $\fl$ of
vertices we obtain after replacing the edges of~$\s$ by chains of edges
(see Figure~\ref{structure}).

%
%

We now deal with well-labeled $g$-trees. We will need the following definition:
\begin{defi}
We call \textit{Motzkin path} a sequence $(M_n)_{0 \le n \le\sigma}$
for some $\sigma\ge0$ such that $M_0=0$ and for $0 \le i \le\sigma
-1$, $M_{i+1} - M_i \in\{-1,0,1\}$. We write $\sigma(M) \de\sigma$
its lifetime.
\end{defi}

Let $\s$ be a scheme. We suppose that we have well-labeled forests
$(\mathfrak{f}
^\mathfrak{e},\mathfrak{l}^\mathfrak{e})\in\fF_{\sigma^\mathfrak
{e}}^{m^\mathfrak{e}}$, $\mathfrak{e}\in\vec E(\s)$, where for
all $\mathfrak{e}$, $\sigma^{\bar\mathfrak{e}}=\sigma^\mathfrak
{e}$, as well as an integer $u \in
\lbracket0,2m^{\mathfrak{e}_*} + \sigma^{\mathfrak{e}_*} -1
\rbracket$. Suppose
moreover that we have a
family of Motzkin\vspace*{-1pt} paths $(\mathfrak{M}^\mathfrak{e})_{\mathfrak
{e}\in\vec E(\s)}$ such that $\mathfrak{M}^\mathfrak{e}
$ is defined on $\lbracket0,\sigma^\mathfrak{e} \rbracket$ and
$\mathfrak
{M}^\mathfrak{e}(\sigma^\mathfrak{e}) = l^{\mathfrak{e}^+}
- l^{\mathfrak{e}^-}$ for some family of integers $(l^v)_{v\in V(\s
)}$ with $l^{\mathfrak{e}
_*^-}=0$. We suppose that the Motzkin paths satisfy the following relation:
\[
\mathfrak{M}^{\bar\mathfrak{e}}(i) = \mathfrak{M}^\mathfrak
{e}(\sigma^\mathfrak{e}-i) - l^\mathfrak{e},\qquad 0\le i \le\sigma
^\mathfrak{e}\qquad
\mbox{where } l^\mathfrak{e}\de l^{\mathfrak{e}^+} -
l^{\mathfrak{e}^-}.
\]
We will say that a quadruple $( \s, (\mathfrak{M}^\mathfrak
{e})_{\mathfrak{e}\in\vec E(\s)}, (\mathfrak{f}
^{\mathfrak{e}},\mathfrak{l}^\mathfrak{e})_{\mathfrak{e}\in\vec
E(\s)},u )$ satisfying these
constraints is \textit{compatible}.

We construct a well-labeled $g$-tree as follows. We begin by suitably
relabeling the forests. For every half-edge $\mathfrak{e}$, first, we
shift the
labels of $\mathfrak{M}^\mathfrak{e}$ by $l^{\mathfrak{e}^-}$ so
that it goes from $l^{\mathfrak{e}^-}$ to
$l^{\mathfrak{e}^+}$. Then we shift all the labels of $(\mathfrak
{f}^\mathfrak{e},\mathfrak{l}^\mathfrak{e})$ tree by
tree according to the Motzkin path: precisely, we change $\mathfrak
{l}^\mathfrak{e}$ into
$w\in\mathfrak{f}^\mathfrak{e}\mapsto l^{\mathfrak{e}^-} +
\mathfrak{M}^\mathfrak{e}(\aaa(w)-1) + \mathfrak{l}^\mathfrak
{e}(w)$. Then we
replace the half-edge $\mathfrak{e}$ with this forest, as in the previous
section. As before, we find the position of the root thanks to $u$.
Finally, we shift all the labels for the root label to be equal to
$0$.\looseness=-1

\begin{prop}\label{decwl}
The above description provides a bijection between the set of all
well-labeled $g$-trees and the set of all compatible quadruples.

Moreover, $g$-trees with $n$ edges correspond to quadruples satisfying
the condition $\sum_{\mathfrak{e}\in\vec E(\s)} (m^\mathfrak
{e}+ \frac1 2 \sigma^\mathfrak{e}
)= n$.
\end{prop}

If we call $(C^\mathfrak{e},L^\mathfrak{e})$ the contour pair of
$(\mathfrak{f}^\mathfrak{e},\mathfrak{l}^\mathfrak{e})$, then we
may retrieve the oldest ancestor of $\mathfrak{f}^\mathfrak{e}(i)$
thanks to $C^\mathfrak{e}$ by the relation
\[
\aaa(\mathfrak{f}^\mathfrak{e}(i)) -1 = \sigma^\mathfrak
{e}- \underline C^\mathfrak{e}(i),
\]
where we use the notation
\[
\underline{X}_s \de\inf_{[0, s]} X
\]
for any process $(X_s)_{s\ge0}$. The function
%
%
\begin{equation}\label{lele}
\Lab^\mathfrak{e}\de\bigl( L^\mathfrak{e}(t) + \mathfrak
{M}^\mathfrak{e}\bigl( \sigma^\mathfrak{e}- \underline C^\mathfrak{e}(t)
\bigr) \bigr)_{0\le t \le2m^\mathfrak{e}+\sigma^\mathfrak{e}},
\end{equation}
then records the labels of the forest $\mathfrak{f}^\mathfrak{e}$,
once shifted tree by
tree according to the Motzkin path $\mathfrak{M}^\mathfrak{e}$. This
function will be used
in Section~\ref{secset}.

Through the Chapuy--Marcus--Schaeffer bijection, a uniform random
quadrangulation corresponds to a uniform random well-labeled $g$-tree.
It can then be decomposed into a scheme, a collection of well-labeled
forests, a collection of Motzkin paths and an integer, as explained
above. The following section exposes the scaling limits of these objects.\vadjust{\goodbreak}


\subsection{Scaling limits}\label{secsl}

Let us define the space $\K$ of continuous real-valued functions on
$\R
_+$ killed at some time
\[
\K\de\bigcup_{x \in\R_+} \C([0, x],\R).
\]
For an element $f \in\K$, we will define its lifetime $\sigma(f)$ as
the only $x$ such that $f \in\C([0, x],\R)$. We endow this space with
the following metric:
\[
d_\K(f,g) \de|\sigma(f) - \sigma(g) | + \sup_{y \ge0} \bigl\lt f
\bigl(y\wedge\sigma(f)\bigr)-g\bigl(y\wedge\sigma(g)\bigr)\bigr\rt.
\]

Throughout this section, $m$ and $\sigma$ will denote positive real
numbers and~$l$ will be any real number.


\subsubsection{Brownian bridges, first-passage Brownian bridges and
Brownian snake}

We define here the Brownian bridge $B_{[0,m]}^{0\to l}$ on $[0,m]$
from $0$ to $l$ and the first-passage Brownian bridge $F_{[0,m]}^{0\to
-\sigma}$ on $[0,m]$ from $0$ to $-\sigma$. Informally, $B_{[0,m]}^{0
\to l}$ and $F_{[0,m]}^{0 \to-\sigma}$ are a standard Brownian motion
$\beta$ on $[0,m]$ conditioned, respectively, on the events $\{\beta_m =
l\}$ and $\{\inf\{s\ge0 \dvtx \beta_s=-\sigma\}=m\}$. Because both
theses events occur with probability $0$, we need to define these
objects properly. There are several equivalent ways to do so; see for
example~\cite{bertoin03ptf,billingsley68cpm,revuz99cma}. We call
$\ga
_a$ the density of a centered Gaussian variable with variance~$a$, as
well as $\ga'_a$ its derivative
\[
\ga_a(x) \de\frac1 {\sqrt{2\pi a}} \exp\biggl(- \frac{x^2} {2a}
\biggr)\quad\mbox{and}\quad\ga'_a(x) = - \frac x a \ga_a (x ).
\]

Let $(\beta_t)_{0\le t\le m}$ be a standard Brownian motion. As
explained in~\cite{fitzsimmons93mbc}, Proposition~1, the law of the
Brownian bridge is characterized by the equation $B_{[0,m]}^{0\to l}(m)
= l$ and the formula
\[
\mathbb{E}\bigl[ f \bigl(\bigl(B_{[0,m]}^{0\to l}(t)\bigr)_{0 \le t \le m'} \bigr) \bigr] = \mathbb
{E}\biggl[ f ( (\beta_t)_{0 \le t \le m'} )\frac{\ga_{m - m'}(l-\beta
_{m'})}{\ga_m(l)} \biggr]
\]
for all bounded measurable functions $f$ on $\K$, for all $0\le m' <
m$. We define the law of the first-passage Brownian bridge in a similar
way, by letting
%
%
\begin{eqnarray}\label{f}
&&\mathbb{E}\bigl[ f \bigl(\bigl(F_{[0,m]}^{0\to-\sigma}(t)\bigr)_{0 \le t \le m'} \bigr)
\bigr]\nonumber\\[-8pt]\\[-8pt]
&&\qquad= \mathbb{E}\biggl[ f ((\beta_t)_{0 \le t \le m'} )\frac{\ga_{m -
m'}'(-\sigma - \beta_{m'})}{\ga_m'(-\sigma)} \mathbh{1}_{\{
{\underline{\beta}}_{m'} > - \sigma\}} \biggr]
\nonumber
\end{eqnarray}
for all bounded measurable functions $f$ on $\K$, for all $0\le m' < m$
and\break $F_{[0,m]}^{0\to-\sigma}(m) = -\sigma$.\vspace*{1pt}

For any real numbers $l_1$, $l_2$, $\sigma_1 > \sigma_2$, we define the
bridge on $[0,m]$ from~$l_1$ to $l_2$ and the first-passage bridge on
$[0,m]$ from $\sigma_1$ to $\sigma_2$ by
\[
B_{[0,m]}^{l_1 \to l_2}(s) \de l_1 + B_{[0,m]}^{0 \to l_2-l_1}
\quad\mbox{and}\quad
F_{[0,m]}^{\sigma_1 \to\sigma_2} \de\sigma_1 + F_{[0,m]}^{0 \to
\sigma_2-\sigma_1}.
\]
See~\cite{bettinelli10slr}, Section 5.1, for a more precise description
of these objects. In particular,
\cite{bettinelli10slr}, Lemmas 10 and 14, show that these objects appear as the limits
of their discrete analogs.

Conditionally given a first-passage Brownian bridge $F =
F_{[0,m]}^{\sigma\to0}$, we define a Gaussian process $
(Z_{[0,m]}(s))_{0\le s \le m}$ with covariance function
\[
\cov\bigl(Z_{[0,m]}(s), Z_{[0,m]}(s')\bigr) = \inf_{[s\wedge s',s\vee
s']} (F - \underline{F}).
\]
The process $( F_{[0,m]}^{\sigma\to0}, Z_{[0,m]})$ has the
law of the so-called Brownian snake's head; see
\cite{duquesne02rtl,legall99sbp} for more details.


\subsubsection{Convergence results}

Recall that $\Sg^*$ is the set of all dominant schemes of genus $g$,
that is, schemes with only vertices of degree $3$. For any $\s\in\Sg$,
we identify an element $(m,\sigma,l,u) \in\R_+^{\vec E(\s)\bs\{
\mathfrak{e}_*\}
}\times(\R_+^*)^{\ori E(\s)}\times\R^{V(\s)\bs\{\mathfrak
{e}^-_*\}}\times\R_+$
with an element of $\R_+^{\vec E(\s)}\times(\R_+^*)^{\vec E(\s
)}\times
\R^{V(\s)}\times\R_+$ by setting:
\begin{itemize}[$\diamond$]
\item[$\diamond$]$m^{\mathfrak{e}_*} \de1 - \sum_{\mathfrak{e}\in\vec E(\s
)\bs\{\mathfrak{e}_*\}} m^\mathfrak{e}$,
\item[$\diamond$]$\sigma^{\bar\mathfrak{e}}\de\sigma^{\mathfrak{e}}$ for
every $\mathfrak{e}\in\ori E(\s)$,
\item[$\diamond$]$l^{\mathfrak{e}_*^-}\de0$.
\end{itemize}

We write
\[
\Delta_\s\de\biggl\{(x_\mathfrak{e})_{\mathfrak{e}\in\vec E(\s)} \in
[0,1]^{\vec E(\s)},
\sum_{\mathfrak{e}\in\vec E(\s)} x_\mathfrak{e}=1 \biggr\},
\]
the simplex of dimension $|\vec E(\s)|-1$. Note that $m$ lies in
$\Delta
_\s$ as long as $m^{\mathfrak{e}_*} \ge0$. We define the probability
$\mu$ by,
for all measurable function $\varphi$ on $\bigcup_{\s\in\Sg} \{\s
\}
\times\Delta_\s\times(\R_+^*)^{\vec E(\s)}\times\R^{V(\s
)}\times[0,1]$,
\begin{eqnarray*}
\mu(\varphi) &=& \frac1\Upsilon\sum_{\s\in\Sg^*} \int_{\mathcal
{S}^\s}
d\mathcal{L}^\s\,
\mathbh{1}_{\{m^{\mathfrak{e}_*}\ge0, u< m^{\mathfrak{e}_*}\}}
\varphi(\s, m,\sigma, l,u
)\\
&&\hphantom{\frac1\Upsilon\sum_{\s\in\Sg^*} \int_{\mathcal
{S}^\s}}
{}\times\prod_{\mathfrak{e}\in{\vec E}(\s)} - \ga'_{m^\mathfrak{e}} (
\sigma^\mathfrak{e})
\prod_{\mathfrak{e}\in\ori{E}(\s)} \ga_{\sigma^\mathfrak{e}}
(l^\mathfrak{e}),
\end{eqnarray*}
where $l^\mathfrak{e}\de l^{\mathfrak{e}^+} - l^{\mathfrak{e}^-}$,
the measure $d\mathcal{L}^\s=
d(m^\mathfrak{e})\, d(\sigma^\mathfrak{e}) \,d(l^v) \,du$ is the
Lebesgue measure on the set
\[
\mathcal{S}^\s\de[0,1]^{\vec E(\s)\bs\{\mathfrak{e}_*\}}\times
(\R_+^*)^{\ori
E(\s)}\times\R^{V(\s)\bs\{\mathfrak{e}^-_*\}}\times[0,1]
\]
and
%
%
\begin{equation}\label{upsi}
\Upsilon= \sum_{\s\in\Sg^*} \int_{\mathcal{S}^\s} d\mathcal
{L}^\s \,\mathbh{1}_{\{m^{\mathfrak{e}_*}\ge0, u< m^{\mathfrak
{e}_*}\}}
\prod_{\mathfrak{e}\in{\vec E}(\s)} - \ga'_{m^\mathfrak{e}} (
\sigma^\mathfrak{e})
\prod_{\mathfrak{e}\in\ori{E}(\s)} \ga_{\sigma^\mathfrak{e}}
(l^\mathfrak{e})
\end{equation}
is a normalization constant. We gave a nonintegral expression for this
constant in~\cite{bettinelli10slr}.\vadjust{\goodbreak}

Let $(\tr_n,\mathfrak{l}_n)$ be uniformly distributed over the set
$\T_n$ of
well-labeled $g$-trees with $n$ vertices. We call $\s_n$ its scheme and
we define, as in Section~\ref{decomp}, $(\mathfrak{f}_n^\mathfrak
{e},\mathfrak{l}_n^\mathfrak{e})_{\mathfrak{e}\in\vec
E(\s_n)}$ its well-labeled forests, $(m_n^\mathfrak{e})_{\mathfrak
{e}\in\vec E(\s_n)}$ and
$(\sigma_n^\mathfrak{e})_{\mathfrak{e}\in\vec E(\s_n)}$,
respectively, their sizes and
lengths, $(l_n^v)_{v\in V(\s_n)}$ the shifted labels of its nodes,
$(\mathfrak{M}
_n^{\mathfrak{e}})_{\mathfrak{e}\in\vec E(\s_n)}$ its Motzkin paths
and $u_n$ the integer
recording the position of the root in the first forest $\mathfrak
{f}_n^{\mathfrak{e}_*}$.
We call $(C_n^\mathfrak{e},L_n^\mathfrak{e})$ the contour pair of the
well-labeled forest
$(\mathfrak{f}_n^\mathfrak{e},\mathfrak{l}_n^\mathfrak{e})$, and we
extend the definition of $\mathfrak{M}_n^\mathfrak{e}$ to
$[0,\sigma_n^\mathfrak{e}]$ by linear interpolation. We then define
the rescaled
versions of these objects [recall that $\gamma\de(8/9)^{1/4}$]
\begin{eqnarray*}
m_{(n)}^\mathfrak{e}&\de&\frac{2m_n^\mathfrak{e}+ \sigma_n^\mathfrak
{e}}{2n},\qquad
\sigma_{(n)}^\mathfrak{e}\de\frac{\sigma_n^\mathfrak{e}}{\sqrt
{2n}},\qquad
l_{(n)}^v \de\frac{l_n^v}{\g},\qquad
u_{(n)} \de\frac{u_n}{2n},
\\
C_{(n)}^\mathfrak{e}&\de&\biggl(\frac{C_n^\mathfrak{e}(2nt)} {\sqrt
{2n}} \biggr)_{0\le t \le
m_{(n)}^\mathfrak{e}},\qquad
L_{(n)}^\mathfrak{e}\de\biggl(\frac{L_n^\mathfrak{e}(2nt)} {\g}\biggr)
_{0\le t \le m_{(n)}^\mathfrak{e}},
\\
\mathfrak{M}_{(n)}^\mathfrak{e}&\de&\biggl(\frac{\mathfrak
{M}_n^\mathfrak{e}(\sqrt{2n} t)}{\g} \biggr)_{0\le t
\le\sigma_{(n)}^\mathfrak{e}}.
\end{eqnarray*}
\begin{rem*}
Throughout this paper, the notation with a parenthesized $n$ will
always refer to suitably rescaled objects, as in the definitions above.
\end{rem*}

We described in~\cite{bettinelli10slr} the limiting law of these objects:
\begin{prop}\label{cvint}
The random vector
\begin{eqnarray*}
&&\bigl( \s_n, \bigl( m_{(n)}^\mathfrak{e}\bigr)_{\mathfrak{e}\in\vec
E(\s_n)}, \bigl(\sigma
_{(n)}^\mathfrak{e}\bigr)_{\mathfrak{e}\in\vec E(\s_n)}, \bigl(
l_{(n)}^v\bigr)_{v\in V(\s
_n)}, u_{(n)}, \\
&&\hspace*{75.7pt}\bigl( C_{(n)}^\mathfrak{e}, L_{(n)}^\mathfrak{e}\bigr)_{\mathfrak
{e}\in\vec E(\s_n)}, \bigl( \mathfrak{M}
_{(n)}^\mathfrak{e}\bigr)_{\mathfrak{e}\in\vec E(\s_n)} \bigr)
\end{eqnarray*}
converges in law toward the random vector
\begin{eqnarray*}
&&\bigl( \s_\infty, (m_\infty^\mathfrak{e})_{\mathfrak{e}\in
\vec E(\s_\infty)}, (\sigma
_\infty^\mathfrak{e})_{\mathfrak{e}\in\vec E(\s_\infty)}, (
l_\infty^v)_{v\in V(\s
_\infty)}, u_\infty,\\
&&\hspace*{80.3pt}(C_\infty^\mathfrak{e}, L_\infty^\mathfrak{e})_{\mathfrak
{e}\in\vec E(\s_\infty)}, (\mathfrak{M}
_\infty^\mathfrak{e})_{\mathfrak{e}\in\vec E(\s_\infty)} \bigr),
\end{eqnarray*}
whose law is defined as follows:
\begin{itemize}[$\diamond$]
\item[$\diamond$] the law of the vector
\[
\mathfrak I_\infty\de\bigl(\s_\infty, (m_\infty^\mathfrak{e})
_{\mathfrak{e}\in\vec
E(\s_\infty)}, (\sigma_\infty^\mathfrak{e})_{\mathfrak
{e}\in\vec E(\s_\infty)}, (
l_\infty^v)_{v\in V(\s_\infty)}, u_\infty\bigr)
\]
is the probability $\mu$,
\item[$\diamond$] conditionally given $\mathfrak I_\infty$:
\begin{itemize}
\item the processes $(C_\infty^\mathfrak{e}, L_\infty^\mathfrak
{e})$, ${\mathfrak{e}\in\vec
E(\s_\infty)}$, and $(\mathfrak{M}_\infty^\mathfrak{e})$,
${\mathfrak{e}\in\ori E(\s_\infty)}$,
are independent;
\item the process $(C_\infty^\mathfrak{e}, L_\infty^\mathfrak
{e})$ has the law of a
Brownian snake's head on $[0,m_\infty^\mathfrak{e}]$ going from
$\sigma_\infty^\mathfrak{e}
$ to $0$
\[
(C_\infty^\mathfrak{e}, L_\infty^\mathfrak{e})\law\bigl(
F_{[0,m_\infty^\mathfrak{e}]}^{\sigma
_\infty^\mathfrak{e}\to0}, Z_{[0,m_\infty^\mathfrak{e}]}\bigr);
\]
\item the process $(\mathfrak{M}_\infty^\mathfrak{e})$ has the
law of a Brownian
bridge on $[0,\sigma_\infty^\mathfrak{e}]$ from $0$ to $l_\infty
^\mathfrak{e}\de l_\infty
^{\mathfrak{e}^+} - l_\infty^{\mathfrak{e}^-}$
\[
(\mathfrak{M}_\infty^\mathfrak{e})\law B_{[0,\sigma_\infty
^\mathfrak{e}]}^{0 \to l_\infty^\mathfrak{e}};
\]
\item the Motzkin paths are linked through the relation
\[
\mathfrak{M}_\infty^{\bar\mathfrak{e}}(s) = \mathfrak{M}_\infty
^\mathfrak{e}(\sigma_\infty^\mathfrak{e}-s) -
l_\infty^\mathfrak{e}.
\]
\end{itemize}
\end{itemize}
\end{prop}

Applying Skorokhod's representation theorem, we may and will assume
that this convergence holds almost surely. As a result, note that
for $n$ large enough, $\s_n=\s_\infty$.


\subsection{Maps seen as quotients of $[0,1]$}\label{secset}

Let $\mathfrak{q}_n$ be uniformly distributed over the set $\Q_n$ of bipartite
quadrangulations of genus $g$ with $n$ faces. Conditionally given~$\mathfrak{q}
_n$, we take $v_n^\bullet$ uniformly over $V(\mathfrak{q}_n)$ so that
$(\mathfrak{q}
_n,v_n^\bullet)$ is uniform over the set $\Qb_n$ of pointed bipartite
quadrangulations of genus $g$ with $n$ faces. Recall that every element
of $\Q_n$ has the same number of vertices, $n+2-2g$. Through the
Chapuy--Marcus--Schaeffer bijection, $(\mathfrak{q}_n,v_n^\bullet)$
corresponds
to a uniform well-labeled $g$-tree with $n$ edges $(\tr_n,\mathfrak
{l}_n)$. The
parameter $\eps_\pm\in\{-1,1\}$ appearing in the bijection will be
irrelevant to what follows.

Recall the notation $\tr_n(0)$, $\tr_n(1), \ldots, \tr_n(2n)$ and
$\mathfrak{q}
_n(0)$, $\mathfrak{q}_n(1), \ldots, \mathfrak{q}_n(2n)$ from
Section~\ref{seccms}. For
technical reasons, it will be more convenient, when traveling along the
$g$-tree, not to begin by its root but rather by the first edge of the
first forest. Precisely, we define
\[
\rr\tr_n(i) \de
\cases{
\tr_n(i - u_n +2n), &\quad if $0 \le i \le u_n$,\cr
\tr_n(i - u_n), &\quad if $u_n \le i \le2n$,}
\]
where $u_n$ is the integer recording the position of the root in the
first forest of $\tr_n$. We define $\rr\mathfrak{q}_n(i)$ in a
similar way, and
endow $\lbracket0,2n \rbracket$ with the pseudo-metric $d_n$ defined by
\[
d_n(i,j) \de d_{\mathfrak{q}_n}(\rr\mathfrak{q}_n(i),\rr
\mathfrak{q}_n(j)).
\]

We define the equivalence relation $\sim_n$ on $\lbracket0,2n
\rbracket$ by
declaring that $i \sim_n j$ if $\rr\mathfrak{q}_n(i)=\rr\mathfrak
{q}_n(j)$, that is, if
$d_n(i,j) =0$. We call $\pi_n$ the canonical projection from
$\lbracket0,2n \rbracket$ to $\lbracket0,2n \rbracket_{/\sim_n}$,
and we slightly abuse notation by
seeing $d_n$ as a~metric on $\lbracket0,2n \rbracket_{/\sim_n}$
defined by
$d_n(\pi
_n(i),\pi_n(j)) \de d_n(i,j)$. In what follows, we will always make the
same abuse with every pseudo-metric. The metric space $(\lbracket0,2n
\rbracket_{/\sim_n},d_n )$ is then isometric to $(V(\mathfrak
{q}_n)\bs\{
v_n^\bullet\},d_{\mathfrak{q}_n} )$, which is at $\dGH$-distance
at most $1$
from the space $(V(\mathfrak{q}_n),d_{\mathfrak{q}_n} )$.\vadjust{\goodbreak}

We extend the definition of $d_n$ to noninteger values by linear
interpolation and define its rescaled version: for $s,t \in[0,1]$, we let
%
%
\begin{equation}\label{rescdn}
d_{(n)}(s,t) \de\frac1 {\g} d_n(2ns,2nt).
\end{equation}


\subsubsection*{Spatial contour function of $(\tr_n,\mathfrak{l}_n)$}

The spatial contour function of the pair $(\tr_n,\mathfrak{l}_n)$ is the
function $\Lab_n\dvtx [0,2n] \to\R$, defined by
\[
\Lab_n(i) \de\mathfrak{l}_n(\rr\tr_n(i))- \mathfrak{l}_n(
\rr\tr_n(0)),\qquad 0\le i \le2n,
\]
and linearly interpolated between integer values. Its rescaled version is
\[
\Lab_{(n)} \de\biggl(\frac{\Lab_n(2nt)} {\g}\biggr)_{0\le t \le1}.
\]

Recall definition~(\ref{lele}) of the process $\Lab_n^\mathfrak{e}$.
We define
its rescaled version by
\[
\Lab_{(n)}^\mathfrak{e}\de\biggl(\frac{\Lab_n^\mathfrak{e}(2nt)} {\g
}\biggr)_{0\le t \le
m_{(n)}^\mathfrak{e}} = \bigl(L_{(n)}^\mathfrak{e}(t) + \mathfrak
{M}_{(n)}^\mathfrak{e}\bigl( \sigma_{(n)}^\mathfrak{e}-
\underline C_{(n)}^\mathfrak{e}(t) \bigr) \bigr)_{0\le t \le
m_{(n)}^\mathfrak{e}}.
\]
Proposition~\ref{cvint} shows that $\Lab_{(n)}^\mathfrak{e}$
converges in the
space $(\K, d_\K)$ toward
\[
\Lab_\infty^\mathfrak{e}\de\bigl( L_\infty^\mathfrak{e}(t) +
\mathfrak{M}_\infty^\mathfrak{e}\bigl( \sigma
_\infty^\mathfrak{e}- \underline C_\infty^\mathfrak{e}(t) \bigr)
\bigr)_{0\le t \le m_\infty
^\mathfrak{e}}.
\]

We can express $\Lab_{(n)}$ in terms of the processes $\Lab
_{(n)}^\mathfrak{e}$'s
by concatenating them. For $f,g\in\K_0$ two functions started at $0$,
we call $f \bullet g \in\K_0$ their concatenation defined by $\sigma(f
\bullet g) \de\sigma(f) + \sigma(g)$ and, for $0\le t \le\sigma(f
\bullet g)$,
\[
f \bullet g(t) \de
\cases{
f(t), &\quad if $0\le t \le\sigma(f)$,\cr
f(\sigma(f)) + g\bigl(t-\sigma(f)\bigr), &\quad if $\sigma(f) \le t \le
\sigma(f) + \sigma(g)$.}
\]
We arrange the half-edges of $\s_n$ according to its facial order,
beginning with the root $\mathfrak{e}_1=\mathfrak{e}_*, \ldots,
\mathfrak{e}_\ks$, so that $\Lab
_{(n)} = \Lab_{(n)}^{\mathfrak{e}_1} \bullet\Lab_{(n)}^{\mathfrak
{e}_2} \bullet\cdots
\bullet\Lab_{(n)}^{\mathfrak{e}_\ks}$. By continuity of the
concatenation, $\Lab
_{(n)}$ converges in $(\K,d_\K)$ toward $\Lab_\infty\de\Lab
_\infty^{\mathfrak{e}
_1} \bullet\Lab_\infty^{\mathfrak{e}_2} \bullet\cdots\bullet\Lab
_\infty^{\mathfrak{e}_\ks
}$, where the half-edges of $\s_\infty$ are arranged in the same way.

\subsubsection*{Upper bound for $d_{(n)}$}

Bound~(\ref{dlemme}) provides us with an upper bound on~$d_{(n)}$. We define
\[
d_n^\circ(i,j) \de\Lab_n(i) + \Lab_n(j) - 2 \max\Bigl(\min_
{\mbox{\fontsize{8.36pt}{8.36pt}\selectfont{$k \mbox{$\in$}
\overrightarrow{\lbracket i,j \rbracket}$}}} \Lab_n(k) ,\min_
{\mbox{\fontsize{8.36pt}{8.36pt}\selectfont{$k \mbox{$\in$}
\overrightarrow
{\lbracket j,i \rbracket}$}}} \Lab_n(k) \Bigr)+2,
\]
we extend it to $[0,2n]$ by linear interpolation and define its
rescaled version~$d_{(n)}^\circ$ as we did for $d_n$ by~(\ref{rescdn}).
We readily obtain that
%
%
\begin{equation}\label{bound2}
d_{(n)}(s,t) \le d^\circ_{(n)}(s,t).
\end{equation}

Moreover, the process $(d^\circ_{(n)}(s,t))_{0 \le s,t \le1}$
converges in $(\C([0,1]^2,\R),\|\cdot\|_\infty)$ toward the
process $(d^\circ_\infty(s,t))_{0 \le s,t \le1}$ defined by
\[
d^\circ_\infty(s,t) \de\Lab_\infty(s) + \Lab_\infty(t) - 2 \max
\Bigl(
\min_{\mbox{\fontsize{8.36pt}{8.36pt}\selectfont{$x \mbox{$\in$}
\overrightarrow{[s,t]}$}}} \Lab_\infty(x),\min_
{\mbox{\fontsize{8.36pt}{8.36pt}\selectfont{$x \mbox{$\in$}
\overrightarrow{[t,s]}$}}} \Lab_\infty(x) \Bigr),
\]
where
%
%
\begin{equation}\label{orast}
\overrightarrow{[s,t]} \de
\cases{
[s,t], &\quad if $s \le t$, \cr
[s,1] \cup[0,t], &\quad if $t < s$.}
\end{equation}


\subsubsection*{Tightness of the processes $d_{(n)}$'s}

In~\cite{bettinelli10slr}, Lemma 19, we showed the tightness of the
processes $d_{(n)}$'s laws thanks to the inequality~(\ref{bound2}). As a
result, from any increasing sequence of integers, we may extract a
(deterministic) subsequence $(n_k)_{k\ge0}$ such that there exists a
function $d_\infty\in\C([0,1]^2,\R)$ satisfying
%
%
\begin{equation}\label{dinfty}
(d_{(n_k)}(s,t) )_{0 \le s,t \le1} \xrightarrow_{k \to\infty}^{(\mathrm d)} (
d_\infty(s,t) )_{0 \le s,t \le1}.
\end{equation}
By Skorokhod's representation theorem, we will assume that this
convergence holds almost surely. We can check that the function
$d_\infty$ is actually a pseudo-metric. We define the equivalence
relation associated with it by saying that $s \sim_\infty t$ if
$d_\infty(s,t)=0$, and we call $\mathfrak{q}_\infty\de[0,1]_{/\sim
_\infty}$. We
proved in~\cite{bettinelli10slr} that
\[
\biggl( V(\mathfrak{q}_{n_k}),\frac1 {\gamma n_k^{1/4}} d_{\mathfrak
{q}_{n_k}} \biggr) \xrightarrow_{k \to\infty}^{(\mathrm d)}
(\mathfrak{q}_\infty,d_\infty)
\]
in the sense of the Gromov--Hausdorff topology.

From now on, we fix such a subsequence $(n_k)_{k\ge0}$. We will always
focus on this particular subsequence in the following, and we will
consider convergences when $n\to\infty$ to hold along this particular
subsequence.


\section{Real $g$-trees}\label{secrgt}

In the discrete setting, it is sometimes convenient to work directly
with the space $\tr_n$ instead of $\lbracket0,2n \rbracket$. In the
continuous
setting, we will see $\mathfrak{q}_\infty$ as a quotient of a
continuous version
of a $g$-tree, which we will call real $g$-tree. In other words, we
will see the identifications $s\sim_\infty t$ as of two different
kinds: some are inherited ``from the $g$-tree structure,'' whereas the
others come ``from the map structure.''

\subsection{Definitions}\label{secdefgtree}

As $g$-trees generalize plane trees in genus $g$, real $g$-trees are
the objects that naturally generalize real trees. We will only use
basic facts on real trees in this work. See, for example,
\cite{legall05rta} for more detail.

We consider a fixed dominant scheme $\s\in\Sg^*$. Let $(m^\mathfrak
{e})_{\mathfrak{e}\in
\vec E(\s)}$ and $(\sigma^\mathfrak{e})_{\mathfrak{e}\in\vec E(\s
)}$ be two families of
positive numbers satisfying $\sum_\mathfrak{e}m^\mathfrak{e}=1$ and
$\sigma^\mathfrak{e}=\sigma
^{\bar\mathfrak{e}}$ for all~$\mathfrak{e}$. As usual, we arrange
the half-edges of $\s$
according to its facial order, $\mathfrak{e}_1=\mathfrak{e}_*, \ldots
, \mathfrak{e}_\ks$. For every
$s\in[0,1)$, there exists a unique $1 \le k \le\ks$ such that
\[
\sum_{i=1}^{k-1} m^{\mathfrak{e}_i} \le s < \sum_{i=1}^{k}
m^{\mathfrak{e}_i}.
\]
We let $\mathfrak{e}(s) \de\mathfrak{e}_k$ and $\langle s \rangle
\de s- \sum_{i=1}^{k-1}
m^{\mathfrak{e}_i} \in[0,m^{\mathfrak{e}(s)})$. By convention, we
set $\mathfrak{e}(1)=\mathfrak{e}_1$ and
$\langle1 \rangle=0$. Beware that these notions depend on the family
$(m^\mathfrak{e})_{\mathfrak{e}\in\vec E(\s)}$. There should be no
ambiguity in what follows.

Let us suppose we have a family $(h^\mathfrak{e})_{\mathfrak{e}\in
\vec E(\s)}$ of
continuous functions $h^\mathfrak{e}\dvtx\break [0, m^\mathfrak{e}] \to\R_+$
such that $h^\mathfrak{e}
(0)=\sigma^\mathfrak{e}$ and $h^\mathfrak{e}(m^\mathfrak{e})=0$. It
will be useful to consider their
concatenation: we define the continuous function $h\dvtx [0,1] \to\R_+$
going from~$\sum_\mathfrak{e}\sigma^\mathfrak{e}$ to $0$ by
%
%
\begin{equation}\label{defconcat}
h \de( h^{\mathfrak{e}_1} - \sigma^{\mathfrak{e}_1} )
\bullet( h^{\mathfrak{e}_2} -
\sigma^{\mathfrak{e}_2} ) \bullet\cdots\bullet\bigl(
h^{\mathfrak{e}_\ks} - \sigma^{\mathfrak{e}
_\ks} \bigr) + \sum_{i=1}^{\ks} \sigma^{\mathfrak{e}_i}.\hspace*{-26pt}
\end{equation}

We define the relation $\simeq$ on $[0,1]$ as the coarsest equivalence
relation for which $s \simeq t$ if one of the following occurs:
%
%
\begin{subequation}
\begin{eqnarray}
\label{eqit21}
&&h(s)=h(t)=\inf_{[s\wedge t, s\vee t]} h;
\\
\label{eqit22}
&&h(s) = \underline h (s),\qquad h(t) = \underline
h (t),\nonumber\\[-8pt]\\[-8pt]
&&\qquad\mathfrak{e}(s)=\overline{\mathfrak{e}(t)}
\quad\mbox{and}\quad
h^{\mathfrak{e}(s)}(\langle s \rangle) =
\sigma^{\mathfrak{e}(t)} - h^{\mathfrak{e}(t)}(\langle
t\rangle);\nonumber
\\
\label{eqit23}
&&\langle s\rangle = \langle t\rangle= 0\quad\mbox{and}\quad
\mathfrak{e}(s)^- = \mathfrak{e}(t)^-.
\end{eqnarray}
\end{subequation}

If we see the $h^\mathfrak{e}$'s as contour functions (in a continuous setting),
the first item identifies numbers coding the same point in one of the
forests. The second item identifies the floors of forests ``facing each
other'': the numbers~$s$ and $t$ should code floor points (two first
equalities) of forests facing each other (third equality) and
correspond to the same point (fourth equality). Finally, the third item
identifies the nodes. We call \textit{real $g$-tree} any space $\TT
\de
[0,1]_{/\simeq}$ obtained by such a construction.\footnote{There should
be a more intrinsic definition for these spaces in terms of compact
metric spaces that are locally real trees. As we will need to use this
construction in what follows, we chose to define them as such for simplicity.}


%
%
\begin{figure}

\includegraphics{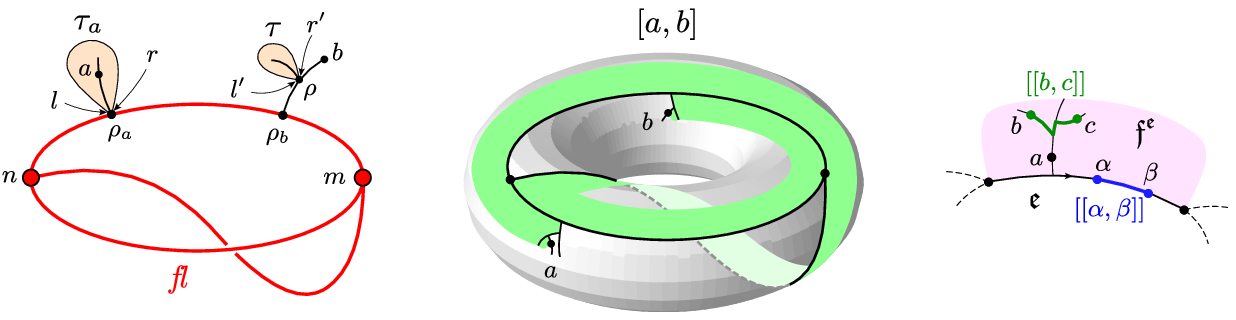}

\caption{\textup{Left.} On this picture, we can see the floor $\fl$,
the two nodes $n$ and $m$, an example of tree $\tau_a$ and an example
of tree $\tau$ to the left of $\lhb\rho_b,b \rhb$ rooted at $\rho$.
\textup{Middle.} The set $[a,b]$. \textup{Right.} On this picture,~$a$
is an ancestor of~$b$ and $c$, and we can see the sets $\lhb b,c \rhb$,
$\lhb\alpha,\beta\rhb$ and $\mathfrak{f}^\mathfrak{e}$.}
\label{not}
\end{figure}

We now define the notions we will use throughout this work (see
Figure~\ref{not}). For $s \in[0,1]$, we write $\TT(s)$ its equivalence
class in the quotient $\TT= [0,1]_{/\simeq}$. Similarly to the
discrete case, the floor of $\TT$ is defined as follows.
\begin{defi}
We call \textit{floor} of $\TT$ the set $\fl\de\TT(\{s \dvtx
h(s) = \underline h(s) \})$.
\end{defi}

For $a=\TT(s)\in\TT\bs\fl$, let $l \de\inf\{t \le s \dvtx
\underline
h(t)=\underline{h}(s) \}$ and $r \de\sup\{t \ge s \dvtx\break \underline
h(t)=\underline{h}(s) \}$. The set $\tau_a \de\TT([l,r])$ is a real
tree rooted at $\rho_a \de\TT(l)=\TT(r)\in\fl$.
\begin{defi}
We call \textit{tree} of $\TT$ a set of the form $\tau_a$ for any
$a\in
\TT\bs\fl$.
\end{defi}

If $a\in\fl$, we simply set $\rho_a \de a$. Let $\tau$ be a tree of
$\TT
$ rooted at $\rho$, and $a,b \in\tau$. We call $\lhb a,b \rhb$ the
range of the unique injective path linking~$a$ to~$b$. In particular,
the set $\lhb\rho,a\rhb$ will be of interest. It represents the
ancestral lineage of~$a$ in the tree~$\tau$. We say that~$a$ is an
\textit{ancestor} of~$b$, and we write $a \preceq b$, if $a \in\lhb
\rho,b \rhb$. We write $a \prec b$ if $a \preceq b$ and $a \neq b$.
%
\begin{defi}
Let $b=\TT(t) \in\TT\bs\fl$ and $\rho\in\lhb\rho_b,b\rhb\bs\{
\rho_b,b\}$. Let $l' \de\inf\{s\le t \dvtx \TT(s)=\rho\}$ and $r'
\de
\sup\{s\le t \dvtx \TT(s)=\rho\}$. Then, provided $l' \neq r'$, we
call the \textit{tree to the left} of $\lhb\rho_b,b\rhb$ rooted at
$\rho$ the set $\TT([l',r'])$.

We define the \textit{tree to the right} of $\lhb\rho_b,b\rhb$ rooted
at $\rho$ in a similar way, by replacing ``$\le$'' with ``$\ge$'' in
the definitions of $l'$ and $r'$.
\end{defi}
\begin{defi}
We call \textit{subtree} of $\TT$ any tree of $\TT$, or any tree to
the left or right of $\lhb\rho_b,b\rhb$ for some $b \in\TT\bs\fl$.
\end{defi}

Note that subtrees of $\TT$ are real trees, and that trees of $\TT$ are
also subtrees of~$\TT$. For a subtree $\tau$, the maximal interval
$[s,t]$ such that $\tau= \TT([s,t])$ is called the \textit{interval
coding} the subtree $\tau$.
\begin{defi}
For $\mathfrak{e}\in\vec E(\s)$, we call the \textit{forest to the
left of $\mathfrak{e}$}
the set $\mathfrak{f}^\mathfrak{e}\de\TT( \overline{ \{
s \dvtx \mathfrak{e}(s)=\mathfrak{e}\}
} )$.
\end{defi}

The \textit{nodes} of $\TT$ are the elements of $\TT(\{s \dvtx
\langle s
\rangle= 0\})$. In what follows, we will identify the nodes of $\TT$
with the vertices of $\s$. In particular, the two nodes $\mathfrak
{e}^-$ and $\mathfrak{e}
^+$ lie in $\mathfrak{f}^\mathfrak{e}$. We extend the definition of
$\lhb a,b \rhb$ to the
floor of $\mathfrak{f}^\mathfrak{e}$: for $a,b \in\mathfrak
{f}^\mathfrak{e}\cap\fl$, let $s, t \in\overline
{ \{r \dvtx \mathfrak{e}(r)=\mathfrak{e}\}} $\vadjust{\goodbreak} be such that
$a=\TT(s)$ and $b= \TT
(t)$. We define
\[
\lhb a, b \rhb\de\TT([s\wedge t,s\vee t]) \cap\fl
\]
the range of the unique\footnote{Note that $\mathfrak{e}^+ \neq
\mathfrak{e}^-$ because $\s
$ is a dominant scheme.} injective path from~$a$ to~$b$ that stays
inside $\mathfrak{f}^\mathfrak{e}$. For clarity, we write the set
$\lhb\mathfrak{e}^-,\mathfrak{e}^+\rhb$
simply as $\lhb\mathfrak{e}\rhb$. Note that, in particular, $\lhb
\mathfrak{e}\rhb= \mathfrak{f}
^\mathfrak{e}\cap\mathfrak{f}^{\bar\mathfrak{e}}=\mathfrak
{f}^\mathfrak{e}\cap\fl$.

Let $a,b\in\TT$. There is a natural way\footnote{Note that, if $a,b
\in\fl$, there are other possible ways to explore the $g$-tree between
them. Indeed, a point of $\fl$ is visited twice---or three times if it
is a node---when we travel around~$\fl$. In particular, this definition
depends on the position of the root in $\s$ for such points. In what
follows, we never use this definition for such points, so there will be
no confusion.} to explore $\TT$ from~$a$ to~$b$. If $\inf\TT^{-1}(a)
\le\sup\TT^{-1}(b)$, then let $t\de\inf\{ r \ge\inf\TT^{-1}(a)
\dvtx b=\TT(r) \}$ and $s\de\sup\{ r \le t \dvtx a=\TT(r) \}$. If
$\sup
\TT^{-1}(b) < \inf\TT^{-1}(a)$, then let $t \de\inf\TT^{-1} (b)$ and
$s\de\sup\TT^{-1}(a)$. We define
%
%
\begin{equation}\label{ab}
{[a,b]} \de\TT(\overrightarrow{[s,t]} ),
\end{equation}
where $\overrightarrow{[s,t]}$ is defined by~(\ref{orast}).


We call $\TT_n$ (resp., $\TT_\infty$) the real $g$-tree obtained from
the scheme $\s_n$ (resp.,~$\s_\infty$) and the family
$(C_{(n)}^\mathfrak{e})_{\mathfrak{e}
\in\vec E(\s_n)}$ [resp., $(C_\infty^\mathfrak{e})_{\mathfrak{e}\in
\vec E(\s_\infty)}$].
For the sake\vspace*{1pt} of consistency with~\cite{bettinelli10slr}, we call $\CC
_{(n)}$ and $\CC_\infty$ the functions obtained by~(\ref{defconcat}) in
this construction. We also call $\simeq_{(n)}$ and $\simeq_\infty$ the
corresponding equivalence relations. When dealing with $\TT_\infty$, we
add an $\infty$ symbol to the notation defined above: for example, the
floor of $\TT_\infty$ will be noted $\fl_\infty$, and its forest to the
left of $\mathfrak{e}$ will be noted $\mathfrak{f}_\infty^\mathfrak
{e}$. It is more natural to
use $\tr_n$ rather than $\TT_n$ in the discrete setting. As $\tr_n$ may
be viewed as a subset of $\TT_n$, we will use for $\tr_n$ the formalism
we defined above simply by restriction. Note that the notions of floor,
forests, trees and nodes are consistent with the definitions we gave in
Section~\ref{decomp} in that case.


Note that, because the functions $C^\mathfrak{e}_\infty$'s are first-passage
Brownian bridges, the probability that there exists $\eps>0$ such that
$C^\mathfrak{e}_\infty(s) > C^\mathfrak{e}_\infty(0)$ for all $s\in
(0,\eps)$ is equal
to $0$. As a result, there are almost surely no\vspace*{1pt} trees rooted at the
nodes of $\TT_\infty$. Moreover, the fact that the forests $\mathfrak
{f}^\mathfrak{e}$
and $\mathfrak{f}^{\bar\mathfrak{e}}$ are independent yields that,
almost surely, we cannot
have a tree in $\mathfrak{f}^\mathfrak{e}$ and a tree in $\mathfrak
{f}^{\bar\mathfrak{e}}$ rooted at the same
point. As a consequence, we see that, almost surely, all the points of
$\TT_\infty$ are of order less than~$3$.

\subsection{Maps seen as quotients of real $g$-trees}

Consistently with the notation $\tr_n(i)$ and $\mathfrak{q}_n(i)$ in
the discrete
setting, we call $\TT_\infty(s)$ [resp., $\mathfrak{q}_\infty(s)$]
the equivalence
class of $s\in[0,1]$ in $\TT_\infty= [0,1]_{/\simeq_\infty}$ (resp.,
in $\mathfrak{q}_\infty= [0,1]_{/\sim_\infty}$).
\begin{lem}\label{coarse}
The equivalence relation $\simeq_\infty$ is coarser than $\sim
_\infty$,
so that we can see $\mathfrak{q}_\infty$ as the quotient of $\TT
_\infty$ by the
equivalence relation on $\TT_\infty$ induced from $\sim_\infty$.
\end{lem}
\begin{pf}
By definition of $\simeq_\infty$, it suffices to show that if $s < t$
satisfy~(\ref{eqit21}), (\ref{eqit22}) or~(\ref{eqit23}), then
$s\sim
_\infty t$. Let us first suppose that $s$ and $t$ satisfy (\ref
{eqit21}), that is,
\[
\CC_\infty(s) = \CC_\infty(t) = \inf_{[s,t]} \CC_\infty.
\]
In a first time, we moreover suppose that $\CC_\infty(r) > \CC
_\infty
(s)$ for all $r\in(s,t)$. Using Proposition~\ref{cvint}, we can find
integers $0 \le s_n < t_n \le2n$ such that $(s_{(n)},t_{(n)}) \de
(s_n/2n, t_n/2n) \to(s,t)$ and $\CC_{(n)}(s_{(n)}) = \CC
_{(n)}(t_{(n)}) = \inf_{[s_{(n)},t_{(n)}]} \CC_{(n)}$. The latter
condition imposes that $\rr\tr_n(s_n)=\rr\tr_n(t_n)$ so that
$d_n(s_n,t_n)=0$ and $s \sim_\infty t$ by~(\ref{dinfty}).

Equation~(\ref{f}) shows that, for every $\mathfrak{e}$, the law of
$C^\mathfrak{e}_\infty$
is absolutely continuous with respect to the Wiener measure on any
interval $[0, m^\mathfrak{e}_\infty- \eps]$, for $\eps>0$. Because local
minimums of Brownian motion are pairwise distinct, this is also true
for any $C^\mathfrak{e}_\infty$, and thus for the whole process $\CC
_\infty$ by
construction. If there exists $r\in(s,t)$ for which $\CC_\infty(r) =
\CC_\infty(s)$, it is thus unique. We may then apply the previous
reasoning to $(s,r)$ and $(r,t)$ and find that $s\sim_\infty r$ and
$r\sim_\infty t$, so that $s\sim_\infty t$.

Let us now suppose that $s$ and $t$ satisfy~(\ref{eqit22}). If there is
$0 \le r < s$ such that $\CC_\infty(r)=\CC_\infty(s)$, then $r
\simeq
_\infty s$ by~(\ref{eqit21}). The same holds with $t$ instead of~$s$.
We may thus restrict our attention to $s$ and $t$ for which $\CC
_\infty
(r) > \CC_\infty(s)$ for all $r\in[0,s)$ and $\CC_\infty(r) > \CC
_\infty(t)$ for all $r\in[0,t)$. Let us call $\mathfrak{e}=\mathfrak
{e}(s) = \overline{ \mathfrak{e}
(t)}$. In order to avoid confusion, we use the notation $\langle\cdot
\rangle_n$ and $\mathfrak{e}_n(\cdot)$ when dealing with the
functions $\CC
_{(n)}^\mathfrak{e}$'s. We know that for $n$ large enough, we have $\s
_n=\s_\infty
$. We only consider such $n$'s in the following. We first find $0 \le
s_n \le2n$ such that $s_{(n)} \de s_n/2n \to s$, $\mathfrak
{e}_n(s_{(n)})=\mathfrak{e}$,
and $\CC_{(n)}(s_{(n)})=\underline\CC_{(n)}(s_{(n)})$. We define
\[
t_{(n)} \de\inf\biggl\{r \in\frac1 {2n} \lbracket0,2n \rbracket \dvtx
\mathfrak{e}
_n(r)=\bar\mathfrak{e}, \CC_{(n)}^{\bar\mathfrak{e}}(\langle
r \rangle_n) = \sigma
_{(n)}^{\mathfrak{e}} - \CC_{(n)}^{\mathfrak{e}}\bigl(\bigl\langle
s_{(n)}\bigr\rangle_n\bigr) \biggr\},
\]
so that $t_{(n)} \simeq_{(n)} s_{(n)}$, and then $d_{(n)}
(s_{(n)},t_{(n)}) =0$. Taking an extraction if needed, we may
suppose that $t_{(n)} \to t'\sim_\infty s$. By construction,
$\mathfrak{e}(t')=\mathfrak{e}
(t)$ and $\CC_\infty(t')=\underline\CC_\infty(t')=\CC_\infty(t)$.
So $t'$ and $t$ fulfill requirement~(\ref{eqit21}) and $t'\sim_\infty
t$ by the above argument. The case of~(\ref{eqit23}) is easier and may
be treated in a~similar way.
\end{pf}

This lemma allows us to define a pseudo-metric and an equivalence
relation on $\TT_\infty$, still denoted by $d_\infty$ and $\sim
_\infty
$, by setting $d_\infty(\TT_\infty(s),\TT_\infty(t)) \de
d_\infty(s,t)$ and declaring $\TT_\infty(s) \sim_\infty\TT_\infty(t)$
if $s \sim_\infty t$. The metric space $(\mathfrak{q}_\infty
,d_\infty)$ is then
isometric to $({\TT_\infty}_{/\sim_\infty},d_\infty)$. We
define $d^\circ_\infty$ on $\TT_\infty$ by letting
\[
d^\circ_\infty(a,b) \de\inf\{d^\circ_\infty(s,t) \dvtx a=\TT
_\infty
(s), b=\TT_\infty(t) \}.\vadjust{\goodbreak}
\]

We will see in Lemma~\ref{min} that there is a.s. only one point where
the function~$\Lab_\infty$ reaches its minimum. On this event, the
following lemma holds.\vspace*{-1pt}

\begin{lem}\label{dssb}
$\!\!\!$Let $s^\bullet$ be the unique point where $\Lab_\infty$ reaches its
minimum.~Then
\[
d_\infty(s,s^\bullet) = \Lab_\infty(s) - \Lab_\infty(s^\bullet).
\]

Moreover, $s\sim_\infty t$ implies $\Lab_\infty(s)=\Lab_\infty(t)$.\vspace*{-1pt}
\end{lem}
\begin{pf}
This readily comes from the discrete setting. Let $0 \le s_n^\bullet
\le2n$ be an integer where $\Lab_n$ reaches its minimum. By extracting
if necessary, we may suppose that $s_n^\bullet/2n$ converges
and its limit is
necessarily $s^\bullet$. Let $0 \le s_n \le2n$ be such that
\mbox{$s_n/2n \to s$}. From the Chapuy--Marcus--Schaeffer bijection,
$d_n(s_n,s_n^\bullet)= \Lab_n(s_n)- \Lab_n(s_n^\bullet)+ 1$.
Letting $n
\to\infty$ after renormalizing yields the first assertion. The second
one follows from the first one and the triangle inequality.\vspace*{-1pt}
\end{pf}

As a result of Lemmas~\ref{coarse} and~\ref{dssb}, we can define
$\Lab
_\infty$ on $\TT_\infty$ by $\Lab_\infty(\TT_\infty(s)
) \de\Lab
_\infty(s)$. When $(a,b) \notin(\fl_\infty)^2$, we have
%
%
\begin{equation}\label{dzero}\quad
d^\circ_\infty(a,b) =\Lab_\infty(a) + \Lab_\infty(b) - 2 \max
\Bigl(
\min_{ x \in[a,b]} \Lab_\infty(x),\min_{ x \in[b,a]} \Lab_\infty(x)
\Bigr),
\end{equation}
where $[a,b]$ was defined by~(\ref{ab}).\vspace*{-1pt}


\section{Points identifications}\label{secip}

This section is dedicated to the proof of the following theorem:\vspace*{-1pt}
\begin{theorem}\label{ip}
Almost surely, for every $a,b \in\TT_\infty$, $a\sim_\infty b$ is
equivalent to $d^\circ_\infty(a,b)=0$.\vspace*{-1pt}
\end{theorem}

We already know that $d^\circ_\infty(a,b)=0$ implies $a\sim_\infty b$
from the bound $d_\infty\le d^\circ_\infty$. We will show the converse
through a series of lemmas. We adapt the approach of Le Gall
\cite{legall07tss} to our setting.\vspace*{-1pt}

\subsection{Preliminary lemmas}\label{secprellem}

Let us begin by giving some information on the process $(\CC_\infty,
\Lab_\infty)$.\vspace*{-1pt}
\begin{lem}\label{min}
The set of points where $\Lab_\infty$ reaches its minimum is a.s. a~singleton.\vspace*{-1pt}
\end{lem}

Let $f\dvtx[0,1] \to\R$ be a continuous function. We say that $s\in[0,1)$
is a~\textit{right-increase point} of $f$ if there exists $t\in(s,1]$
such that $f(r) \ge f(s)$ for all $s\le r \le t$. A \textit
{left-increase point} is defined in a symmetric way. We call $\IP(f)$
the set of all (left or right) increase points of $f$.\vspace*{-1pt}
\begin{lem}\label{pc}
A.s., $\IP(\CC_\infty)$ and $\IP(\Lab_\infty)$ are disjoint sets.\vspace*{-1pt}
\end{lem}

As the proofs of these lemmas are rather technical and unrelated to
what follows, we postpone them to Section~\ref{secprle}.\vadjust{\goodbreak}


\subsection{Key lemma}

\begin{rem*}
In what follows, every discrete path denoted by the letter ``$\wp$''
will always be a path in the \textit{map}, never in the tree, that is, a
path using the edges of the map.
\end{rem*}

Let $\tau$ be a subtree of $\tr_n$ and $\wp=(\wp(0),\wp(1),\ldots
,\wp
(r))$ be a path in $\mathfrak{q}_n$ that avoids the base
point $v_n^\bullet$. We
say that the arc $(\wp(0),\wp(1))$ enters the subtree $\tau$ from the
left (resp., from the right) if $\wp(0) \notin\tau$, $\wp(1) \in
\tau$
and $\mathfrak{l}_n(\wp(1)) - \mathfrak{l}_n(\wp(0))=-1$ [resp.,
$\mathfrak{l}_n(\wp(1)) - \mathfrak{l}
_n(\wp(0))=1$]. We say that the path~$\wp$ \textit{passes through} the
subtree $\tau$ between times $i$ and $j$, where $0 < i \le j < r$, if:
\begin{itemize}[$\diamond$]
\item[$\diamond$]$\wp(i-1) \notin\tau$; $\wp(\lbracket i,j \rbracket)
\subseteq\tau$; $\wp(j+1)
\notin\tau$,
\item[$\diamond$]$\mathfrak{l}_n(\wp(i)) - \mathfrak{l}_n(\wp(i-1)) =
\mathfrak{l}_n(\wp(j+1)) - \mathfrak{l}_n(\wp(j))$.
\end{itemize}
The first condition states that $\wp$ ``visits'' $\tau$, whereas the
second one ensures that it really goes ``through.'' It enters and
exits $\tau$ going ``in the same direction.''

We say that a vertex $a_n\in\tr_n$ converges toward a point $a\in\TT
_\infty$ if there exists a sequence of integers $s_n \in\lbracket
0,2n \rbracket$
coding $a_n$ [i.e., $a_n=\rr\tr_n(s_n)$] such that $s_n/2n$ admits a
limit $s$ satisfying $a=\TT_\infty(s)$. Let $\lbracket l_n,r_n
\rbracket$ be the
intervals coding subtrees $\tau_n \subseteq\tr_n$. We say that the
subtree $\tau_n$ converges toward a~subtree $\tau\subseteq\TT
_\infty$
if the sequences $l_n/2n$ and $r_n/2n$ admit limits $l$ and $r$ such
that the interval coding $\tau$ is $[l,r]$. The following lemma is
adapted from Le Gall~\cite{legall07tss}, end of Proposition 4.2.
\begin{lem}\label{legalllem}
With full probability, the following occurs. Let $a, b \in\TT_\infty$
be such that $\Lab_\infty(a)=\Lab_\infty(b)$. We suppose that there
exists a subtree $\tau$ rooted at~$\rho$ such that $\inf_{\tau}
\Lab
_\infty< \Lab_\infty(a) < \Lab_\infty(\rho)$. We further suppose that
we can find vertices $a_n$, $b_n \in\tr_n$ and subtrees $\tau_n$
in $\tr_n$ converging, respectively, toward~$a$,~$b$, $\tau$ and
satisfying the following property: for infinitely many $n$'s, there
exists a geodesic path $\wp_n$ in $\mathfrak{q}_n$ from $a_n$
to $b_n$ that
avoids the base point~$v_n^\bullet$ and passes through the
subtree $\tau_n$.

Then, $a\not\sim_\infty b$.
\end{lem}
\begin{pf}
The idea is that if~$a$ and~$b$ were identified, then all the points in
the discrete subtrees close (in a certain sense) to the geodesic path
would be close to~$a$ in the limit. Fine estimates on the sizes of
balls yield the result. We proceed to the rigorous proof.

We reason by contradiction and suppose that $a \sim_\infty b$. We only
consider integers $n$ for which the hypothesis holds. We call $\rho_n$
the root of $\tau_n$, and we set, for $\eps>0$,
\[
\mathcal U_\infty^\eps\de\biggl\{y\in\tau \dvtx \Lab_\infty(y) <
\Lab
_\infty(a) + \eps; \forall x \in\lhb\rho,y\rhb, \Lab_\infty
(x) >
\Lab_\infty(a)+\frac\eps8 \biggr\}.
\]
We first show that $\mathcal U_\infty^\eps\subseteq B_{\infty
}(a,2\eps
)$, where $B_{\infty}(a,2\eps)$ denotes the closed ball of radius
$2\eps
$ centered at~$a$ in the metric space $(\mathfrak{q}_\infty,d_\infty
)$. Let $y\in
\mathcal U_\infty^\eps$. We can find $y_n\in\tau_n\bs\{\rho_n\}$
converging toward $y$. For $n$ large enough, we have
\begin{eqnarray*}
d_{\mathfrak{q}_n}(a_n,b_n) &\le&\frac\eps{32} n^{1/4} ,\qquad\sup
_{c\in\wp
_n}|\mathfrak{l}_n(c)-\mathfrak{l}_n(a_n)| \le\frac\eps{32}
n^{1/4},\\
\mathfrak{l}_n(y_n) &\le&\mathfrak{l}_n(a_n) + \frac32 \eps
n^{1/4}, \qquad\forall x
\in\lhb\rho_n,y_n\rhb\qquad \mathfrak{l}_n(x) \ge\mathfrak
{l}_n(a_n) + \frac\eps{16} n^{1/4}.
\end{eqnarray*}
The first inequality comes from the fact that $a\sim_\infty b$. The
second inequality is a consequence of the first one. The third
inequality holds because $(\mathfrak{l}_n(y_n) - \mathfrak
{l}_n(v_n^\bullet))/\g\to
\Lab_\infty(y)$ and $(\mathfrak{l}_n(a_n) - \mathfrak
{l}_n(v_n^\bullet))/\g\to\Lab
_\infty(a)$. Finally, the fourth inequality follows by compactness of
$\lhb\rho,y \rhb$.

From now on, we only consider such $n$'s. We call $t_n\de\sup\{t \dvtx
y_n=\rr\tr_n(t)\}$ the last integer coding $y_n$, and $\lbracket
l_n,r_n \rbracket$
the interval coding $\tau_n$. We also call $i \le j$ two integers such
that $\wp_n$ passes through $\tau_n$ between times $i$ and $j$. For the
sake of simplicity, we suppose that $\wp_n$ enters $\tau_n$ from the
left.\footnote{The case where $\wp_n$ enters $\tau_n$ from the right may
be treated by considering the path $h\mapsto\wp_n(d_{\mathfrak{q}
_n}(a_n,b_n)-h)$ instead of $\wp_n$.} Notice that the path $\wp_n$
does not intersect $\lhb\rho_n, y_n \rhb$, because the labels on
$\lhb
\rho_n, y_n \rhb$ are strictly greater than the labels on $\wp_n$.
Let $k$ be the largest integer in $\lbracket i-1,j \rbracket$ such
that $\wp_n(k)$
belongs to the set $\{\wp_n(i-1)\}\cup\rr\tr_n(\lbracket l_n,t_n
\rbracket)$. Then
$\wp_n(k+1) \in\{\wp_n(j+1)\}\cup\rr\tr_n(\lbracket t_n,r_n
\rbracket)$. Moreover,
$\mathfrak{l}_n(\wp_n(k+1))=\mathfrak{l}_n(\wp_n(k))-1$: otherwise,
all the vertices in
$[\wp_n(k+1),\wp_n(k)]$ would have labels greater than $\mathfrak
{l}_n(\wp
_n(k))$, and it is easy to see that this would prohibit $\wp_n$ from
exiting $\tau_n$ by going ``to the right,'' in the sense that we would
not have $\mathfrak{l}_n(\wp_n(j+1)) = \mathfrak{l}_n(\wp_n(j))-1$.
As a result, when
performing the Chapuy--Marcus--Schaeffer bijection for the arc linking
$\wp_n(k)$ to $\wp_n(k+1)$, we have to visit $y_n$. Then, going through
consecutive successors of $t_n$, we are bound to hit $\wp_n(k+1)$, so
that $d_{\mathfrak{q}_n}(y_n,\wp_n) \le\mathfrak{l}_n(y_n) -
\mathfrak{l}_n(\wp_n(k+1))$. This
yields that $d_{\mathfrak{q}_n}(a_n,y_n) \le d_{\mathfrak
{q}_n}(a_n,b_n) + d_{\mathfrak{q}_n}(y_n,\wp
_n) \le2\eps \g$, and, by taking the limit, $d_\infty(a,y) \le
2\eps$.

We conclude thanks to two lemmas, whose proofs are postponed to
Section~\ref{secprle}. They are derived from similar results in the
planar case:~\cite{legall07tss}, Lem\-ma~2.4, and
\cite{legall08glp}, Corollary 6.2. We call $\lambda$ the volume measure on $\mathfrak
{q}_\infty$,
that is, the image of the Lebesgue measure on $[0,1]$ by the canonical
projection from $[0,1]$ to $\mathfrak{q}_\infty$.
\begin{lem}\label{leme2}
Almost surely, for every $\eta> 0$ and every subtree $\tau$ rooted
at~$\rho$, the condition $\inf_{\tau} \Lab_\infty<\Lab_\infty
(\rho
)-\eta$ implies that
\begin{eqnarray*}
&&\liminf_{\eps\to0} \eps^{-2} \lambda\biggl( \biggl\{ y\in\tau \dvtx
\Lab_\infty(y) < \Lab_\infty(\rho) - \eta+ \eps;\\
&&\hspace*{66.4pt}\forall x \in\lhb\rho,y\rhb, \Lab_\infty(x) > \Lab_\infty
(\rho)-\eta
+\frac\eps8 \biggr\}\biggr) >0.
\end{eqnarray*}
\end{lem}
\begin{lem}\label{leme4}
Let $\delta\in(0,1]$. For every $p\ge1$,
\[
\mathbb{E}\biggl[ \biggl(\sup_{\eps>0}\biggl(\sup_{x\in\mathfrak{q}_\infty} \frac
{\lambda(B_\infty (x,\eps))}{\eps^{4-\delta}}\biggr)\biggr)^p \biggr] < \infty.
\]
\end{lem}

We apply Lemma~\ref{leme2} to $\tau$ and $\eta=\Lab_\infty(\rho
)-\Lab
_\infty(a)>0$, and we find that, for $\eps$ small enough,
\[
\lambda(\mathcal U_\infty^\eps) \ge\eps^{5/2}.
\]
The inclusion $\mathcal U_\infty^\eps\subseteq B_{\infty}(a,2\eps)$
yields that
\[
S\de\sup_{\eps>0}\biggl(\sup_{x\in\mathfrak{q}_\infty} \frac
{\lambda(B_\infty(x,\eps
))}{\eps^{7/2}}\biggr)= \infty.
\]
Lemma~\ref{leme4} applied to $\delta=1/2$ and $p=1$ yields that $S$ is
integrable, so that $S<\infty$ a.s. This is a contradiction.
\end{pf}

\subsection{Set overflown by a path}

We call $\fl_n$ the floor of $\tr_n$. Let $i\in\lbracket0,2n
\rbracket$, and let
$\suc(i)$ be its successor in $(\tr_n,\mathfrak{l}_n)$, defined
by~(\ref{suc}).
%
%
\begin{figure}[b]

\includegraphics{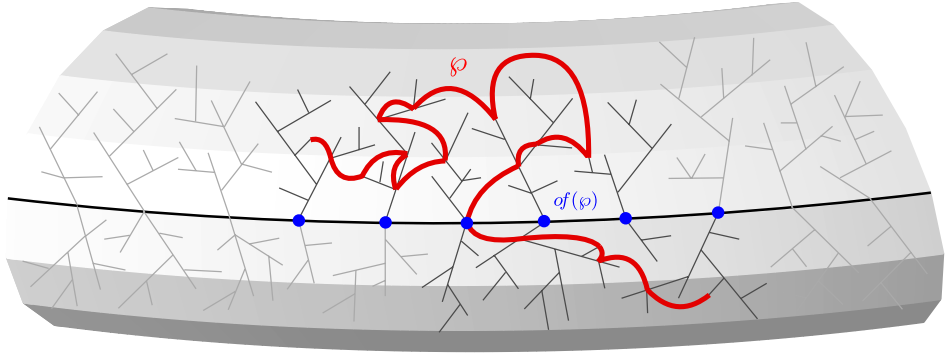}

\caption{The set overflown by the path $\wp$ is the set of (blue)
large dots.}
\label{of}
\end{figure}
We moreover suppose that $\suc(i) \neq\infty$. We say that the arc
linking $\tr_n(i)$ to $\tr_n(\suc(i))$ \textit{overflies} the set
\[
\tr_n (\overrightarrow{\lbracket i,\suc(i) \rbracket} )\cap\fl_n,
\]
where $\overrightarrow{\lbracket i,\suc(i) \rbracket}$ was defined
by~(\ref{oraij}).
We define the set overflown by a path $\wp$ in $\mathfrak{q}_n$ that
avoids the
base point $v_n^\bullet$ as the union of the sets its arcs overfly;
see Figure~\ref{of}.
We denote it by $\of(\wp) \subseteq\fl_n$.
\begin{lem}\label{lemof}
Let $a \sim_\infty b \in\TT_\infty$ and $\alpha$, $\beta\in
\mathfrak{f}_\infty
^\mathfrak{e}\cap\fl_\infty$. We suppose that, for $n$ sufficiently large,
there exist vertices $\alpha_n$, $\beta_n \in\mathfrak
{f}_n^\mathfrak{e}\cap\fl_n$
and $a_n$, $b_n \in\tr_n$ converging, respectively, toward $\alpha$,
$\beta$,~$a$ and~$b$. If, for infinitely many $n$'s, there exists a
geodesic path $\wp_n$ from $a_n$ to $b_n$ that overflies $\lhb\alpha
_n, \beta_n \rhb$, then for all $c \in\lhb\alpha, \beta\rhb$,
\[
\Lab_\infty(c) \ge\Lab_\infty(a)=\Lab_\infty(b).
\]

Moreover, if there exists $c \in\lhb\alpha, \beta\rhb$ for which
$\Lab_\infty(c) = \Lab_\infty(a)$, then \mbox{$a\sim_\infty c$}.
\end{lem}

\begin{pf}
Let $c\in\lhb\alpha, \beta\rhb$. We can find vertices $c_n \in
\lhb
\alpha_n, \beta_n \rhb$ converging to~$c$. By definition, there is an
arc of $\wp_n$ that overflies $c_n$. Say it links a~vertex labeled $l$
to a vertex $v$ labeled $l-1$. From the Chapuy--Marcus--Schaeffer
construction, we readily obtain that $\mathfrak{l}_n(c_n) \ge l$.
Using the
fact that $\mathfrak{l}_n(a_n)-l \le d_{\mathfrak{q}_n}(a_n,b_n)$, we find
\[
\mathfrak{l}_n(c_n) \ge\mathfrak{l}_n(a_n) - d_{\mathfrak{q}_n}(a_n,b_n).
\]

Moreover, we can construct a path from $c_n$ to $v$ going through
consecutive successors of $c_n$. As a result, $d_{\mathfrak
{q}_n}(c_n,\wp_n) \le
\mathfrak{l}_n(c_n)-l+1$, so that
\[
d_{\mathfrak{q}_n}(c_n,a_n) \le\mathfrak{l}_n(c_n) - \mathfrak
{l}_n(a_n) + 2 d_{\mathfrak{q}_n}(a_n,b_n) + 1.
\]

Both claims follow by taking limits in these inequalities after
renormalization, and by using the fact that $d_{\mathfrak
{q}_n}(a_n,b_n) = o(n^{1/4})$.
\end{pf}

\subsection{Points identifications}

We proceed in three steps. We first show that points of $\fl_\infty$
are not identified with any other points, then that points cannot be
identified with their strict ancestors, and finally Theorem~\ref{ip}.

\subsubsection{Floor points are not identified with any other
points}\label{secfp}

\begin{lem}\label{lemfl}
A.s., if $a\in\fl_\infty$ and $b\in\TT_\infty$ are such that
$a\sim
_\infty b$, then $a=b$.
\end{lem}
\begin{pf}
Let $a \in\fl_\infty$ and $b\in\TT_\infty\bs\{a\}$ be such that
$a\sim_\infty b$. We first suppose that~$a$ is not a node. There exists
$\mathfrak{e}\in\vec E(\s_\infty)$ such that $a\in\mathfrak
{f}_\infty^\mathfrak{e}\cap\mathfrak{f}_\infty
^{\bar\mathfrak{e}}$, and we can find $s$, $t$ satisfying $a=\TT
_\infty(s) = \TT
_\infty(t)$, $\mathfrak{e}(s)=\mathfrak{e}$ and $\mathfrak{e}(t)
=\bar\mathfrak{e}$. Without loss of
generality, we may suppose that $s<t$. Until further notice, we will
moreover suppose that $\rho_b \notin\lhb\mathfrak{e}\rhb$.

We restrict ourselves to the case $\s_n=\s_\infty$, which happens
for $n$ sufficiently large. We can find $a_n \in\fl_n$ and $b_n \in
\tr
_n$ converging toward~$a$ and~$b$ and satisfying \mbox{$\rho_{b_n} \notin
\lhb\mathfrak{e}\rhb$}. Let $\wp_n$ be a geodesic path (in
$\mathfrak{q}_n$, for $d_{\mathfrak{q}
_n}$) from $a_n$ to $b_n$. It has to overfly at least $\lhb a_n,
\mathfrak{e}^-
\rhb$ or $\lhb a_n, \mathfrak{e}^+ \rhb$. Indeed, every pair $(x,y)
\in\lhb a_n,
\mathfrak{e}^- \rhb\times\lhb a_n, \mathfrak{e}^+ \rhb$ breaks
$\tr_n$ into connected
components, and the points $a_n$ and $b_n$ do not belong to the same of
these components. There has to be an arc of $\wp_n$ that links a point
belonging to the component containing~$a_n$ to one of the other
components. Such an arc overflies $x$ or $y$.\vadjust{\goodbreak}

Let us suppose that, for infinitely many $n$'s, $\wp_n$ overflies $\lhb
a_n, \mathfrak{e}^- \rhb$. Then, Lemma~\ref{lemof} ensures that
$\Lab_\infty(c)
\ge\Lab_\infty(a) = \Lab_\infty(b)$ for all $c \in\lhb a,
\mathfrak{e}^- \rhb$.
Properties of Brownian snakes show that the labels on $\lhb a ,
\mathfrak{e}^-
\rhb$ are Brownian. Precisely, we may code $\lhb\mathfrak{e}\rhb$
by the
interval $[0,\sigma^\mathfrak{e}]$ as follows. For $x \in[0,\sigma
^\mathfrak{e}]$, we
define $T_x \de\inf\{ r \ge\langle s \rangle \dvtx \CC_\infty(r)=
\CC_\infty(\langle s \rangle) - x \}$. Then $\lhb\mathfrak{e}\rhb
= \TT_\infty
(\{T_x, 0 \le x \le\sigma_\mathfrak{e}\})$, and
\[
\bigl(\Lab_\infty(T_x)-\Lab_\infty(\langle s \rangle)\bigr)_{0 \le
x \le
\sigma_\mathfrak{e}}=(\mathfrak{M}^\mathfrak{e}_\infty(x)
)_{0 \le x \le\sigma_\mathfrak{e}},
\]
where, conditionally given $\mathfrak I_\infty$, the process
$\mathfrak{M}^\mathfrak{e}
_\infty$ (defined during Proposition~\ref{cvint}) has the law of a
certain Brownian bridge. Using the fact that local minimums of Brownian
motion are distinct, we can find $d\in\lhb a, \mathfrak{e}^- \rhb\bs
\{a\}$
such that $\Lab_\infty(c) > \Lab_\infty(a)$ for all $c \in\lhb a, d
\rhb\bs\{a\}$.

Because $a\in\fl_\infty$, $s$ and $t$ are both increase points of
$\CC
_\infty$ and thus are not increase points of $\Lab_\infty$, by
Lemma~\ref{pc}. As a result, there exist two trees $\tau^1 \subseteq
\mathfrak{f}
_\infty^\mathfrak{e}$ and $\tau^2 \subseteq\mathfrak{f}_\infty
^{\bar\mathfrak{e}}$ rooted at $\rho
^1$, $\rho^2 \in\lhb a, d \rhb\bs\{a\}$ satisfying $\inf_{\tau^i}
\Lab
_\infty< \Lab_\infty(a) < \Lab_\infty(\rho^i)$ (see Figure~\ref
{tau1tau2}).

%
%
\begin{figure}

\includegraphics{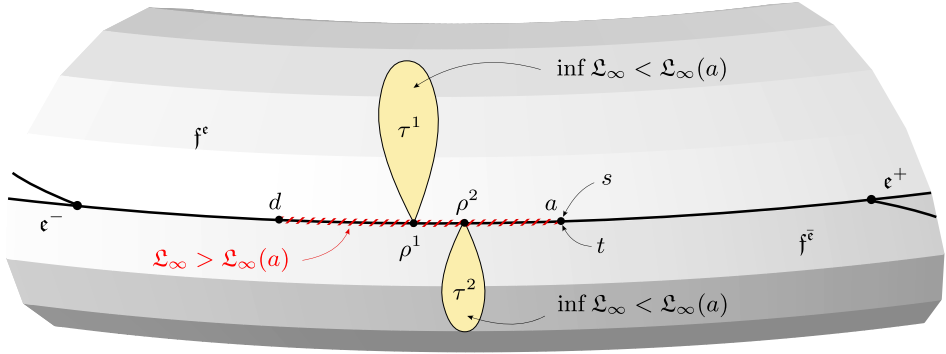}

\caption{The trees $\tau^1$ and $\tau^2$.}
\label{tau1tau2}
\end{figure}

%
%
\begin{figure}

\includegraphics{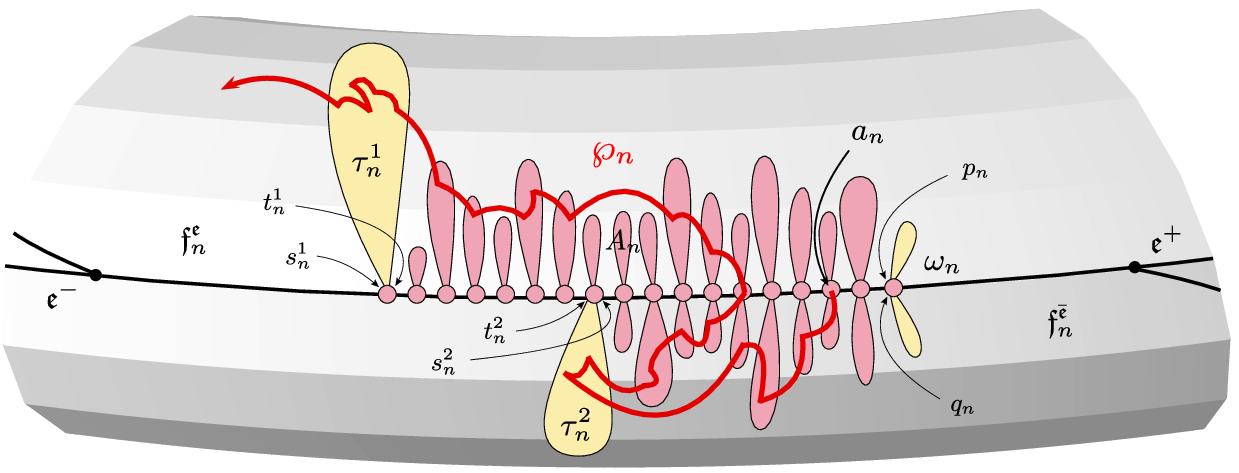}

\caption{The path $\wp_n$ passing through the tree $\tau_n^1$.}
\label{taun12}
\end{figure}

Similarly, if for infinitely many $n$'s, $\wp_n$ overflies $\lhb a_n,
\mathfrak{e}
^+ \rhb$, then we can find two trees $\tau^3 \subseteq\mathfrak
{f}_\infty^\mathfrak{e}$
and $\tau^4 \subseteq\mathfrak{f}_\infty^{\bar\mathfrak{e}}$
rooted at $\rho^3$, $\rho^4
\in\lhb a, \mathfrak{e}^+ \rhb\bs\{a\}$ satisfying $\inf_{\tau^i}
\Lab_\infty<
\Lab_\infty(a) < \Lab_\infty(\rho^i)$, and $\Lab_\infty(c) >
\Lab_\infty
(a)$ for all $c\in\lhb\rho^3,\rho^4 \rhb$. Three cases may occur:
\begin{longlist}
\item for $n$ large enough, $\wp_n$ does not overfly $\lhb a_n,
\mathfrak{e}^+
\rhb$ (and therefore overflies $\lhb a_n, \mathfrak{e}^- \rhb$);
\item for $n$ large enough, $\wp_n$ does not overfly $\lhb a_n,
\mathfrak{e}^-
\rhb$ (and therefore overflies $\lhb a_n, \mathfrak{e}^+ \rhb$);
\item for infinitely many $n$'s, $\wp_n$ overflies $\lhb a_n, \mathfrak
{e}^+ \rhb
$, and for infinitely many $n$'s, $\wp_n$ overflies $\lhb a_n,
\mathfrak{e}^- \rhb$.
\end{longlist}

In case (i), the trees $\tau^1$ and $\tau^2$ are well defined. Let
$\tau_n^1 \subseteq\mathfrak{f}_n^\mathfrak{e}$, $\tau_n^2
\subseteq\mathfrak{f}_n^{\bar\mathfrak{e}}$ be
trees rooted at $\rho_n^1$, $\rho_n^2 \in\lhb a_n, \mathfrak{e}^-
\rhb$
converging to $\tau^1$ and $\tau^2$. We claim that, for $n$\vadjust{\goodbreak}
sufficiently large, $\wp_n$ passes through $\tau_n^1$ or $\tau_n^2$.
First, notice that, for $n$ large enough, $\wp_n \cap\lhb\rho_n^1,
\rho_n^2 \rhb= \varnothing$. Otherwise, for infinitely many $n$'s, we
could find $\alpha_n \in\wp_n \cap\lhb\rho_n^1, \rho_n^2 \rhb$, and,
up to extraction, we would have $\alpha_n \to\alpha\in\lhb\rho^1,
\rho^2 \rhb\subseteq\lhb a, d \rhb\bs\{a\}$. Furthermore,
$d_{\mathfrak{q}
_n}(a_n, \alpha_n) \le d_{\mathfrak{q}_n}(a_n,b_n)$ so that $a \sim
_\infty\alpha
$, and $\Lab_\infty(a) = \Lab_\infty(\alpha)$ by Lemma~\ref{dssb},
which is impossible. For $n$ even larger, it holds that $\inf_{\tau
_n^i} \mathfrak{l}_n < \inf_{\wp_n} \mathfrak{l}_n$. Roughly
speaking,\vspace*{1pt} $\wp_n$ cannot
go from a tree located at the right of $\tau_n^1$ (resp., at the left of
$\tau^2_n$) to a tree located at its left in $\mathfrak{f}^\mathfrak
{e}_n$ (resp., to a tree
located at its right in $\mathfrak{f}_n^{\bar\mathfrak{e}}$) without
entering it. Then $\wp
_n$ has to enter $\tau_n^1$ from the right or $\tau_n^2$ from the left
and pass through one of these trees (see Figure~\ref{taun12}).

More precisely, we call $\lbracket s_n^1,t_n^1 \rbracket$ and
$\lbracket s_n^2,t_n^2 \rbracket$
the sets coding the subtrees~$\tau_n^1$ and~$\tau_n^2$. Let $\omega_n
\in\lhb a_n, \mathfrak{e}^+ \rhb$ be a point that is not overflown
by $\wp_n$,
$p_n \de\inf\{t_n^1 \le r \le2n \dvtx \omega_n = \rr\tr_n(r) \}$
and $q_n \de\sup\{0 \le r \le s_n^2 \dvtx \omega_n = \rr\tr_n(r)
\}
$. Then, we let
\[
A_n \de\rr\tr_n(\overrightarrow{\lbracket t_n^1,p_n \rbracket}
\cup
\overrightarrow{\lbracket q_n,s_n^2 \rbracket} ).
\]

We call $\wp_n(i-1)$ the last point of $\wp_n$ belonging to $A_n$. Such
a point exists because $a_n \in A_n$ and $b_n \notin A_n$. The remarks
in the preceding paragraphs yield that neither $\wp_n(i-1)$ nor $\wp
_n(i)$ belong to $\lhb\rho_n^1,\rho_n^2 \rhb$, and, because of the way
arcs are constructed in the Chapuy--Marcus--Schaeffer bijection, we see
that \mbox{$\wp_n(i) \in\tau_n^1 \cup\tau_n^2$}. Without loss of generality,
we may assume that $\wp_n(i) \in\tau_n^1$. Because $\wp_n$ does not
overfly $\omega_n$, it enters $\tau_n^1$ from the right at time $i$,
that is, $\mathfrak{l}_n(\wp_n(i)) = \mathfrak{l}_n(\wp_n(i-1))+1$.
Let $\wp_n(j+1)$ be
the first\vspace*{1pt} point after $\wp_n(i)$ not belonging to $\tau_n^1$. It exists
because $b_n \notin\tau_n^1$. Then, because $\wp_n(j+1) \notin A_n$
and $\wp_n$ does not overfly $\omega_n$, we see that
$\mathfrak{l}_n(\wp_n(j+1)) = \mathfrak{l}_n(\wp_n(j))+1$, so
that $\wp_n$ passes
through $\tau_n^1$ between times $i$ and $j$.

In case (ii), we apply the same reasoning with $\tau^3$ and $\tau^4$
instead of $\tau^1$ and $\tau^2$. In case (iii), the four trees
$\tau
^1$, $\tau^2$, $\tau^3$ and $\tau^4$ are well defined, and we obtain
that $\wp_n$ has to pass through one of their discrete approximations.
We then conclude by Lemma~\ref{legalllem} that $a\not\sim_\infty b$,
which contradicts our hypothesis.

We treat the case where $\rho_b \in\lhb\mathfrak{e}\rhb\bs\{a\}$
in a similar
way, simply by replacing~$\mathfrak{e}^+$ (resp., $\mathfrak{e}^-$)
by $\rho_{b}$ if $\rho
_{b} \in\lhb a, \mathfrak{e}^+ \rhb$ (resp., $\rho_{b} \in\lhb a,
\mathfrak{e}^- \rhb$).
When~$a$ is a node, we apply the same arguments, finding up to six
trees (one for each forest containing~$a$). Finally, if $\rho_b = a$,
then~$a$ is a strict ancestor of~$b$. This will be a particular case of
Lemma~\ref{lemad}.
\end{pf}

\subsubsection{Points are not identified with their strict ancestors}

\begin{lem}\label{lemad}
A.s., for every $a,b \in\TT_\infty$ such that $\rho_a =\rho_b$ and $a
\prec b$, we have $a \not\sim_\infty b$.
\end{lem}

The proof of this lemma uses the same kind of arguments we used in
Section~\ref{secfp}, is slightly easier than the proof of Lemma \ref
{lemfl} and is very similar to Le Gall's proof for Proposition 4.2
in~\cite{legall07tss}, so that we leave the details to the reader.

\subsubsection{Points~$a$,~$b$ are only identified when $d_\infty^\circ(a,b) =0$}

\begin{lem}\label{lemst}
A.s., for every tree $\tau\subseteq\TT_\infty$ rooted at $\rho\in
\fl
_\infty$ and all $a,b \in\tau\bs\{\rho\}$ satisfying $a\sim
_\infty b$,
we have $d_\infty^\circ(a,b) =0$.
\end{lem}
\begin{pf}
Let $\tau\subseteq\TT_\infty$ be a tree rooted at $\rho\in\fl
_\infty
$ and $a,b \in\tau\bs\{\rho\}$ satisfying $a \neq b$ and $a\sim
_\infty
b$. By Lemma~\ref{lemad}, we know that $a \not\prec b$ and $b \not
\prec
a$. As a~consequence, we have either $s < t$ for all $(s,t) \in\TT
_\infty^{-1}(a) \times\TT_\infty^{-1}(b)$ or $s > t$ for all $(s,t)
\in\TT_\infty^{-1}(a) \times\TT_\infty^{-1}(b)$. Without loss of
generality, we will assume that the first case occurs. Let us suppose
that $d_\infty^\circ(a,b) > 0$. By Lemma~\ref{dssb}, we know that
$\Lab
_\infty(a)=\Lab_\infty(b)$, and by~(\ref{dzero}), we have both
$\inf
_{[a,b]} \Lab_\infty< \Lab_\infty(a)$ and $\inf_{[b,a]} \Lab
_\infty<
\Lab_\infty(a)$. As a result, there are two subtrees $\tau^1
\subseteq
[a,b]$ and $\tau^2 \subseteq[b,a]$ rooted at $\rho^1 \in\lhb a,b
\rhb
\bs\{a,b\}$ and $\rho^2 \in(\lhb\rho,a \rhb\cup\lhb\rho,b \rhb
\cup\fl_\infty) \bs\{a,b\}$ satisfying $\inf_{\tau^i} \Lab
_\infty<
\Lab_\infty(a)$.

Let\vspace*{1pt} $\tau_n \subseteq\tr_n$ be a tree rooted at $\rho_n$ and $a_n$,
$b_n \in\tr_n$ be points converging to $\tau$,~$a$, and~$b$. Let
$\tau
_n^1 \subseteq[a_n,b_n]$ and $\tau_n^2 \subseteq[b_n,a_n]$ be
subtrees rooted at $\rho_n^1 \in\lhb a_n,b_n \rhb\bs\{a_n,b_n\}$ and
$\rho_n^2 \in(\lhb\rho_n,a_n \rhb\cup\lhb\rho_n,b_n \rhb\cup
\fl
_n) \bs\{a_n,b_n\}$ converging\vspace*{1pt} toward $\tau^1$ and~$\tau^2$. We
consider a geodesic path $\wp_n$ from $a_n$ to $b_n$. Recall that
$a\sim
_\infty b$ implies that $d_{\mathfrak{q}_n}(a_n,b_n)=o(n^{1/4})$.

Because every point in $\lhb\rho, \rho^1 \rhb$ is a strict ancestor
to~$a$ or~$b$, for $n$ large enough, $\wp_n$ does not intersect $\lhb
\rho_n, \rho_n^1 \rhb$. Otherwise, we could find an accumulation point
$\alpha$ identified with~$a$ and~$b$, such that $\alpha\prec a$ or
$\alpha\prec b$ (possibly both), and this would contradict Lemma \ref
{lemad}. If $\rho^2 \in\tau$, for $n$ large, $\wp_n$ does not
intersect $\lhb\rho_n, \rho_n^2 \rhb$ either. The same reasoning
yields that $\wp_n$ does not intersect $\fl_n$ for $n$ sufficiently
large, because of Lemma~\ref{lemfl}.

Let $\lbracket s_n^1,t_n^1 \rbracket$ and $\lbracket s_n^2,t_n^2
\rbracket$ be the sets coding
the subtrees $\tau_n^1$ and $\tau_n^2$. We let
\[
A_n \de\rr\tr_n(\overrightarrow{\lbracket t_n^2,s_n^1 \rbracket}
)\quad\mbox{and}\quad B_n
\de\rr\tr_n(\overrightarrow{\lbracket t_n^1,s_n^2 \rbracket} ).\vadjust{\goodbreak}
\]
By convention, if $\rho_n^2 \notin\mathfrak{f}_n^\mathfrak{e}$, we
set $\lhb\rho_n, \rho
_n^2 \rhb\de\varnothing$. It is easy to see that $a_n \in A_n$, $b_n
\in B_n$, $ A_n \cap B_n \subseteq\lhb\rho_n, \rho_n^1 \rhb\cup
\lhb
\rho_n, \rho_n^2 \rhb\cup\fl_n$ and $A_n \cup B_n \cup\tau_n^1
\cup
\tau_n^2 = \tr_n$.

We conclude as in the proof of Lemma~\ref{lemfl}. We call $\wp_n(i-1)$
the last point of $\wp_n$ belonging to $A_n$. Such a point exists
because $a_n \in A_n$ and $b_n \notin A_n$. The remarks in the
preceding paragraphs yield that, for $n$ large enough, neither $\wp
_n(i-1)$ nor $\wp_n(i)$ belong to $A_n \cap B_n$. For $n$ even larger,
$\inf_{\tau_n^j} \mathfrak{l}_n < \inf_{\wp_n} \mathfrak{l}_n$,
and because of the way
arcs are constructed in the Chapuy--Marcus--Schaeffer bijection, we see
that $\wp_n(i) \in\tau_n^1 \cup\tau_n^2$. The path $\wp_n$ either
enters~$\tau_n^1$ from the left or enters $\tau_n^2$ from the right.
Without loss of generality, we may suppose that $\wp_n(i) \in\tau
_n^1$. Let $\wp_n(i'+1)$ be the first point after $\wp_n(i)$ not
belonging to $\tau_n^1$. Then $\wp_n(i'+1) \in B_n \cup\tau_n^2$.
If $\wp_n$ passes through $\tau_n^1$ between times $i$ and $i'$, we are
done. Otherwise, $\wp_n(i'+1) \in\tau_n^2$ because of the condition
$\inf_{\tau_n^2} \mathfrak{l}_n < \inf_{\wp_n} \mathfrak{l}_n$
(informally,\vspace*{-1pt} $\wp_n$
cannot pass over $\tau_n^2$ without entering it). We consider the first
point $\wp_n(i''+1)$ after $\wp_n(i')$ not belonging to $\tau_n^2$, and
reiterate the argument. Because $\wp_n$ is a finite path, we see
that~$\wp_n$ will eventually pass through $\tau_n^1$ or $\tau_n^2$;
see Figure~\ref{fig10}.

%
%
\begin{figure}

\includegraphics{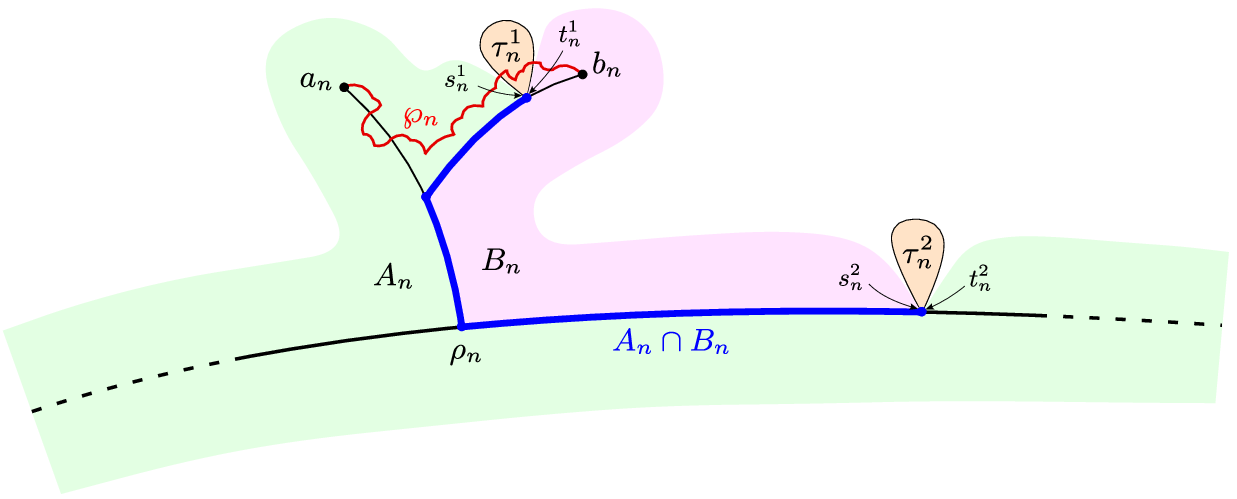}

\caption{The path $\wp_n$ passing through the subtree
$\tau_n^1$.}\label{fig10}
\end{figure}

If $\wp_n$ passes through $\tau_n^1$ (resp., $\tau_n^2$) for infinitely
many $n$'s, a reasoning similar to the one we used in the proof of
Lemma~\ref{lemof} yields that $\Lab_\infty(\rho^1) > \Lab_\infty(a)$
[resp., $\Lab_\infty(\rho^2) > \Lab_\infty(a)$]. We conclude by
Lemma~\ref{legalllem} that $a \sim_\infty b$. This is a contradiction.
\end{pf}
\begin{lem}\label{lemdt}
A.s., for all $a,b \in\TT_\infty\bs\fl_\infty$ such that $\rho_a
\neq\rho_b$ and $a\sim_\infty b$, we have $d_\infty^\circ(a,b) =0$.
\end{lem}
\begin{pf}
The proof of this lemma is very similar to that of Lemma~\ref{lemst}.
Let $a,b \in\TT_\infty\bs\fl_\infty$ be such that $\rho_a \neq
\rho
_b$ and $a\sim_\infty b$. Here again, we may suppose that $s < t$ for
all $(s,t) \in\TT_\infty^{-1}(a) \times\TT_\infty^{-1}(b)$, and we
can find two subtrees $\tau^1 \subseteq[a,b]$ and $\tau^2 \subseteq
[b,a]$ rooted at $\rho^1, \rho^2 \in(\lhb\rho_a,a \rhb\cup\lhb
\rho
_b,b \rhb\cup\fl_\infty)\bs\{a,b\}$ satisfying $\inf_{\tau^i}
\Lab
_\infty< \Lab_\infty(a)$. As before, we consider the discrete
approximations $a_n$, $b_n$, $\tau_n^1 = \rr\tr_n(\lbracket
s_n^1,t_n^1 \rbracket)$
and $\tau_n^2 = \rr\tr_n(\lbracket s_n^2,t_n^2 \rbracket)$ of~$a$,
$b$, $\tau^1$\vadjust{\goodbreak}
and $\tau^2$. Let $\wp_n$ be a geodesic path from $a_n$ to $b_n$. We
still define
\[
A_n \de\rr\tr_n(\overrightarrow{\lbracket t_n^2,s_n^1 \rbracket}
)\quad\mbox{and}\quad B_n
\de\rr\tr_n(\overrightarrow{\lbracket t_n^1,s_n^2 \rbracket} ),
\]
and we see by the same arguments as in Lemma~\ref{lemst} that, for $n$
sufficiently large, $\wp_n$ does not intersect $A_n \cap B_n$. We then
conclude exactly as before.~%
\end{pf}

Theorem~\ref{ip} follows from Lemmas~\ref{lemfl},~\ref{lemad}, \ref
{lemst} and~\ref{lemdt}. A straightforward consequence of Theorem \ref
{ip} is that, if the equivalence class of $a=\TT_\infty(s)$ for~$\sim
_\infty$ is not trivial, then $s$ is an increase point of $\Lab
_\infty
$. By Lemma~\ref{pc}, the equivalence class of~$a$ for $\simeq_\infty$
is then trivial. Such points may be called \textit{leaves} by analogy
with tree terminology.


\section{1-regularity of quadrangulations}\label{secsurf}

The goal of this section is to prove Theorem~\ref{cvq2}. To that end,
we use the notion of regular convergence, introduced by Whyburn
\cite{whyburn35rcm}.

\subsection{1-regularity}

Recall that $(\MM,\dGH)$ is the set of isometry classes of compact
metric spaces, endowed with the Gromov--Hausdorff metric. We say that a
metric space $(\X,\delta)$ is a \textit{path metric space} if any two
points $x,y \in\X$ may be joined by a path isometric to a real
segment---necessarily of length $\delta(x,y)$. We call $\PM$ the set of
isometry classes of path metric spaces. By
\cite{burago01cmg}, Theorem 7.5.1, $\PM$ is a closed subset of $\MM$.
\begin{defi}
We say that a sequence $(\X_n)_{n \ge1}$ of path metric spaces is
\textit{$1$-regular} if for every $\eps> 0$, there exists $\delta>0$
such that for $n$ large enough, every loop of diameter less than
$\delta
$ in $\X_n$ is homotopic to $0$ in its $\eps$-neighborhood.
\end{defi}

This definition is actually slightly stronger than Whyburn's original
definition~\cite{whyburn35rcm}. See the discussion in the second
section of~\cite{miermont08sphericity} for more details. We also chose
here not to restrict the notion of 1-regularity only to converging
sequences of path metric spaces, as it was done in
\cite{miermont08sphericity,whyburn35rcm}, because the notion of 1-regularity
(as stated here) is not directly related to the convergence of the
sequence of path metric spaces. Our main tool is the following theorem,
which is a simple consequence of Begle~\cite{begle44rc}, Theorem 7.
\begin{prop}\label{begle}
Let $(\X_n)_{n \ge1}$ be a sequence of path metric spaces all
homeomorphic to the $g$-torus $\TTT_g$. Suppose that $\X_n$ converges
toward $\X$ for the Gromov--Hausdorff topology, and that the sequence
$(\X_n)_{n \ge1}$ is 1-regular. Then $\X$ is 
homeomorphic to $\TTT_g$ as well.
\end{prop}

\subsection{Representation as metric surfaces}

In order to apply Proposition~\ref{begle}, we construct a path metric
space $(\cS_n,\delta_n)$ homeomorphic to $\TTT_g$, and an embedded
graph that is a representative of the map $\mathfrak{q}_n$, such that the
restriction of $(\cS_n,\delta_n)$ to the embedded graph is isometric to
$( V(\mathfrak{q}_n), d_{\mathfrak{q}_n} )$. We use the
method provided by Miermont
in~\cite{miermont08sphericity}, Section 3.1.\vadjust{\goodbreak}

We write $F(\mathfrak{q}_n)$ the set of faces of $\mathfrak{q}_n$.
Let $(X_f,D_f)$, $f \in
F(\mathfrak{q}_n)$ be $n$ copies of the hollow bottomless unit cube
\[
X_f \de[0,1]^3\bs\bigl((0,1)^2 \times[0,1) \bigr)
\]
endowed with the intrinsic metric $D_f$ inherited from the Euclidean
metric. (The distance between two points $x$ and $y$ is the Euclidean
length of a~minimal path in $X_f$ linking $x$ to $y$.)

Let $f\in F(\mathfrak{q}_n)$, and let $e_1$, $e_2$, $e_3$ and $e_4$ be
the four
half-edges incident to~$f$, ordered according to the counterclockwise
order. For $0 \le t \le1$, we define:
\begin{eqnarray*}
c_{e_1}(t) &=& (t,0,0) \in X_f;\\
c_{e_2}(t) &=& (1,t,0) \in X_f;\\
c_{e_3}(t) &=& (1-t,1,0) \in X_f;\\
c_{e_4}(t) &=& (0,1-t,0) \in X_f.
\end{eqnarray*}
In this way, we associate with every half-edge $e\in\vec E(\mathfrak
{q}_n)$ a
path along one of the four edges of the square $\partial X_f$,
where $f$ is the face located to the left of~$e$.\looseness=-1

We then define the relation $\approx$ as the coarsest equivalence
relation for which $c_e(t) \approx c_{\bar e}(1-t)$ for all $e\in\vec
E(\mathfrak{q}_n)$ and $t \in[0,1]$. This corresponds to gluing the spaces
$X_f$'s along their boundaries according to the map structure
of $\mathfrak{q}_n$.
The topological quotient $\cS_n \de(\coprod_{f\in F(\mathfrak{q}_n)}
X_f)_{/\approx}$ is a two-dimensional CW-complex satisfying the
following. Its $1$-skeleton $\EE_n = (\coprod_{f\in F(\mathfrak
{q}_n)} \partial
X_f)_{/\approx}$ is an embedding of $\mathfrak{q}_n$ with faces $X_f
\bs\partial
X_f$. To the edge $\{e,\bar e\} \in E(\mathfrak{q}_n)$ corresponds the edge
of $\cS_n$ made of the equivalence class of the points in $c_e([0,1])$.
Its $0$-skeleton $\V_n$ is in one-to-one correspondence with
$V(\mathfrak{q}_n)$.
Its vertices are the equivalence classes of the corners of the
squares $\partial X_f$.

We endow the space $\coprod_{f\in F(\mathfrak{q}_n)} X_f$ with the largest
pseudo-metric $\delta_n$ compatible with $D_f$, $f\in F(\mathfrak
{q}_n)$ and
$\approx$, in the sense that $\delta_n(x,y) \le D_f(x,y)$ for $x,y\in
X_f$ and $\delta_n(x,y)=0$ whenever $x \approx y$. Its quotient---still
noted $\delta_n$---then defines a pseudo-metric on $\cS_n$ (which
actually is a true metric, as we will see in Proposition~\ref{surf}).
As usual, we define $\delta_{(n)} \de\delta_n /\g$ its rescaled version.

We rely on the following proposition. It was actually stated in
\cite{miermont08sphericity} for the two-dimensional sphere but readily
extends to the $g$-torus.
\begin{prop}[(\cite{miermont08sphericity}, Proposition 1)]\label{surf}
The space $(\cS_n,\delta_n)$ is a path metric space homeomorphic
to $\TTT_g$. Moreover, the restriction of $\cS_n$ to $\V_n$ is
isometric to $(V(\mathfrak{q}_n), d_{\mathfrak{q}_n})$, and any
geodesic path in $\cS_n$
between points in $\V_n$ is a concatenation of edges of $\cS_n$.
Finally, $\dGH( (V(\mathfrak{q}_n), d_{\mathfrak{q}_n}), (\cS
_n,\delta_n) ) \le
3$, so that, by Proposition~\ref{cvq},
\[
\bigl(\cS_{n_k},\delta_{(n_k)} \bigr) \tolk(\mathfrak{q}_\infty,
d_\infty)
\]
in the sense of the Gromov--Hausdorff topology.\vadjust{\goodbreak}
\end{prop}

\subsection{\texorpdfstring{Proof of Theorem \protect\ref{cvq2}}{Proof of Theorem 2}}\label{seccvreg}

We prove here that $(\mathfrak{q}_\infty,d_\infty)$ is a.s.
homeomorphic to $\TTT
_g$ by means of Propositions~\ref{begle} and~\ref{surf}.
To this end,
we only need to show that the sequence $(\cS
_{n_k},\delta
_{(n_k)} )_k$ is 1-regular. At first, we only consider simple loops
made of edges. We proceed in two steps: Lemma~\ref{eps0} shows that
there are no noncontractible ``small'' loops; then Lemma~\ref{1reg}
states that the small loops are homotopic to $0$ in their $\eps$-neighborhood.
\begin{lem}\label{eps0}
A.s., there exists $\eps_0 >0$ such that for all $k$ large enough, any
noncontractible simple loop made of edges in $\cS_{n_k}$ has diameter
greater than $\eps_0$.
\end{lem}
\begin{pf}
The basic idea is that a noncontractible loop in $\cS_n$ has to
intersect $\fl_n$ and to ``jump'' from a forest to another one. At the
limit, the loop transits from a forest to another by visiting two
points that $\sim_\infty$ identifies. If the loops vanish at the limit,
then these two identified points become identified with a point in $\fl
_\infty$, creating an increase point for both $\Lab_\infty$ et $\CC
_\infty$. We proceed to the rigorous proof.

We argue by contradiction and assume that, with positive probability,
along some (random) subsequence of the sequence $(n_k)_{k\ge0}$, there
exist noncontractible simple loops $\wp_n$ made of edges in $\cS_n$
with diameter tending to~$0$ (with respect to the rescaled metric
$\delta_{(n)}$). We reason on this event.

Because $\wp_n$ is noncontractible, it has to intersect $\fl_n$: if
not, $\wp_n$ would entirely be drawn in the unique face of $\s_n$,
which is homeomorphic to a disk, by definition of a map. It would thus
be contractible, by the Jordan curve theorem. Let $a_n \in{\wp_n \cap
{\fl}_n}$. Up to further extraction, we may suppose that $a_n \to a \in
\fl
_\infty$. Notice that every time $\wp_n$ intersects $\fl_n$, it has to
be ``close'' to $a_n$. Precisely, if $b_n\in\wp_n\cap\fl_n$ tends
to~$b$, then $\delta_{(n)}(a_n,b_n) \le\diam(\wp_n) \to0$, which
yields $a \sim_\infty b$, and $a=b$ by Lemma~\ref{lemfl}. Moreover,
for $n$ sufficiently large, the base point $v_n^\bullet\notin\wp_n$:
otherwise, for infinitely many $n$'s, $(\mathfrak{l}_{n}(a_n) - \min
\mathfrak{l}_n
+1)/\g\le\diam(\wp_n) \to0$, so that $\Lab_\infty$ would reach its
minimum at~$a$, and we know by Lemma~\ref{min} that this is not the case.

Let us first suppose that~$a$ is not a node of $\TT_\infty$. There
exists $\mathfrak{e}\in\vec E(\s_\infty)$ such that $a \in
\mathfrak{f}^\mathfrak{e}_\infty\cap\mathfrak{f}
^{\bar\mathfrak{e}}_\infty$ and for $n$ large enough, $a_n \in
\mathfrak{f}^\mathfrak{e}_n \cap\mathfrak{f}
^{\bar\mathfrak{e}}_n$. For $n$ even larger, the whole loop $\wp_n$
``stays in $\mathfrak{f}
^\mathfrak{e}_n \cup\mathfrak{f}^{\bar\mathfrak{e}}_n$.''
Precisely, for all $\mathfrak{e}' \in\vec E(\s
_\infty) \bs\{\mathfrak{e},\bar\mathfrak{e}\}$, we have $\wp_n
\cap\mathfrak{f}_n^{\mathfrak{e}'} = \varnothing
$. Otherwise, since $\vec E(\s_\infty)$ is finite, there would exist
$\mathfrak{e}
'\notin\{\mathfrak{e},\bar\mathfrak{e}\}$ such that for infinitely
many $n$'s, we can find
$c_n \in\wp_n \cap\mathfrak{f}_n^{\mathfrak{e}'}$. Up to
extraction, $c_n \to c \in\mathfrak{f}
_\infty^{\mathfrak{e}'}$, so that $c \neq a$ ($a$ is not a node) and
$c \sim
_\infty a$, which is impossible, by Lemma~\ref{lemfl}.

%
%
\begin{figure}

\includegraphics{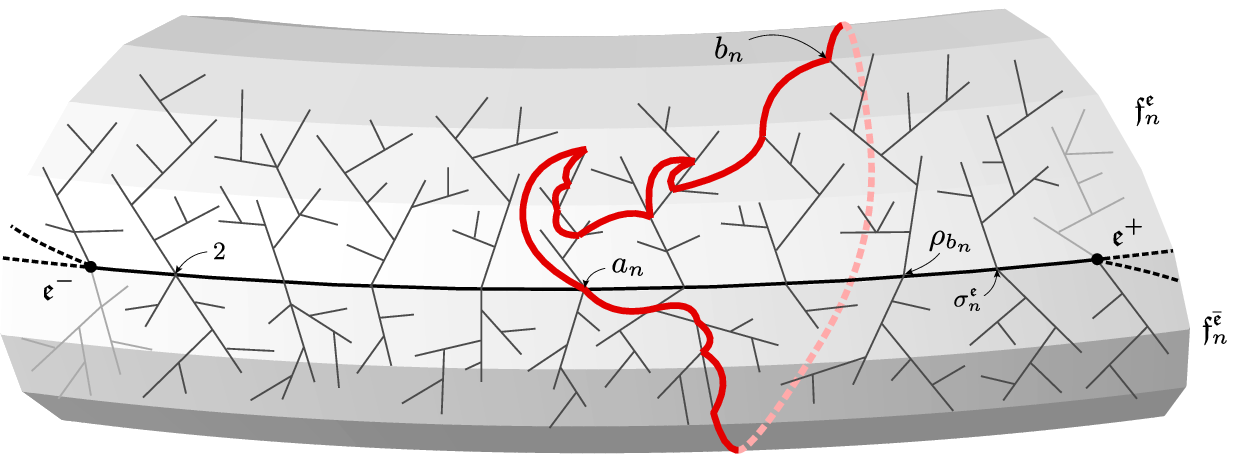}

\caption{A noncontractible loop intersecting $\fl_n$ at $a_n$ and
``jumping'' from $\mathfrak{f}_n^\mathfrak{e}$ to $\mathfrak
{f}_n^{\bar\mathfrak{e}}$ at~$b_n$.}
\label{archmag}
\end{figure}

We claim that there exists an arc of $\wp_n$ linking a point $b_n \in
\mathfrak{f}
_n^\mathfrak{e}$ to some point in $\mathfrak{f}_n^{\bar\mathfrak
{e}}$ that overflies either $\lhb\rho
_{b_n}, \mathfrak{e}^+ \rhb$ or $\lhb\mathfrak{e}^-,\rho_{b_n}
\rhb$ (see Figure \ref
{archmag}). Let us suppose for a moment that this does not hold. In
particular, there is no arc linking a point in $\mathfrak
{f}_n^\mathfrak{e}\bs\fl_n$ to a
point in $\mathfrak{f}_n^{\bar\mathfrak{e}} \bs\fl_n$. It will be
more convenient here to
write $\wp_n$ as
$(a_n=v_1, \alpha_1, v_2, \alpha_2, \ldots, v_{r-1}, \alpha_{r-1},
v_r =
a_n)$ where the $v_i$'s are vertices, and the $\alpha_i$'s are arcs. Let
$i \de\inf\{ j \in\lbracket2,r \rbracket \dvtx v_j \in\fl_n \}$\vadjust{\goodbreak}
be the index of
the first time $\wp_n$ returns to $\fl_n$. Then $v_2, \ldots, v_{i-1}$
belong to the same set $\mathfrak{f}_n^\mathfrak{e}\bs\fl_n$ or
$\mathfrak{f}_n^{\bar\mathfrak{e}} \bs\fl
_n$, and $(\alpha_1, v_2, \alpha_2, \ldots, v_{i-1}, \alpha_{i-1})$ is
thus drawn inside the face of $\s_n$. As a result, the path
$(v_1,\alpha
_1, v_2, \ldots, v_{i-1}, \alpha_{i-1}, v_i)$ is homotopic to the
segment $\lhb v_1,v_i \rhb$. Repeating the argument for every
``excursion'' away from~$\fl_n$, we see that $\wp_n$ is homotopic to a
finite concatenation of segments all included in the topological
segment $\lhb 2, \sigma_n^\mathfrak{e}\rhb$, where we used the
notation of
Section~\ref{secfor} for the forest $\mathfrak{f}_n^\mathfrak{e}$;
see Figure \ref
{archmag}. It follows that $\wp_n$ is contractible, which is a contradiction.

We consider the case where the arc from the previous paragraph
overflies $\lhb\rho_{b_n}, \mathfrak{e}^+ \rhb$. The other case is
treated in a
similar way. From the construction of the Chapuy--Marcus--Schaeffer
bijection, we can find integers $s_n \le t_n$ such that $b_n=\rr\tr
_n(s_n)$, $\mathfrak{e}^+=\rr\tr_n(t_n)$ and for all $s_n \le r \le
t_n$, $\Lab
_n(r) \ge\Lab_n(s_n)$. Up to further extraction, we may suppose that
$s_n/2n \to s$ and $t_n/2n\to t$. Therefore, for all $s \le r \le t$,
$\Lab_\infty(r) \ge\Lab_\infty(s)$. Moreover, the fact that
$b_n\to a
\neq\mathfrak{e}^+$ yields $s <t$, so that $s$ is an increase point
for $\Lab
_\infty$. But $\TT_\infty(s)=a$ and $s$ has to be an increase point for
$\CC_\infty$. By Lemma~\ref{pc}, this cannot happen.

If~$a$ is a node, there are three half-edges $\mathfrak{e}_1$,
$\mathfrak{e}_2$ and $\mathfrak{e}_3$
such that $a = \mathfrak{e}_1^+ =\mathfrak{e}_2^+=\mathfrak{e}_3^+$.
A reasoning similar to what
precedes yields the existence of an arc of $\wp_n$ linking a
point $b_n$ in one of the three sets $\mathfrak{f}^{\mathfrak{e}_i}
\cup\mathfrak{f}^{\bar\mathfrak{e}
_{i+1}}$, $i=1,2,3$ (where we use the convention $\mathfrak
{e}_4=\mathfrak{e}_1$) to a
point lying in another one of these three sets that overflies either,
if $b_n \in\mathfrak{f}_\infty^{\mathfrak{e}_i}$, $\lhb\rho
_{b_n}, a\rhb\cup\lhb\mathfrak{e}
_{i+1} \rhb$ or $\lhb\mathfrak{e}_i^-,\rho_{b_n} \rhb$, or, if
$b_n \in\mathfrak{f}_\infty
^{\mathfrak{e}_{i+1}}$, $\lhb\rho_{b_n}, \mathfrak{e}_{i+1}^+ \rhb
$ or $\lhb\mathfrak{e}_i\rhb\cup
\lhb a,\rho_{b_n} \rhb$. We conclude by similar
arguments.~%
\end{pf}

We now turn our attention to contractible loops. Let $\wp$ be a
contractible simple loop in $\cS_n$ made of edges. Then $\wp$ splits
$\cS_n$ into two domains. Only one of these is homeomorphic to a
disk.\footnote{This is a consequence of the Jordan--Sch\"{o}nflies
theorem, applied in the universal cover of $\cS_n$, which is either the
plane when $g=1$, or the unit disk when $g\ge2$; see, for example,~\cite{epstein66cmi},
Theorem~1.7.} We call it the \textit{inner
domain} of~$\wp$, and we call the other one the \textit{outer domain}
of $\wp$. In particular, these domains are well defined for loops whose
diameter is smaller than $\eps_0$, when $n$ is large enough.
\begin{lem}\label{1reg}
A.s., for all $\eps>0$, there exists $0 < \delta< \eps\wedge\eps_0$
such that for all $k$ sufficiently large, the inner domain of any
simple loop made of edges in $\cS_{n_k}$ with diameter less than
$\delta
$ has diameter less than $\eps$.
\end{lem}
\begin{pf}
We adapt the method used by Miermont in~\cite{miermont08sphericity}.
The idea is that a contractible loop separates a whole part of the map
from the base point. Then the labels in one of the two domains it
separates are larger than the labels on the loop. In the $g$-tree, this
corresponds to having a part with labels larger than the labels on the
``border.'' In the continuous limit, this creates an increase point for
both $\CC_\infty$ and $\Lab_\infty$.

Suppose that, with positive probability, there exists $0 < \eps<\eps
_0$ for which, along some (random) subsequence of the sequence
$(n_k)_{k\ge0}$, there exist contractible simple loops $\wp_n$ made of
edges in $\cS_n$ with diameter tending to $0$ (with respect to the
rescaled metric $\delta_{(n)}$) and whose inner domains are of diameter
larger than $\eps$. Let us reason on this event. First, notice
that,
because $g \ge1$, the outer domain of $\wp_n$ contains at least one
noncontractible loop, so that its diameter is larger than $\eps_0 >
\eps
$ by Lemma~\ref{eps0}.

Let $s^\bullet$ be the unique point where $\Lab_\infty$ reaches its
minimum, and $s_n^\bullet$ be an integer where $\Lab_n$ reaches its
minimum. We call $w^\bullet_n \de\rr\tr_n(s_n^\bullet)$ the
corresponding point in the $g$-tree. This is a vertex at $\delta
_{n}$-distance $1$ from $v^\bullet_n$. Let us take~$x_n$ in the domain
that does not contain $w_n^\bullet$, such that the distance
between~$x_n$ and $\wp_n$ is maximal. (If $w_n^\bullet\in\wp_n$, we
take $x_n$ in either of the two domains according to some convention.)
Let $y_n \in\wp_n \cap( \lhb\rho_{w^\bullet_n}, w^\bullet_n \rhb
\cup
\fl_n \cup\lhb\rho_{x_n},x_n \rhb)$ be such that there exists an
injective path\footnote{Depending on the case, the path $\mathfrak p_n$
will be of one of the following forms:
\begin{itemize}[$\diamond$]
\item[$\diamond$]$\lhb x_n,y_n \rhb$, with $y_n\in\lhb\rho_{x_n},x_n \rhb$;
\item[$\diamond$]$\lhb x_n,\rho_{x_n} \rhb\cup\lhb\rho_{x_n},y_n \rhb$, with
$y_n\in\fl_n$;
\item[$\diamond$]$\lhb x_n,\rho_{x_n} \rhb\cup\lhb\rho_{x_n},e_1^+\rhb\cup
\lhb e_2 \rhb\cup\cdots\cup\lhb e_k \rhb\cup\lhb e_k^+,y_n \rhb$
for some half-edges $e_1$, $e_2, \ldots, e_k$ of $\s_n$ satisfying
$e_i^+=e_{i+1}^-$, with $y_n\in\fl_n$;
\item[$\diamond$]$\lhb x_n,\rho_{x_n} \rhb\cup\lhb\rho_{x_n},e_1^+\rhb\cup
\lhb e_2 \rhb\cup\cdots\cup\lhb e_k \rhb\cup\lhb e_k^+,\rho
_{w_n^\bullet} \rhb\cup\lhb\rho_{w_n^\bullet}, y_n \rhb$ for some
half-edges $e_1, e_2, \ldots, e_k$ of $\s_n$ satisfying
$e_i^+=e_{i+1}^-$, with $y_n\in\lhb\rho_{w^\bullet_n}, w^\bullet_n
\rhb$.
\end{itemize}
} $\mathfrak p_n$ in $\tr_n$ from $x_n$ to~$y_n$ that intersects $\wp
_n$ only at $y_n$. In other words, when going from $x_n$ to $w^\bullet
_n$ along some injective path, $y_n$ is the first vertex belonging
to $\wp_n$ we meet; see Figure~\ref{xyzn}. Such a point exists because
$x_n$ and $w_n^\bullet$ do not belong to the same of the two components
delimited by $\wp_n$. Up to further extraction, we suppose that
$s_n^\bullet/ 2n \to s^\bullet$, $x_n \to x$ and $y_n \to y$. We call
$\mathfrak p \subseteq\lhb\rho_{w^\bullet}, w^\bullet\rhb\cup\fl
_\infty\cup\lhb\rho_{x},x \rhb$ the injective path corresponding
to $\mathfrak p_n$ in the limit, that is, the path defined
as $\mathfrak p_n$ ``without the subscripts $n$.'' Because the distance
between two points in the same domain as $x_n$ is smaller than $2
\delta
_{(n)}(x_n,\wp_n) + \diam(\wp_n)$, we obtain that $\delta
_{(n)}(x_n,y_n) \ge\eps/4$, as soon as $\diam(\wp_n) \le\eps/2$. In
particular, we see that $x \neq y$, and that the path $\mathfrak p$ is
not reduced to a single point.

%
%
\begin{figure}

\includegraphics{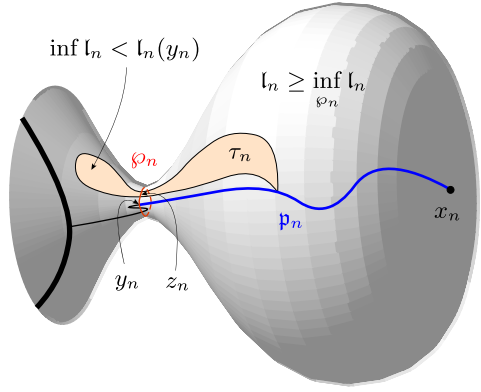}

\caption{The path $\wp_n$ intersects $\tau_n$. This figure represents
the case where $y_n\in\lhb\rho_{x_n},x_n \rhb$.}
\label{xyzn}
\end{figure}

Let us first suppose that $y \neq w^\bullet\de\TT_\infty(s^\bullet)$.
(In particular, $w_n^\bullet\notin\wp_n$ for $n$ large, so that there
is no ambiguity on which domain to chose $x_n$.) In that case, $y \in
(\lhb\rho_{w^\bullet}, w^\bullet\rhb\cup\fl_\infty\cup\lhb
\rho
_{x},x \rhb) \bs\{x,w^\bullet\}$, so that the points\vspace*{1pt} in $\TT_\infty
^{-1}(y)$ are increase points of $\CC_\infty$. By Lemma~\ref{pc}, we
can find a subtree\footnote{Here again, we need to distinguish between
some cases:
\begin{itemize}[$\diamond$]
\item[$\diamond$] if $y \in\lhb\rho_{x},x \rhb$, then $\mathfrak p = \lhb x, y
\rhb$, and $\tau$ is a tree to the left or right of $\lhb\rho_{x},x
\rhb$ rooted at some point in $\lhb{x},y \rhb\bs\{x,y\}$;
\item[$\diamond$] if $y \in\fl_\infty\bs\{\rho_{x}\}$, then $\tau$ is a tree
of $\TT_\infty$ rooted on $(\mathfrak p \cap\fl_\infty) \bs\{y\}$;
\item[$\diamond$] if $y \in\lhb\rho_{w^\bullet}, w^\bullet\rhb\bs\{\rho
_{w^\bullet}\}$, then $\tau$ is a tree to the left or right of $\lhb
\rho_{w^\bullet},y \rhb$.
\end{itemize}}
 $\tau$, not containing $y$, satisfying $\inf_{\tau} \Lab_\infty<
\Lab
_\infty(y)$ and rooted on the path $\mathfrak p$.

We consider a discrete approximation $\tau_n$ rooted on $\mathfrak
p_n$. Because the loop~$\wp_n$ is contractible, all the labels of the
points in the same domain as $x_n$ are larger than $\inf_{\wp_n}
\mathfrak{l}
_n$. Indeed, the labels represent the distances (up to an additive
constant) in $\mathfrak{q}_n$ to the base point, and every geodesic
path from
such a point to the base point has to intersect $\wp_n$. For $n$ large
enough, it holds that $\inf_{\tau_n} \mathfrak{l}_n < \inf_{\wp_n}
\mathfrak{l}_n $. As
a consequence, $\tau_n$ cannot entirely be included in the domain
containing $x_n$. Therefore, the set $\wp_n\cap\tau_n$ is not empty,
so that we can find $z_n \in\wp_n\cap\tau_n$. Up to extraction, we
may suppose that $z_n \to z$.

On one hand, $\delta_{(n)}(y_n,z_n) \le\diam(\wp_n)$, so that $y
\sim
_\infty z$. On the other hand, $z \in\tau$ and $y \notin\tau$, so
that $y \neq z$. Because $y$ is not a leaf, this contradicts
Theorem~\ref{ip}.

When $y = w^\bullet$, we use a different argument. Let $a_n=\rr\tr
_n(\alpha_n)$ and $b_n=\rr\tr_n(\beta_n)$ be, respectively, in the inner
and outer domains of $\wp_n$, such that their distance to $\wp_n$ is
maximal. Because $a_n$ and $b_n$ do not belong to the same domain, we
can find
\[
t_n^1 \in\overrightarrow{\lbracket\alpha_n,\beta_n \rbracket}
\quad\mbox{and}\quad t_n^2
\in
\overrightarrow{\lbracket\beta_n,\alpha_n \rbracket}
\]
such that $\rr\tr_n(t_n^1)$, $\rr\tr_n(t_n^2) \in\wp_n$. Up to
extraction, we suppose that
\[
\frac{\alpha_n}{2n} \to\alpha,\qquad\frac{\beta_n}{2n} \to\beta,\qquad
\frac{t_n^1}{2n} \to t^1 \in\overrightarrow{[\alpha,\beta]}
\quad\mbox{and}\quad
\frac{t_n^2}{2n} \to t^2 \in\overrightarrow{[\beta,\alpha]}.
\]
Because $\diam(\wp_n) \to0$, we have $\TT_\infty(t^1)=\TT_\infty
(t^2)=w^\bullet$. Moreover, the argument we used to prove that $x\neq
y$ yields that $\TT_\infty(\alpha) \neq w^\bullet$ and $\TT_\infty
(\beta
) \neq w^\bullet$. As a result, we obtain that $t^1 \neq t^2$. This
contradicts Lemma~\ref{min}.
\end{pf}

It remains to deal with general loops that are not necessarily made of
edges. We reason on the set of full probability where Lemmas~\ref{eps0}
and~\ref{1reg} hold. We fix $0 < \eps< \diam(\mathfrak{q}_\infty)
/4$. Let $\eps
_0$ be as in Lemma~\ref{eps0} and $\delta$ as in Lemma~\ref{1reg}. For
$k$ sufficiently large, the conclusions of both lemmas hold, together
with the inequality \mbox{$\delta \gamma n_k^{1/4} \ge12$}. Now, take
any loop $\loo$ drawn in~$\cS_{n_k}$ with diameter less than $\delta
/2$. Consider the union of the closed faces\footnote{We call \textit
{closed face} the closure of a face.} visited by $\loo$. The boundary
of this union consists in simple loops made of edges in~$\cS_{n_k}$.
Let us call $\Lambda$ the set of these simple loops.

Because\vspace*{-1pt} every face of $\cS_{n_k}$ has a diameter smaller than
$3/\gamma
n_k^{1/4}$, we see that for all $\lambda\in\Lambda$, $\diam
(\lambda
) \le\diam(\loo)+ 6/\gamma n_k^{1/4} \le\delta$. Then, by
Lemma~\ref{eps0},~$\lambda$ is contractible and, by Lemma~\ref{1reg},
its inner domain is of diameter less than~$\eps$. By definition, for
all $\lambda\in\Lambda$, $\loo$ entirely lies either inside the inner
domain of $\lambda$, or inside its outer domain. We claim that there
exists one loop in $\Lambda$ such that $\loo$ lies in its inner domain.
Then, it will be obvious that $\loo$ is homotopic to $0$ in its $\eps
$-neighborhood.

Let us suppose that $\loo$ lies in the outer domain of every loop
$\lambda\in\Lambda$. Then, every face of $\cS_{n_k}$ is either visited
by $\loo$, or included in the inner domain of some loop $\lambda\in
\Lambda$. As a result, we obtain that $\diam(\mathfrak{q}_\infty)
\le\diam(\loo
)+ 2 \sup_{\lambda\in\Lambda}\diam(\lambda) + 6/\gamma n_k^{1/4}
\le3\delta$. This is in contradiction with our choice of~$\delta$.


\section{Transfering results from the planar case through Chapuy's
bijection}\label{secprle}

In order to prove Lemmas~\ref{min},~\ref{pc},~\ref{leme2} and \ref
{leme4}, we rely on similar results for the Brownian snake driven by a
normalized excursion $(\mathbh{e},Z)$. This means that $\mathbh{e}$
has the law of a
normalized Brownian excursion, and, conditionally given $\mathbh{e}$, the
process $Z$ is a Gaussian process with covariance
\[
\cov(Z_x,Z_y)=\inf_{[x\wedge y,x\vee y]}\mathbh{e}.
\]
We first focus on the proofs of Lemmas~\ref{min} and~\ref{pc}.
Lemmas 3.1 and 3.2 in~\cite{legall08slb} state that, a.s., $Z$ reaches
its minimum at a unique point, and that, a.s., $\IP(\mathbh{e})$ and
$\IP(Z)$
are disjoint sets. We will use a bijection due to Chapuy~\cite{chapuy08sum}
to transfer these results to our case.

\subsection{Chapuy's bijection}

\begin{figure}

\includegraphics{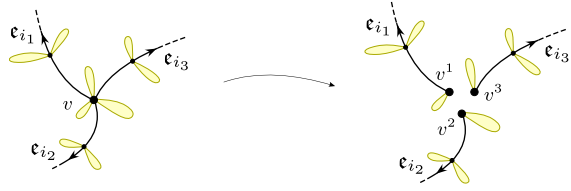}

\caption{Slicing an intertwined node $v$.}
\label{slice}
\end{figure}

Chapuy's bijection consists in ``opening'' $g$-trees into plane trees.
We briefly describe it here. See~\cite{chapuy08sum} for more details.
Let $\tr$ be a $g$-tree whose scheme $\s$ is dominant. Such a $g$-tree
will be called \textit{dominant} in the following. As usual, we arrange
the half-edges of $\s$ according to its facial order: $\mathfrak
{e}_1=\mathfrak{e}_*, \ldots
, \mathfrak{e}_\ks$. Let $v$ be one of the nodes of $\tr$. We can
see it as a
vertex of~$\s$. Let us call $\mathfrak{e}_{i_1}$, $\mathfrak
{e}_{i_2}$ and $\mathfrak{e}_{i_3}$ the
three half-edges starting from~$v$ (i.e., $v=\mathfrak
{e}_{i_1}^-=\mathfrak{e}_{i_2}^-=\mathfrak{e}
_{i_3}^-$), where $i_1 < i_2 < i_3$. We say that $v$ is \textit
{intertwined} if the half-edges $\mathfrak{e}_{i_1}$, $\mathfrak
{e}_{i_2}$, $\mathfrak{e}_{i_3}$ are
arranged according to the counterclockwise order around $v$ (see
Figure~\ref{slice}). When $v$ is intertwined, we may \textit{slice} it:
we define a new map, denoted by $\tr\bbslash v$, by slicing the
node $v$ into three new vertices $v^1$, $v^2$ and $v^3$ (see
Figure \ref
{slice}).

%
%

The map obtained by such an operation turns out to be a dominant
$(g-1)$-tree. After repeating $g$ times this operation, we are left
with a plane tree. In that regard, we call \textit{opening sequence}
of $\tr$ a $g$-uple $(v_1,\ldots, v_g)$ such that $v_g$ is an
intertwined node of $\tr$, and for all $1 \le i \le g-1$, the
vertex $v_i$ is an intertwined node of $\tr\bbslash v_g \bbslash\cdots
\bbslash v_{i+1}$. We can show that every $g$-tree has exactly $2g$
intertwined nodes, and thus $2^g g!$ opening sequences.

To reverse the slicing operation, we have to intertwine and glue back
the three vertices together. We then need to record which vertices are
to be glued together. This motivates the following definition: we call
\textit{tree with $g$ triples} a pair $(\mathbf{t}, (c_1, \ldots,
c_g))$, where:
\begin{itemize}[$\diamond$]
\item[$\diamond$]$\mathbf{t}$ is a (rooted) plane tree;
\item[$\diamond$] for $1 \le i \le g$, $c_i = \{v_i^1, v_i^2, v_i^3 \} \subseteq
V(\mathbf{t})$ is a set of three vertices of $\mathbf{t}$;
\item[$\diamond$] the vertices $v_i^j$, $1 \le i \le g$, $1\le j \le3$, are
pairwise distinct;
\item[$\diamond$] the vertices of the tree
\[
\bigcup_{i,i',j,j'} \biglhb v_i^j, v_{i'}^{j'} \bigrhb
\]
have degree at most $3$, and the $v_i^j$'s have degree exactly $1$ in
that tree. (As in the case of $g$-trees, the set $\lhb a, b\rhb$
represents the range of the unique path linking~$a$ and~$b$ in the tree.)
\end{itemize}
Let $\tr$ be a $g$-tree together with an opening sequence $(v_1,\ldots,
v_g)$. For all $1 \le i \le g$, let us call $c_i$ the triple of
vertices obtained from the slicing of $v_i$, as well as $\mathbf{t}\de
\tr
\bbslash v_g \bbslash\cdots\bbslash v_{1}$ the plane tree. We define
$\Phi(\tr,(v_1,\ldots, v_g)) \de(\mathbf{t}, (c_1, \ldots, c_g))$. Then
$\Phi$
is a bijection from the set of all dominant $g$-tree equipped with an
opening sequence into the set of all trees with $g$ triples.

Now, when the $g$-tree is well-labeled, we can do the same slicing
operation, and the three vertices we obtain all have the same label. We
call \textit{well-labeled tree with $g$ triples} a tree with $g$
triples $(\mathbf{t}, (c_1, \ldots, c_g))$ carrying a labeling function
$\llll
\dvtx V(\mathbf{t}) \to\Z$ such that:
\begin{itemize}[$\diamond$]
\item[$\diamond$]$\llll(e^-) = 0$, where $e$ is the root of $\mathbf{t}$;
\item[$\diamond$] for every pair of neighboring vertices $v\sim v'$, we have $\llll
(v) - \llll(v') \in\{-1,0,1\}$;
\item[$\diamond$] for all $1\le i \le g$, we have $\llll(v_i^1)=\llll
(v_i^2)=\llll(v_i^3)$.
\end{itemize}
We call $\W_n$ the set of all well-labeled trees with $g$ triples
having $n$ edges. The bijection $\Phi$ then extends to a bijection
between dominant well-labeled $g$-trees equipped with an opening
sequence and well-labeled trees with $g$ triples.

\subsection{Contour pair of an opened $g$-tree}\label{secop}

The contour pair of an opened $g$-tree can be obtained from the contour
pair of the $g$-tree itself (and vice versa). The labeling function is
basically the same, but read in a different order. The contour function
is slightly harder to recover, because half of the forests are to be
read with the floor directed ``upward'' instead of ``downward.''
Because we will deal at the same time with $g$-trees and plane trees in
this section, we will use a Gothic font for objects related to
$g$-trees, and a~boldface font for objects related to plane trees. In
the following, we use the notation of Section~\ref{decomp}.

Let $(\tr,\mathfrak{l})$ be a well-labeled dominant $g$-tree with
scheme $\s$
and $(\mathbf{t},\llll)$ be one of the $2^g g!$ corresponding opened
well-labeled trees. The intertwined nodes of the $g$-tree correspond to
intertwined nodes of its scheme, so that the opening sequence used to
open $(\tr,\mathfrak{l})$ into $(\mathbf{t},\llll)$ naturally
corresponds to an
opening sequence of $\s$. Let $\sss$ be the tree obtained by
opening $\s
$ along this opening sequence. We identify the half-edges of $\s$ with
the half-edges of $\sss$, and arrange them according to the facial
order of $\sss$: $\eee_1=\mathfrak{e}_*, \eee_2, \ldots, \eee
_\ks$. (Beware
that this is not the usual arrangement according to the facial order
of $\s$.) Now, the plane tree $\mathbf{t}$ is obtained by replacing every
half-edge $\eee$ of $\sss$ with the corresponding forest $\mathfrak
{f}^{\eee}$ of
Proposition~\ref{decwl}, as in Section~\ref{secdec}.

We call $(C^\mathfrak{e}, L^\mathfrak{e})$ the contour pair of
$(\mathfrak{f}^\mathfrak{e},\mathfrak{l}^\mathfrak{e})$, we let
$\CC^\mathfrak{e}\de C^\mathfrak{e}- \sigma^\mathfrak{e}$ and we
define $\Lab^\mathfrak{e}$ by~(\ref{lele}).
For any edge $\{\eee_i,\eee_j\}\neq\{\mathfrak{e}_*,\bar\mathfrak
{e}_*\}$ with $i < j$, we
will visit the forest $\mathfrak{f}^{\eee_i}$ while ``going up'' and\vadjust{\goodbreak}
the forest
$\mathfrak{f}^{\eee_j}$ while ``coming down'' when we follow the
contour of $\mathbf{t}
$. Precisely, we define
%
%
\begin{equation}\label{op1}
\CCC^{\eee_i} \de\CC^{\eee_i} - 2 \underline\CC^{\eee_i}
\quad\mbox{and}\quad
\CCC^{\eee_j} \de\CC^{\eee_j}.
\end{equation}
The first function is the concatenation of the contour functions of the
trees in $\mathfrak{f}^{\eee_i}$ with an extra ``up step'' between every
consecutive trees. The second one is the concatenation of the contour
functions of the trees in $\mathfrak{f}^{\eee_j}$ with an extra
``down step''
between every consecutive trees. It is merely the contour function
of $\mathfrak{f}^{\eee_j}$ shifted in order to start at $0$. What
%
%
\begin{figure}

\includegraphics{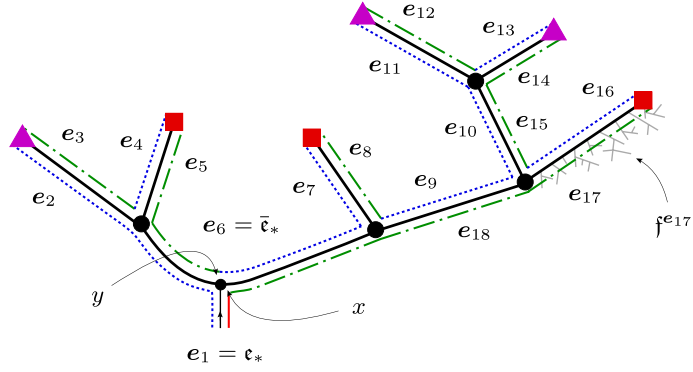}

\caption{Opening of a $2$-tree. The squares form one triple and the
triangles the other one. The (blue) short dashes correspond to the
upward-directed floors and the (green) long dashes to the
downward-directed floors. The (red) solid line on the right of the root
corresponds to the part of the tree containing the root that has to be
visited at the end. The forest $\mathfrak{f}^{\eee_{17}}$ is also
represented on
this figure.}
\label{opened}
\end{figure}
happens to the
forests $\mathfrak{f}^{\mathfrak{e}_*}$ and $\mathfrak{f}^{\bar
\mathfrak{e}_*}$ is a little more intricate. Let
us first call (see Figure~\ref{opened})
%
%
\begin{eqnarray}\label{op2}
x &\de& \inf\{s \dvtx \CC^{\mathfrak{e}_*}(s)=\underline\CC
^{\mathfrak{e}_*}(u )\},
\nonumber\\[-8pt]\\[-8pt]
y &\de& \inf\{s \dvtx \CC^{\bar\mathfrak{e}_*}(s)= - \sigma
^{\mathfrak{e}_*} - \underline
\CC^{\mathfrak{e}_*}(u )\}.
\nonumber
\end{eqnarray}
When visiting the forest $\mathfrak{f}^{\bar\mathfrak{e}_*}$, the
floor is directed
downward up to time $y$ and then upward:
%
%
\begin{equation}\label{op3}\qquad
\CCC^{\bar\mathfrak{e}_*} \de( \CC^{\bar\mathfrak{e}_*}
(s) )_{0 \le s \le y}
\bullet\Bigl( \CC^{\bar\mathfrak{e}_*} (y + s )- 2\inf
_{[y,y + s]} \CC
^{\bar\mathfrak{e}_*} +\CC^{\bar\mathfrak{e}_*}(y )
\Bigr)_{0 \le s \le m^{\bar\mathfrak{e}_*} - y}.
\end{equation}
Finally, the forest $\mathfrak{f}^{\mathfrak{e}_*}$ is visited twice.
The first time (when
beginning the contour), it is visited between times $u$
and $m^{\mathfrak{e}_*}$,
and the floor is directed upward:
%
%
\begin{equation}\label{op4}
\CCC^{\mathfrak{e}_*,1} \de\Bigl( \CC^{\mathfrak{e}_*} (u + s
)- 2\inf_{[u,u + s]}
\CC^{\mathfrak{e}_*} +\CC^{\mathfrak{e}_*}(u )\Bigr)_{0 \le
s \le m^{\mathfrak{e}_*} - u}.
\end{equation}
The second time (when finishing the contour), we visit it between
times $0$ and $x$ with the floor directed downward, then we visit a
part of the tree containing the root between times $x$ and $u$:
%
%
\begin{equation}\label{op5}\qquad
\CCC^{\mathfrak{e}_*,2} \de( \CC^{\mathfrak{e}_*} (s)
)_{0 \le s \le x} \bullet
\Bigl( \CC^{\mathfrak{e}_*} (x + s )- 2\inf_{[x+s,u]} \CC
^{\mathfrak{e}_*} +
\underline{\CC}^{\mathfrak{e}_*}(u )\Bigr)_{0 \le s \le u - x}.
\end{equation}
%
The contour pair of $(\mathbf{t},\llll)$ is then given by
%
%
\begin{equation}\label{contopen}
\cases{
\CCC\de\CCC^{\eee_1,1} \bullet\CCC
^{\eee_2} \bullet\CCC^{\eee_3} \bullet\cdots\bullet\CCC^{\eee
_\ks}
\bullet\CCC^{\eee_1,2},\cr
\LLL\de\Lab^{\eee_1,1} \bullet\Lab
^{\eee_2} \bullet\Lab^{\eee_3} \bullet\cdots\bullet\Lab^{\eee
_\ks}
\bullet\Lab^{\eee_1,2},}
\end{equation}
where
\[
\Lab^{\eee_1,1} \de\bigl(\Lab^{\eee_1}(u+s) - \Lab^{\eee
_1}(u)\bigr)_{0
\le s \le m^{\eee_1} - u} \quad\mbox{and}\quad
\Lab^{\eee_1,2} \de(\Lab^{\eee_1}(s) )_{0 \le s \le u}.
\]

\subsection{Opened uniform well-labeled $g$-tree}

As in Section~\ref{secsl}, we let $(\tr_n,\mathfrak{l}_n)$ be uniformly
distributed over the set $\T_n$ of well-labeled $g$-trees with $n$
edges, and, applying Skorokhod's representation theorem, we assume that
the convergence of Proposition~\ref{cvint} holds almost surely. Let
$(\ii_n)_{n\in\N}$ be a sequence of i.i.d. random variables uniformly
distributed over $\lbracket1,2^g g! \rbracket$ and independent of
$(\tr
_n,\mathfrak{l}
_n)_{n\in\N}$. With any dominant scheme $\s\in\Sg_*$ and integer
$\ii
\in\lbracket1,2^g g! \rbracket$, we associate a deterministic
opening sequence.
When $(\tr_n,\mathfrak{l}_n)$ is dominant, we may then define $(\mathbf{t}_n,\llll
_n)$ as the opened tree of $(\tr_n,\mathfrak{l}_n)$ according to the opening
sequence determined by the integer~$\ii_n$. In this case, we call
$(\CCC
_n,\LLL_n)$ the contour pair of $(\mathbf{t}_n,\llll_n)$. When $(\tr
_n,\mathfrak{l}
_n)$ is not dominant, we simply set $(\CCC_n,\LLL_n)=(\oo_{2n},\oo
_{2n})$, where we write $\oo_{\zeta}\dvtx x \in[ 0, \zeta] \mapsto
0$. We
also let
\[
\CCC_{(n)} \de\biggl(\frac{\CCC_n(2nt)} {\sqrt{2n}} \biggr)_{0\le t \le
1}\quad\mbox{and}\quad
\LLL_{(n)} \de\biggl(\frac{\LLL_n(2nt)} {\g}\biggr)_{0\le t \le1}
\]
be the rescaled versions of $\CCC_n$ and $\LLL_n$.

We now work at fixed $\omega$ for which Proposition~\ref{cvint} holds,
$\s_\infty\in\Sg_*$, and such that for all $\ii\in\lbracket1,2^g
g! \rbracket$,
$|\{n\in\N \dvtx \ii_n=\ii\} | = \infty$. Note that the set of such
$\omega$'s is of full probability. For $n$ large enough, $\s_n=\s
_\infty
\in\Sg_*$, so that $(\mathbf{t}_n,\llll_n)$ is well defined. For
all $n$
such that $\s_n=\s_\infty$ and $\ii_n=\ii$, we always open the $g$-tree
$(\tr_n,\mathfrak{l}_n)$ according to the same opening sequence, so
that the
ordering $\eee_1, \eee_2, \ldots, \eee_\ks$ of the half-edges
of $\s
_n$ is always the same. As a result, we obtain that
\[
\bigl(\CCC_{(n)},\LLL_{(n)}\bigr) \mathop{\hbox to 1cm{\rightarrowfill}}_{n \to\infty}^{n\dvtx \ii
_n=\ii} (\CCC
_\infty^\ii,\LLL_\infty^\ii),
\]
where $(\CCC_\infty^\ii,\LLL_\infty^\ii)$ is defined by~(\ref{op1})--(\ref{op5}) and~(\ref{contopen})
when replacing every occurrence of $\CC^\mathfrak{e}$ by $\CC_\infty
^\mathfrak{e}\de
C_\infty^\mathfrak{e}- \sigma_\infty^\mathfrak{e}$ and every
occurrence of $\Lab^\mathfrak{e}$ by
$\Lab_\infty^\mathfrak{e}$. Note that $(\CCC_{(n)},\LLL_{(n)})$
has exactly $2^g
g!$ a priori distinct accumulation points, each corresponding to one of
the ways of opening the real $g$-tree $\TT_\infty$.

Now, because every $\Lab_\infty^\mathfrak{e}$ goes from $0$ to $0$,
it is easy to
see that for all $\ii$, the points where $\Lab_\infty$ reaches its
minimum are in one-to-one correspondence with the points where $\LLL
_\infty^\ii$ reaches its minimum. Moreover, we can see that if $\CC
_\infty$ and $\Lab_\infty$ have a common increase point, then at least\vadjust{\goodbreak}
one of the pairs $(\CCC_\infty^\ii,\LLL_\infty^\ii)$ will also
have a
common increase point. Indeed, let us suppose that $\CC_\infty$
and $\Lab_\infty$ have a common increase point. Then, there exists
$\mathfrak{e}
\in\vec E(\s_\infty)$ such that $\CC_\infty^\mathfrak{e}$ and
$\Lab_\infty^\mathfrak{e}$
have a common increase point $s\in[0, m^\mathfrak{e}_\infty]$. We
use the
following lemma:
\begin{lem}\label{iptransfert}
Let $f\dvtx[0,m] \to\R$ be a function.
\begin{itemize}[$\diamond$]
\item[$\diamond$] If $s \in[0,m)$ is an increase point of $f$, then $s$ is an
increase point of $f-2\underline{f}$ as well.
\item[$\diamond$] If $s \in(0,m]$ is an increase point of $f$, then $s$ is an
increase point of $r\mapsto f(r)-2\inf_{[r,m]}{f}$.
\end{itemize}
\end{lem}

We postpone the proof of this lemma and finish our argument. If $s <
m^\mathfrak{e}_\infty$, then $s$ is a common increase point of $\CCC
_\infty^\mathfrak{e}$
and $\LLL_\infty^\mathfrak{e}$ thanks to Lemma~\ref{iptransfert}.
When $\mathfrak{e}=\mathfrak{e}_*$,
this fact remains true if we define $\CCC_\infty^\mathfrak{e}\de
\CCC_\infty^{\mathfrak{e}
,2} \bullet\CCC_\infty^{\mathfrak{e},1}$. Note that~$x$ is an
increase point of
$\CCC_\infty^\mathfrak{e}$, even if $0$ is not an increase point of
the second
function defining $\CCC_\infty^{\mathfrak{e},2}$ in~(\ref{op5}). In
this case,
for all $\ii$, $\CCC_\infty^\ii$ and $\LLL_\infty^\ii$ have a~common
increase point.

Let us now suppose that $s = m^\mathfrak{e}_\infty$, and let us fix
$\ii\in\lbracket1,2^g g! \rbracket$. We consider the opening
corresponding to $\ii$. If $\eee
_i=\mathfrak{e}$ is visited while coming down in the contour of the
opened tree,
then we conclude as above. If both~$\eee_i$ and~$\eee_{i+1}$ are
visited while going up, then $0$ will be an increase point of $\CCC
_\infty^{\eee_{i+1}}$, so that
$\CCC_\infty^\ii$ and $\LLL_\infty^\ii$ will still have a common
increase point. In the remaining case where $\eee_i$ is visited while
going up and $\eee_{i+1}$ is visited while coming down (i.e., $\eee
_{i+1}=\bar{\eee_{i}}$), we cannot conclude that $\CCC_\infty^\ii$ and
$\LLL_\infty^\ii$ have a common increase point. This, however, only
happens when the node $\mathfrak{e}^+$ belongs to the opening
sequence. But when
we pick an opening sequence, we can always choose not to pick a given
node, because at each stage of the process, we have at least $2$
intertwined nodes. This implies that at least one of the opening
sequences will not contain $\mathfrak{e}^+$, and the corresponding
pair $(\CCC
_\infty^\ii, \LLL_\infty^\ii)$ will have a common increase point.
\begin{pf*}{Proof of Lemma~\ref{iptransfert}}
Let $s \in[0,m)$ be an increase point of $f$. If $s$ is a
right-increase point of $f$, then $f(r) \ge f(s)$ when $s \le r \le t$
for some $t > s$. For such $r$'s, $\underline f(r) = \underline f(s)$,
so that $f(r) - 2 \underline f(r) \ge f(s) - 2 \underline f(s)$,
and $s$ is a~right-increase point of $f-2 \underline f$.

If $s$ is a left-increase point of $f$, then $f(r) \ge f(s)$ when $t
\le r \le s$ for some \mbox{$t < s$}. If $f(s) > \underline f(s)$, then, using
the fact that $\underline f(s) = \underline f(r) \wedge\inf_{[r,s]}
f$, we obtain that $\underline f(r) = \underline f(s)$ when $t \le r
\le s$ and conclude as above that $s$ is a left-increase point of $f-2
\underline f$. Finally, if $f(s) = \underline f(s)$, then for all $r
\ge s$, we have $f(r) - 2 \underline f(r) = (f(r) - \underline
f(r)) - \underline f(r) \ge0 - \underline f(s) = f(s) - 2
\underline f(s)$, and because $s <m$, we conclude that $s$ is a
right-increase point of $f-2 \underline f$.

We obtain the second assertion of the lemma by applying the first one
to $m-s$ and the function $x \mapsto f(m-x)$.
\end{pf*}

\subsection{Uniform well-labeled tree with $g$ triples}

Conditionally on the event $D_n \de\{ (\CCC_n,\LLL_n)\neq(\oo
_{2n},\oo
_{2n})\}$, the distribution of $(\CCC_n,\LLL_n)$ is that of the contour
pair of a uniform well-labeled tree with $g$ triples. We use this fact
to see that the law of $(\CCC_{(n)},\LLL_{(n)})$ converges weakly
toward a law absolutely continuous with respect to the law of
$(\mathbh{e}
,Z)$. Let $(\tau_n, \lambda_n)$ be uniformly distributed over the
set $\mathcal{T}^0_n$ of all well-labeled plane trees with $n$ edges.
We call
$(\Gamma_n, \Lambda_n)$ the contour pair of $(\tau_n,\lambda_n)$ and
define as usual the rescaled versions of both functions,
%
%
\begin{equation}\label{scuncp}
\Gamma_{(n)} \de\biggl(\frac{\Gamma_n(2n t)} {\sqrt{2n}} \biggr)_{0 \le t
\le1} \quad\mbox{and}\quad\Lambda_{(n)} \de\biggl(\frac{\Lambda_n(2nt)} {\g}
\biggr)_{0
\le t \le1}.
\end{equation}

For all $n \ge1$, $k \in\Z$ and $x \in\R$, we define
\[
X_n(k) \de|\{ v \in\tau_n \dvtx \lambda_n(v)=k
\}
| \quad\mbox{and}\quad X_{(n)}(x) \de\frac1 n \g X_n(
\lfloor
\g x \rfloor),
\]
respectively, the profile and rescaled profile of $(\tau_n,\lambda_n)$.
We let $\I$ be the one-dimensional ISE (random) measure defined by
\[
\langle\I,h \rangle\de\int_0^1 dt\,  h(Z_t)
\]
for every nonnegative measurable function $h$. By
\cite{bousquetmelou06dia}, Theorem 2.1, it is known that $\I$ a.s. has a
continuous density $\fISE$ with compact support. In other words,
$\langle\I,h \rangle= \int_\R dx\,  h(x) \fISE(x)$ for every
nonnegative measurable function $h$.
\begin{prop}\label{jd}
The triple $(\Gamma_{(n)}, \Lambda_{(n)}, X_{(n)} )$
converges weakly toward the triple $(\mathbh{e}, Z, \fISE)$ in the
space $\C
([0, 1],\R)^2\times\mathcal{C}_c(\R)$ endowed with the product topology.
\end{prop}
\begin{pf}
It is known that the pair $(\Gamma_{(n)}, \Lambda_{(n)} )$
converges weakly to $(\mathbh{e}, Z)$: in~\cite{chassaing04rpl}, Theorem 5,
Chassaing and Schaeffer proved this fact with $\lfloor2nt \rfloor$
instead of $2nt$ in the definition~(\ref{scuncp}). The claim as stated
here easily follows by using the uniform continuity of $(\mathbh{e}, Z)$.
Using~\cite{bousquetmelou06dia}, Theorem 3.6, and the fact that $\fISE$
is a.s. uniformly continuous~\cite{bousquetmelou06dia}, Theorem 2.1, we
also obtain that the sequence $X_{(n)}$ converges weakly to $\fISE$. As
a result, the sequences of the laws of the processes $\Gamma_{(n)}$,
$\Lambda_{(n)}$ and $X_{(n)}$ are tight. The sequence $(\nu_n)$ of the
laws of $(\Gamma_{(n)}, \Lambda_{(n)}, X_{(n)} )$ is then
tight as well, and, by Prokhorov's lemma, the set $\{\nu_n, n \ge0\}$
is relatively compact. Let $\nu$ be an accumulation point of the
sequence $(\nu_n)$. There exists a subsequence along which $
(\Gamma
_{(n)}, \Lambda_{(n)}, X_{(n)} )$ converges weakly toward a
random variable $(\mathbh{e}',Z',f')$ with law $\nu$. Thanks to Skorokhod's
theorem, we may and will assume that this convergence holds almost
surely along this subsequence. We know that
\[
(\mathbh{e}',Z') \stackrel{(\mathrm{d})}{=}(\mathbh{e}, Z) \quad\mbox{and}\quad f' \stackrel{(\mathrm{d})}{=}\fISE.
\]
It remains to see that $f'$ is the density of the occupation measure
of $Z'$, that is,
%
%
\begin{equation}
\label{dom}
\int_0^1 dt\,  h(Z_t') = \int_\R dx  \, h(x) f'(x)
\end{equation}
for all $h$ continuous with compact support. First, notice that
\begin{eqnarray*}
\frac1 n \sum_{k \in\Z} X_n(k) h (\gi k )&=&\frac1 n \int
_{\R}
dx  \,X_n (\lf x \rf)h (\gi\lf x \rf)\\[-2pt]
&=&\int_{\R} dx  \,X_{(n)} (x) h (\gi\lf\g x \rf)\\[-2pt]
&\to&\int_\R dx  \,f'(x) h(x)
\end{eqnarray*}
by dominated convergence, a.s. as $n \to\infty$ along the subsequence
we consider. It is convenient to introduce now the notation $\la s \ra
_n$ defined as follows: for $s \in[ 0, {2n})$, we set
\[
\la s \ra_n \de
\cases{
\lc s \rc, &\quad if $ \Gamma_n(\lc s \rc) - \Gamma_n(\lf s \rf) =
1$,\cr
\lf s \rf, &\quad if $ \Gamma_n(\lc s \rc) - \Gamma_n(\lf s \rf) = -1$.}
\]
Then, if we denote by $\tau_n(i)$ the $i$th vertex of the facial
sequence of $\tau_n$, and by $\rho_n$ the root of $\tau_n$, we obtain
that the time the process $( \tau_n(\la s\ra_n) )_{s \in[0,
2n)}$ spends at each vertex $v \in\tau_n\bs\{\rho_n\}$ is
exactly $2$. So we have
\begin{eqnarray*}
&&
\frac1 n \sum_{k \in\Z} X_n(k) h (\gi k )\\[-2pt]
&&\qquad= \frac1 n \sum_{v
\in\tau_n\bs\{\rho_n \}} h (\gi\lambda_n(v) )+ \frac1 n
h(0)\\[-2pt]
&&\qquad= \frac1 {2n} \int_0^{2n} ds \, h (\gi\Lambda_n (\la s \ra_n
))+ \frac1 n h(0)\\[-2pt]
&&\qquad= \int_0^1 ds \, h (\gi\Lambda_n(\la2ns \ra_n ))+
\frac
1 n h(0)\\[-2pt]
&&\qquad\to\int_0^1 dt \, h(Z_t')
\end{eqnarray*}
a.s. along the subsequence considered. We used the fact that
\[
\gi\Lambda_n(\la2ns \ra_n )\to Z'_s,
\]
which is obtained by using the uniform continuity of $Z'$.

This proves that $(\mathbh{e}',Z',f')$ has the same law as $(\mathbh
{e}, Z, \fISE)$.
Thus the only accumulation point $\nu$ of the sequence $(\nu_n)$ is the
the law of the process $(\mathbh{e}, Z, \fISE)$. By relative
compactness of
the set $\{\nu_n, n \ge0\}$, we obtain the weak convergence of the
sequence $(\nu_n)$ toward $\nu$.\vadjust{\goodbreak}
\end{pf}

We define
\[
W \de\frac{(\int\fISE^3 )^g} {\mathbb{E}[ (\int\fISE^3 )^g ]}.
\]
This quantity is well defined~\cite{chapuy08sum}, Lemma 10. We also
define the law of the pair $(\CCC_\infty, \LLL_\infty)$ by the
following formula: for every bounded Borel function~$\varphi$ on $\C
([0, 1],\R)^2$,
%
%
\begin{equation}\label{defw}
\mathbb{E}[ \varphi(\CCC_\infty, \LLL_\infty) ] = \mathbb{E}[ W
\varphi (\mathbh{e},Z) ].
\end{equation}

\begin{prop}\label{cs2}
The pair\vspace*{1pt} $(\CCC_{(n)}, \LLL_{(n)} )$ converges weakly toward
the pair $(\CCC_\infty, \LLL_\infty)$ in the space $( \C([0,
1],\R
)^2, \| \cdot\|_\infty)$ of pair of continuous real-valued
functions on $[0, 1]$ endowed with the uniform topology.
\end{prop}
\begin{pf}
Let $f$ be a bounded continuous function on $\C([0, 1],\R)^2$. We have
\begin{eqnarray*}
\mathbb{E}\bigl[ f\bigl( \CCC_{(n)},\LLL_{(n)} \bigr) \bigr] &=& \Pb(D_n)
\mathop{\sum_{(\tau,\lambda) \in\mathcal{T}^0_n}}_{(\tau
,\lambda) \leftrightarrow(\CCC,\LLL )} f(\CCC,\LLL)  \Pb\bigl((\tau
_n,\lambda_n)
=(\tau
,\lambda) | D_n \bigr)\\
&&{} + \Pb(\overline{D}_n )f( \oo_{2n},\oo_{2n} ),
\end{eqnarray*}
where we used the notation $(\tau,\lambda) \leftrightarrow(\CCC
,\LLL)$ to
mean that the well-labeled tree $(\tau,\lambda)$ is coded by the
contour pair $(\CCC,\LLL)$. It was shown in~\cite{chapuy08sum}, Lemma 8,
that the number of well-labeled trees with $g$ triples having $n$ edges
is equivalent to the number of well-labeled plane trees having $n$
edges, together with $g$ triples of vertices (not necessarily distinct
and not arranged) such that all the vertices of the same triple have
the same label. More precisely, we have
\[
\Pb\bigl((\tau_n,\lambda_n) = (\tau,\lambda) | D_n \bigr)
= \frac1 {|\W_n|} {\biggl(\sum_{k \in\Z} \lt\{v \in\tau \dvtx
\lambda(v) = k \}\rt^3 \biggr)^g} + O (n^{- 1/4 } ).
\]
And, because $f$ is bounded and $\Pb(D_n) \to1$, we obtain that
\[
\mathbb{E}\bigl[ f\bigl( \CCC_{(n)},\LLL_{(n)} \bigr) \bigr] \sim\frac{|\mathcal
{T}^0_n|} {|\W_n|}
\mathbb{E}\biggl[ \biggl(\sum_{k \in\Z} X_n(k)^3 \biggr)^g f \bigl(\Gamma_{(n)}, \Lambda
_{(n)} \bigr) \biggr].
\]
Using the asymptotic formulas $|\mathcal{T}^0_n| \sim\sqrt\pi
12^n n^{-
3/2}$, as well as $|\W_n| \sim c_g 12^n\times n^{(5g-3)/2}$ for some
positive constant $c_g$ only depending on $g$
(\cite{chapuy08sum}, Lemma 8), as well as the computation
\[
n^{-5/2} \sum_{k \in\Z} X_n(k)^3 = n^{-5/2} \int_{\R} dx \, X_n (
\lf
x \rf)^3= \gamma^{-2} \int_\R dx\,  X_{(n)}(x)^3,
\]
we see that there exists a positive constant $c$ such that
\[
\mathbb{E}\bigl[ f\bigl( \CCC_{(n)},\LLL_{(n)} \bigr) \bigr]
\sim c \mathbb{E}\biggl[ \biggl(\int
_\R dx \,  X_{(n)}(x)^3\biggr)^g f \bigl(\Gamma_{(n)}, \Lambda_{(n)} \bigr) \biggr].\vadjust{\goodbreak}
\]

Now, let $\eps> 0$. Thanks to~\cite{chapuy08sum},\vspace*{1pt} Lemma 10, we see
that both quantities $\mathbb{E}[ (\int\fISE^3 )^g ]$ and $\sup_{n}
\mathbb{E}[ ( \int X_{(n)}^3 )^{g+1} ]$ are finite. Then, using the
fact that
\[
\mathbb{E}\biggl[ \biggl(\int X_{(n)}^3 \biggr)^g \mathbh{1}_{\{\int X_{(n)}^3 > L\}} \biggr]
\le\frac1 L  \mathbb{E}\biggl[ \biggl(\int X_{(n)}^3 \biggr)^{g+1} \biggr],
\]
we obtain that, for $L$ sufficiently large,
\[
\sup_n \mathbb{E}\biggl[ \biggl(\int_\R dx \, X_{(n)}(x)^3 \biggr)^g f \bigl(\Gamma_{(n)},
\Lambda_{(n)} \bigr) \mathbh{1}_{\{\int X_{(n)}^3 > L\}} \biggr] <
\eps
\]
and
\[
\mathbb{E}\biggl[ \biggl(\int\fISE^3 \biggr)^g f (\mathbh{e},Z )\mathbh{1}_{\{\int
\fISE^3 > L\}} \biggr] < \eps.
\]
Thanks to the Proposition~\ref{jd}, for $n$ sufficiently large,
\begin{eqnarray*}
&&\biggl| \mathbb{E}\biggl[ \biggl(\int_\R dx  X_{(n)}(x)^3 \biggr)^g f \bigl(\Gamma_{(n)},
\Lambda_{(n)} \bigr) \mathbh{1}_{\{\int X_{(n)}^3 \le L\}} \biggr]\\
&&\hspace*{48.7pt}{} - \mathbb
{E}\biggl[ \biggl(\int \fISE^3 \biggr)^g f (\mathbh{e},Z )\mathbh{1}_{\{\int\fISE^3
\le L\}} \biggr] \biggr| <
\eps.
\end{eqnarray*}
This yields the existence of a constant $C$ such that
\[
\mathbb{E}\bigl[ f \bigl(\CCC_{(n)}, \LLL_{(n)} \bigr) \bigr] \ton C  \mathbb{E}\biggl[
\biggl(\int \fISE ^3 \biggr)^g f (\mathbh{e},Z ) \biggr],
\]
and we compute the value of $C$ by taking $f \equiv1$.
\end{pf}

Thanks to~(\ref{defw}), we see that the properties that hold almost
surely for the pair $(\mathbh{e},Z)$ also hold almost surely for
$(\CCC_\infty
,\LLL_\infty)$. We may now conclude thanks to
\cite{legall08slb}, Lemma 3.1, that
\begin{eqnarray*}
&&\Pb\bigl( \exists s \neq t \dvtx \Lab_\infty(s)=\Lab_\infty
(t)=\min
\Lab_\infty\bigr)\\
&&\qquad\le\frac1 {2^g g!} \sum_{\ii=1}^{2^g g!} \Pb\bigl( \exists s
\neq t
\dvtx
\LLL^\ii_\infty(s)=\LLL^\ii_\infty(t)=\min\LLL^\ii_\infty
\bigr)\\
&&\qquad= \Pb\bigl( \exists s \neq t \dvtx \LLL_\infty(s)=\LLL_\infty
(t)=\min\LLL_\infty\bigr) =0,
\end{eqnarray*}
and, by~\cite{legall08slb}, Lemma 3.2,
\begin{eqnarray*}
\Pb \bigl(\IP(\CC_\infty) \cap\IP(\Lab_\infty) \neq\varnothing
\bigr) &\le&
\sum_{\ii=1}^{2^g g!}
\Pb\bigl( \IP(\CCC_\infty^\ii) \cap\IP(\LLL_\infty^\ii) \neq
\varnothing
\bigr)\\
&=& 2^g g!  \Pb\bigl( \IP(\CCC_\infty) \cap\IP(\LLL_\infty)
\neq
\varnothing\bigr)\\
&=&0.
\end{eqnarray*}

This concludes the proof of Lemmas~\ref{min} and~\ref{pc}.

\subsection{Remaining proofs}

\subsubsection{\texorpdfstring{Proof of Lemma \protect\ref{leme2}}{Proof of Lemma 12}}

Chapuy's bijection may naturally be transposed in the continuous
setting. Let $\ii\in\lbracket1,2^g g! \rbracket$ be an integer
corresponding to an
opening sequence, and $\mathbf T_\infty^{\ii}$ the real tree coded
by $\CCC_\infty^\ii$. The interval $[0,1]$ may be split into $2g+1$
intervals coding the two halves of $\mathfrak{f}^{\mathfrak
{e}_*}_\infty$ and the other
forests of $\TT_\infty$. Through the continuous analog of Chapuy's
bijection, these intervals are reordered into an order corresponding to
the opening sequence. We call $\varphi^\ii\dvtx[0,1] \to[0,1]$ the
bijection accounting for this reordering. It is a cadlag function with
derivative $1$ satisfying $\Lab_\infty(s)=\LLL^\ii_\infty(\varphi
^\ii
(s))$ for all $s\in[0,1]$.

In order to see that Lemma~\ref{leme2} is a consequence of
\cite{legall07tss}, Lemma 2.4, let us first see what happens to subtrees of
$\TT_\infty$ through the continuous analog of Chapuy's bijection. It is
natural to call root of $\TT_\infty$ the point $\partial\de\TT
_\infty
(u_\infty)$, where the real number $u_\infty$ was defined in
Proposition~\ref{cvint} as the limit of the integer coding the root
in $\tr_n$, properly rescaled. Using classical properties of the
Brownian motion together with Proposition~\ref{cvint}, it is easy to
see that, almost surely, $\partial$ is a leaf of $\TT_\infty$, so
that $\tau_{\partial}$ is well defined. Any subtree of~$\TT_\infty$ not
included in $\tau_{\partial}$ (these subtrees require extra care, we
will treat them separately) is transformed through Chapuy's bijection
into some subtree of the opened tree $\mathbf T_\infty^\ii$ (i.e., into
some tree to the left or right of some branch of $\mathbf T_\infty^\ii
$). This is easy to see when the subtree is not rooted at a node of
$\TT
_\infty$, and we saw at the end of Section~\ref{secdefgtree} that,
almost surely, all the subtrees are rooted outside the set of nodes of
$\TT_\infty$.

We reason by contradiction to rule out these subtrees. We call $\Leb$
the Lebesgue measure on $[0,1]$. Let us suppose that there exist $\eta
>0$, and some subtree $\tau$, coded by $[l,r]$, not included in $\tau
_{\partial}$, such that $\inf_{[l,r]} \Lab_\infty<\Lab_\infty
(l)-\eta
$, and
%
%
\begin{eqnarray}\label{liminf}
&&\liminf_{\eps\to0} \eps^{-2} \Leb\biggl( \biggl\{ s \in[l,r] \dvtx
\Lab_\infty(s) < \Lab_\infty(l) - \eta+ \eps;
\nonumber\\
&&\hspace*{72pt} \forall x \in
[\CC
_\infty(l), \CC_\infty(s)],\\
&&\hspace*{72pt}\Lab_\infty\bigl(\sup\{t\le s \dvtx \CC_\infty(t)=x\}\bigr) > \Lab_\infty
(l)-\eta
+\frac\eps8 \biggr\} \biggr) =0.
\nonumber
\end{eqnarray}

Note that, by definition of $\CCC_\infty^\ii$, the function $s
\mapsto
\CC_\infty(s) - \CCC^\ii_\infty(\varphi^\ii(s))$ is constant on
$[l,r]$. Let us call $l'\de\varphi^\ii(l)$ and $r'\de\varphi^\ii(r)$.
It is easy to see that~(\ref{liminf}) remains true when replacing,
respectively, $l$, $r$, $\CC_\infty$ and $\Lab_\infty$ with~$l'$, $r'$,
$\CCC^\ii_\infty$ and~$\LLL^\ii_\infty$. Thanks to
Proposition \ref
{cs2}, the conclusion of~\cite{legall07tss}, Lemma~2.4, is also true for
the opened tree $\mathbf T_\infty^\ii$, and the fact that $[l',r']$
codes a~subtree of the opened tree yields a~contradiction.

We then use a re-rooting argument to conclude. With positive
probability, $\tau_{\partial}$~is no longer the tree containing the
root in the uniformly re-rooted $g$-tree. Let us suppose that, with
positive probability, there exists a subtree of $\TT_\infty$ included
in $\tau_{\partial}$, satisfying the hypotheses but not the conclusion
of Lemma~\ref{leme2}. Then, with positive probability, there will exist
a subtree not included in the tree containing the root of the uniformly
re-rooted $g$-tree, satisfying the hypotheses but not the conclusion of
Lemma~\ref{leme2}. The fact that the uniformly re-rooted $g$-tree has
the same law as $\TT_\infty$ yields a contradiction.

\subsubsection{\texorpdfstring{Proof of Lemma \protect\ref{leme4}}{Proof of Lemma 13}}

Using the same arguments as in~\cite{legall08glp}, we can see that
Lemma~\ref{leme4} is a consequence of the following lemma (see
\cite{legall08glp}, Corollary 6.2):
\begin{lem}
For every $p\ge1$ and every $\delta\in(0,1]$, there exists a
constant $c_{p,\delta} < \infty$ such that, for every $\eps> 0$,
\[
\mathbb{E}\biggl[ \biggl(\int_0^1 \mathbh{1}_{\{\Lab_\infty(s) \le\min\Lab
_\infty+ \eps \}} \,ds\biggr)^p \biggr] \le c_{p,\delta} \eps^{4p-\delta}.
\]
\end{lem}
\begin{pf}
This readily comes from~\cite{legall08glp}, Lemma 6.1, stating that for
every $p\ge1$ and every $\delta\in(0,1]$, there exists a constant
$c'_{p,\delta} < \infty$ such that, for every $\eps> 0$,
\[
\mathbb{E}\biggl[ \biggl(\int_0^1 \mathbh{1}_{\{Z_s \le\min Z + \eps\}} \,ds\biggr)^p
\biggr] \le
c'_{p,\delta} \eps^{4p-\delta}.
\]
Obviously, this still holds for $\delta\in(1,2]$. Using the link
between $\Lab_\infty$ and $\LLL_\infty$, as well as
Proposition \ref
{cs2}, we see that, for $p\ge1$ and $\delta\in(0,1]$,
\begin{eqnarray*}
\mathbb{E}\biggl[ \biggl(\int_0^1 \mathbh{1}_{\{\Lab_\infty(s) \le\min\Lab
_\infty+ \eps \}} \,ds\biggr)^p \biggr] &=& \mathbb{E}\biggl[ \biggl(\int_0^1 \mathbh{1}_{\{
\LLL_\infty(s) \le\min\LLL _\infty+ \eps\}}\, ds\biggr)^p \biggr]\\
&=& \mathbb{E}\biggl[ W \biggl(\int_0^1 \mathbh{1}_{\{Z_s \le\min Z +
\eps\}}\,
ds\biggr)^p \biggr]\\
&\le&(\mathbb{E}[ W^2 ] c'_{2p,2\delta})^{1/2} \eps^{4p-\delta
}\\
&=&c_{p,\delta} \eps^{4p-\delta},
\end{eqnarray*}
where $c_{p,\delta} \de(\mathbb{E}[ W^2 ] c'_{2p,2\delta})^{
1/2} <
\infty$, by~\cite{chapuy08sum}, Lemma 10.
\end{pf}

\section*{Acknowledgment}
The author is sincerely grateful to Gr\'{e}gory Miermont for the
precious advice and support he provided during the accomplishment of
this work.



%
\printaddresses

\end{document}